\documentclass[final]{siamltex}

\usepackage{graphicx}
\usepackage{url}
\usepackage{epstopdf}

\usepackage{amsmath,amssymb}
\usepackage{subfig}
\usepackage{color}
\usepackage{algorithm}
\usepackage{bm}

\graphicspath{{./}{./figs/}}

  \def\clap#1{\hbox to 0pt{\hss#1\hss}}

\providecommand{\mat}[1]{\bm{#1}}%
\renewcommand{\vec}[1]{\mathbf{#1}}

\providecommand{\mA}{\ensuremath{\mat{A}}}

\providecommand{\mC}{\ensuremath{\mat{C}}}

\providecommand{\mI}{\ensuremath{\mat{I}}}

\providecommand{\mM}{\ensuremath{\mat{M}}}

\providecommand{\mQ}{\ensuremath{\mat{Q}}}

\providecommand{\mW}{\ensuremath{\mat{W}}}

\providecommand{\vd}{\ensuremath{\vec{d}}}
\providecommand{\ve}{\ensuremath{\vec{e}}}

\providecommand{\vm}{\ensuremath{\vec{m}}}
\providecommand{\vn}{\ensuremath{\vec{n}}}

\providecommand{\vs}{\ensuremath{\vec{s}}}

\providecommand{\vw}{\ensuremath{\vec{w}}}
\providecommand{\vx}{\ensuremath{\vec{x}}}
\providecommand{\vy}{\ensuremath{\vec{y}}}
\providecommand{\vz}{\ensuremath{\vec{z}}}

\newcommand{\hmC}{\hat{\mC}}

\newcommand{\hLambda}{\hat{\Lambda}}

\newcommand{\hmW}{\hat{\mW}}
\newcommand{\hlambda}{\hat{\lambda}}

\newcommand{\hvy}{\hat{\vy}}
\newcommand{\hvz}{\hat{\vz}}
\newcommand{\hpi}{\hat{\pi}}

\newcommand{\sC}{\mathcal{C}}

\newcommand{\bmat}[1]{\begin{bmatrix}#1\end{bmatrix}}

\newcommand{\prior}{\rho_{\operatorname{pri}}}
\newcommand{\priory}{\rho_{\operatorname{pri}}}
\newcommand{\priorz}{\rho_{\operatorname{pri}}}
\newcommand{\lik}{\rho_{\operatorname{lik}}}

\newcommand{\post}{\rho_{\operatorname{pos}}}
\newcommand{\cpost}{c_{\operatorname{pos}}}
\newcommand{\cpi}{c_\pi}
\newcommand{\hg}{\hat{g}}
\renewcommand{\hpi}{\hat{\pi}}
\newcommand{\chpi}{c_{\hat{\pi}}}
\newcommand{\geps}{g_\varepsilon}
\newcommand{\pieps}{\pi_\varepsilon}

\newcommand{\hgeps}{\hat{g}_\varepsilon}
\newcommand{\hpieps}{\hat{\pi}_\varepsilon}
\newcommand{\chpieps}{c_{\hat{\pi}_\varepsilon}}

\newcommand{\xtrue}{\vx_{\text{true}}}
\newcommand{\Mtil}{\tilde{\mM}}

\renewenvironment{algorithm}[1][\relax]{\refstepcounter{algorithm}%
\addcontentsline{loa}{algorithm}%
    {\protect\numberline{Algorithm~\thealgorithm}{\ignorespaces#1}}%
\par\vspace{1\baselineskip}%
\expandafter\ifx#1\relax
\parindent0pt {\scshape\bfseries Algorithm~\thealgorithm.}\\%%%
\else
\parindent0pt {\scshape\bfseries
Algorithm~\thealgorithm.}\enspace{\bfseries#1.}\\%%%
\fi}
{\vspace{1\baselineskip}\par}
\title{ACCELERATING MCMC WITH ACTIVE SUBSPACES}

\author{
Paul G.~Constantine\thanks{Ben L.~Fryrear Assistant Professor, Department of Applied Mathematics and Statistics, Colorado School of Mines, Golden, CO 80401 (\texttt{paul.constantine@mines.edu}).}
\and
Carson Kent\thanks{Ph.D. student, Institute for Computational and Mathematical Engineering, Stanford University, Stanford, CA 94305 (\texttt{crkent@stanford.edu}).}
\and
Tan Bui-Thanh\thanks{Assistant Professor, Department of Aerospace Engineering and Engineering Mechanics, Institute for Computational Engineering and Sciences, University of Texas at Austin, Austin, TX 78705 (\texttt{tanbui@ices.utexas.edu}).}
}

\begin{document}
\maketitle

\begin{abstract}
The Markov chain Monte Carlo (MCMC) method is the computational workhorse for Bayesian inverse problems. However, MCMC struggles in high-dimensional parameter spaces, since its iterates must sequentially explore the high-dimensional space. This struggle is compounded in physical applications when the nonlinear forward model is computationally expensive. One approach to accelerate MCMC is to reduce the dimension of the state space. \emph{Active subspaces} are part of an emerging set of tools for subspace-based dimension reduction. An active subspace in a given inverse problem indicates a separation between a low-dimensional subspace that is informed by the data and its orthogonal complement that is constrained by the prior. With this information, one can run the sequential MCMC on the active variables while sampling independently according to the prior on the inactive variables. However, this approach to increase efficiency may introduce bias. We provide a bound on the Hellinger distance between the true posterior and its active subspace-exploiting approximation. And we demonstrate the active subspace-accelerated MCMC on two computational examples: (i) a two-dimensional parameter space with a quadratic forward model and one-dimensional active subspace and (ii) a 100-dimensional parameter space with a PDE-based forward model and a two-dimensional active subspace.
\end{abstract}



\pagestyle{myheadings}
\thispagestyle{plain}
\markboth{PAUL G.~CONSTANTINE, CARSON KENT, AND TAN BUI-THANH}{ACCELERATING MCMC WITH ACTIVE SUBSPACES}

\section{Introduction}
\label{sec:intro}

Modern computing enables powerful Bayesian inference methods to quantify uncertainty in complex science and engineering models. The workhorse of these methods is Markov chain Monte Carlo (MCMC), which constructs a Markov chain whose steps produce correlated samples from the conditional posterior density of the parameters given data. Several excellent references introduce and develop the theory and practice of MCMC in the context of statistical inverse problems~\cite{brooks2011handbook, kaipio2005statistical,Stuart2010}. However, standard MCMC remains an inefficient tool when (i) the dimension of the parameter space is large and/or (ii) the forward model in the likelihood is expensive. Recent papers have proposed strategies to increase MCMC's efficiency by introducing structure-exploiting proposal densities~\cite{Apte2007,DRAM2006,vrugt2009accelerating,Girolami11,Martin2012,Bui2014}. When the forward model numerically approximates a PDE solution, the differential operators may enable MCMC variants whose state space dimension is independent of the PDE discretization~\cite{Bui2012,Bui2014,Martin2012,Cui2014}. 

MCMC struggles in high dimensions because the Markov chain must explore the high-dimensional space to find regions of relatively high posterior probability. The recent paper by Cui, et al.~\cite{Cui2014} seeks to reduce the Markov chain's state space dimension by identifying a low-dimensional \emph{likelihood-informed subspace}. When an inverse problem is ill-posed, the data does not inform all parameters; Cui, et al.~use the eigendecomposition of the negative log-likelihood's Hessian---averaged over the parameter space with respect to the posterior---to identify a data-informed subspace; its orthogonal complement is informed by the prior. This separation leads to an efficient MCMC method when the forward model is twice differentiable, the prior is Gaussian, and the noise is Gaussian. They have extended these ideas to develop optimal posterior covariance approximations for linear forward models~\cite{spantini2015optimal} and dimension-independent methods for infinite dimensional problems~\cite{cui2014dimension}.

Our recent work in \emph{active subspaces} resembles the subspace-based dimension reduction of Cui, et al.~\cite{Cui2014}, though the context is broader than statistical inverse problems. The active subspace is defined for a scalar-valued, multivariate function; the active subspace basis consists of the principal eigenvectors of the average outer product of the gradient with itself. When a function admits an active subspace, one can exploit the low-dimensional structure to enable approximation, integration, and optimization in high dimensions~\cite{Constantine2014,constantine2015computing,asm2015}. We have applied this technique to help study aerospace models~\cite{constantine2014exploiting,Lukaczyk2014}, hydrological models~\cite{Jefferson2015}, and solar cell models~\cite{constantine2015discovering}. We review active subspaces and the resulting approximations in section \ref{sec:as}.

In this paper, we seek to discover the active subspace from the negative log-likelihood in a Bayesian inverse problem. If the negative log-likelihood admits an active subspace, then we use the active subspace to construct a function that approximates the Bayesian posterior. The approximate posterior enables an efficient MCMC that exploits the low-dimensional structure; the Markov chain operates on the active variables (i.e., the coordinates of the active subspace), while the inactive variables are drawn independently from their derived prior. In contrast to the likelihood-informed subspace of Cui, et al.~\cite{Cui2014}, the active subspace is defined by the gradient of the negative log-likelihood instead of a Gauss-Newton-based approximation of its Hessian. Additionally, the matrix that defines the active subspace is an average with respect to the prior as opposed to the posterior. This choice is computationally appealing, since one can draw independent samples from the prior and estimate integrals with standard Monte Carlo. Consequently, computing the active subspace is a preprocessing step that occurs before the MCMC.

In section \ref{sec:postapprox}, we analyze the approximation error in the posterior, measured by the Hellinger distance. We detail the MCMC algorithm with the active subspace in section \ref{sec:mcmc}. We demonstrate the approach in two synthetic Bayesian inverse problems: (i) a quadratic forward model with a two-dimensional parameter space and (ii)  a forward model governed by Poisson's equation, where the data are the PDE solution---perturbed by noise---at points on the boundary, and the parameters are the coefficients of the differential operator. We describe the models and show results in section \ref{sec:experiment}. In section \ref{sec:conclusion}, we summarize the method and discuss some practical limitations and future work. 

%Bayesian DOE~\cite{Long2013,Long2015}.

%\begin{itemize}
%\item Section \ref{sec:intro}, intro
%\item Section \ref{sec:as}, active subspaces, Section \ref{sec:approx}, approximating functions, Section \ref{sec:compas}, computing active subspaces, 
%\item Section \ref{sec:postapprox}, approximating the Bayesian posterior, Section \ref{sec:prior}, integrating with prior, Section \ref{sec:posterior}, integrating with posterior, Section \ref{sec:linear}, 
%\item Section \ref{sec:mcmc}, MCMC with active subspace
%\item Section \ref{sec:experiment}, experiment
%\item Section \ref{sec:conclusion}, conclusions
%\item Appendix \ref{sec:proofa}, proof of first theorem, Appendix \ref{sec:proofb}, proof of second theorem 
%\end{itemize}

\section{Active subspaces}
\label{sec:as}

The following description of active subspaces can be found in~\cite{Constantine2014,asm2015}; we include this section to keep the paper reasonably self contained. Let $f=f(\vx)$ be a function from $\mathbb{R}^m$ to $\mathbb{R}$; the input vector $\vx\in\mathbb{R}^m$ has $m$ independent components. Let $\rho:\mathbb{R}\rightarrow\mathbb{R}_+$ be a given probability density function. Assume that $\rho=\rho(\vx)$ and $\vx$ are such that
\begin{equation}
\int \vx\,\rho\,d\vx \;=\; 0, \qquad
\int \vx\,\vx^T\,\rho\,d\vx \;=\; \mI,
\end{equation}
where $\mI$ is the $m\times m$ identity matrix\footnote{Any random vector with a full rank covariance can be shifted and linearly transformed to satisfy the zero-mean and identity covariance assumptions.}. Also assume that $f$ is differentiable with gradient vector $\nabla f(\vx)\in\mathbb{R}^m$, whose components are square-integrable with respect to $\rho$. Define the $m\times m$ symmetric, positive semidefinite matrix $\mC$ and its eigenvalue decomposition as
\begin{equation}
\label{eq:C}
\mC \;=\; \int \nabla f(\vx)\,\nabla f(\vx)^T\,\rho(\vx)\,d\vx 
\;=\; \mW\Lambda\mW^T,
\end{equation}
where $\mW$ is the orthogonal matrix of eigenvectors, and $\Lambda$ is the diagonal matrix of eigenvalues ordered in descending order. The $i$th eigenvalue $\lambda_i$ satisfies
\begin{equation}
\label{eq:lambda}
\lambda_i \;=\; 
\vw_i^T\,\mC\,\vw_i \;=\;
\int \big( \vw_i^T\nabla f(\vx) \big)^2\,\rho(\vx)\,d\vx.
\end{equation}
In words, \eqref{eq:lambda} means that the $i$th eigenvalue measures the average, squared directional derivative of $f$ along the corresponding eigenvector $\vw_i$. Thus, $\lambda_i=0$ if and only if $f$ is constant along the direction $\vw_i$ in $\mathbb{R}^m$. 

To define the active subspace, assume $\lambda_n>\lambda_{n+1}$ for some $n<m$ and partition the eigenvalues and eigenvectors as
\begin{equation}
\label{eq:part}
\Lambda \;=\; \bmat{\Lambda_1 & \\ & \Lambda_2},\qquad
\mW \;=\; \bmat{\mW_1 & \mW_2},
\end{equation}
where $\Lambda_1$ contains the first $n$ eigenvalues, and the columns of $\mW_1$ are the first $n$ eigenvectors. The \emph{active subspace} is the span of the columns of $\mW_1$. However, the active subspace is not necessarily a subset of the domain of $f$---even when the domain of $f$ is $\mathbb{R}^m$. Instead, the columns of $\mW_1$ are a set of directions; perturbing $\vx$ along these directions changes $f(\vx)$ more, on average, than perturbing $\vx$ along the directions corresponding to the columns of $\mW_2$. Any $\vx\in\mathbb{R}^m$ can be written
\begin{equation}
\vx \;=\; \mW_1\mW_1^T\vx+\mW_2\mW_2^T\vx
\;=\; \mW_1\vy + \mW_2\vz,
\end{equation}
where $\vy=\mW_1^T\vx$ are the \emph{active variables}, and $\vz=\mW_2^T\vx$ are the \emph{inactive variables}. The density function $\rho$ begets a joint density between the active and inactive variables,
\begin{equation}
\label{eq:joint}
\rho(\vx) \;=\; \rho(\mW_1\vy + \mW_2\vz) \;=\; \rho(\vy,\vz),
\end{equation}
which leads to marginal and conditional densities under the standard construction. If $\rho$ is a standard Gaussian density on $\vx$ then, due to the orthogonal columns of $\mW_1$ and $\mW_2$, the marginal and conditional densities on $\vy$ and $\vz$ are standard Gaussian densities on $\mathbb{R}^n$ and $\mathbb{R}^{m-n}$, respectively.  

\subsection{Approximation with the active subspace}
\label{sec:approx}

If $\lambda_1,\dots,\lambda_n$ are much larger than $\lambda_{n+1},\dots,\lambda_m$, then we can approximate $f$ by a function of $n<m$ linear combinations of $\vx$. To construct this approximation, define $g:\mathbb{R}^n\rightarrow\mathbb{R}$ by the conditional average of $f$ given $\vy$,
\begin{equation}
\label{eq:cexp}
g(\vy) \;=\; \int f(\mW_1\vy+\mW_2\vz)\,\rho(\vz|\vy)\,d\vz,
\end{equation}
where $\rho(\vz|\vy)$ is the conditional density of $\vz$ given $\vy$. Under this construction, we have the following bound on the root-mean-squared error of the approximation $f(\vx)\approx g(\mW_1^T\vx)$.
\begin{theorem}[Theorem 3.1~\cite{Constantine2014}]
\label{thm:cexperr}
\begin{equation}
\label{eq:cexperr}
\left(\int (f(\vx) - g(\mW_1^T\vx))^2\,\rho\,d\vx\right)^{1/2}
\;\leq\;
C\,(\lambda_{n+1}+\cdots+\lambda_m)^{1/2},
\end{equation}
where $C$ is the Poincar\'{e} constant associated with the density $\rho$. 
\end{theorem}

If $\rho$ is a standard Gaussian density, then $C$ in \eqref{eq:cexperr} is 1~\cite{chen1982inequality}. The conditional expectation in \eqref{eq:cexp} is not useful for computation, since evaluating $g(\vy)$ involves computing an $(n-m)$-dimensional integral. To move toward a useful computational tool, we introduce the Monte Carlo approximation, $\hg\approx g$, defined as
\begin{equation}
\label{eq:mccexp}
\hg(\vy) \;=\; \frac{1}{M}\sum_{i=1}^M f(\mW_1\vy + \mW_2\vz_i),
\end{equation}
where $\vz_i$ are drawn independently according to the conditional density $\rho(\vz|\vy)$. The approximation $f(\vx)\approx \hg(\mW_1^T\vx)$ admits the following root-mean-squared error estimate.
\begin{theorem}[Theorem 3.2~\cite{Constantine2014}]
\label{thm:mccexperr}
\begin{equation}
\label{eq:mccexperr}
\left(\int (f(\vx) - \hat{g}(\mW_1^T\vx))^2\,\rho\,d\vx\right)^{1/2}
\;\leq\;
C(1+M^{-1/2})(\lambda_{n+1}+\cdots+\lambda_m)^{1/2},
\end{equation}
where $C$ is the constant from Theorem \ref{thm:cexperr}.
\end{theorem}

If $f$ is such that the eigenvalues $\lambda_{n+1}=\cdots=\lambda_m=0$, then the Monte Carlo estimate is exact for any number $M>0$ of samples. Another way to see this is that $\lambda_{n+1}=\cdots=\lambda_m=0$ implies $f$ is constant along the directions corresponding to $\mW_2$'s columns, and the average of a constant is the constant.

\subsection{Computing the active subspace with Monte Carlo}
\label{sec:compas}
We assume that the dimension $m$ is sufficiently large that Monte Carlo is the most practical choice to estimate the matrix $\mC$ in \eqref{eq:C}. (The two-parameter example in section \ref{sec:2parm} uses a tensor product Gauss-Hermite quadrature rule.) Our recent work analyzes the Monte Carlo approximation~\cite{constantine2015computing}. Let $\vx_j$ with $j=1,\dots,N$ be drawn independently according to the density $\rho$. For each $\vx_j$, compute the gradient $\nabla f_j=\nabla f(\vx_j)$. Then approximate
\begin{equation}
\label{eq:mc}
\mC \;\approx\; \hmC 
\;=\; \frac{1}{N}\sum_{j=1}^N \nabla f_j\,\nabla f_j^T 
\;=\; \hmW\hLambda\hmW^T.
\end{equation}
Let the estimated eigenvalues $\hLambda$ and eigenvectors $\hmW$ be partitioned as in \eqref{eq:part}. Let $\varepsilon$ be the error in the estimated active subspace,
\begin{equation}
\label{eq:suberr}
\varepsilon \;=\; \|\mW_1\mW_1^T - \hmW_1\hmW_1^T\| \;=\; \|\hmW_1^T\mW_2\|,
\end{equation}
where $\|\cdot\|$ is the matrix 2-norm; see~\cite[Section 2.5.3]{golub2013}. In~\cite{constantine2015computing}, we show that when the number $N$ of samples in \eqref{eq:mc} is greater than a problem dependent lower bound, the relative error in the estimated eigenvalues $\hLambda$ falls below a user-specified tolerance with high probability; the lower bound grows as the log of the dimension $m$. With sufficient samples, the subspace error $\varepsilon$ in \eqref{eq:suberr} satisfies~\cite[Corollary 3.7]{constantine2015computing}
\begin{equation}
\label{eq:suberrbnd}
\varepsilon \;\leq\;
\frac{4\,\lambda_1\,\delta}{\lambda_n - \lambda_{n+1}},
\end{equation}
where $\delta$ is the user-specified error tolerance for the eigenvalue estimates. In practice, if we observe a gap between $\hlambda_n$ and $\hlambda_{n+1}$, then the Monte Carlo procedure gives a good estimate of the $n$-dimensional subspace. 

\subsection{Approximation with the estimated active subspace}
Practical versions of the approximations in \eqref{eq:cexp} and \eqref{eq:mccexp} use the estimated eigenvectors $\hmW_1$. The estimated active and inactive variables are $\hvy=\hmW_1^T\vx$ and $\hvz=\hmW_2^T\vx$, respectively. And the joint density is similar to \eqref{eq:joint}, since $\hmW$ is an orthogonal matrix. The conditional average using the estimated eigenvectors is
\begin{equation}
\label{eq:epscexp}
\geps(\hvy) \;=\; \int f(\hmW_1\hvy+\hmW_2\hvz)\,\rho(\hvz|\hvy)\,d\hvz.
\end{equation}
The root-mean-squared error in the approximation $f(\vx)\approx\geps(\hmW_1^T\vx)$ is given in the next theorem.
\begin{theorem}[Theorem 3.5~\cite{Constantine2014}]
\label{thm:epscexperr}
\begin{equation}
\label{eq:epscexperr}
\begin{aligned}
&\left(\int (f(\vx) - \geps(\hmW_1^T\vx))^2 \,\rho\,d\vx\right)^{1/2}\\
&\qquad\leq\;
C\left(
\varepsilon(\lambda_1+\cdots+\lambda_n)^{1/2} + 
(\lambda_{n+1}+\cdots+\lambda_m)^{1/2}
\right),
\end{aligned}
\end{equation}
where $C$ is from Theorem \ref{thm:cexperr}, and $\varepsilon$ is the subspace error from \eqref{eq:suberr}.
\end{theorem}

When the error $\varepsilon$ in the estimated subspace (see \eqref{eq:suberr}) is not zero, the error estimate includes contributions from the larger eigenvalues. The Monte Carlo estimate $\hgeps$ of the conditional average $\geps$ is
\begin{equation}
\label{eq:epsmccexp}
\hgeps(\hvy) \;=\; 
\frac{1}{M}\sum_{i=1}^M f(\hmW_1\hvy + \hmW_2\hvz_i), 
\end{equation}
where $\hvz_i$ are drawn independently from the conditional density $\rho(\hvz|\hvy)$. The next theorem bounds the root-mean-squared error in the Monte Carlo approximation. 
\begin{theorem}[Theorem 3.6~\cite{Constantine2014}]
\label{thm:epsmccexperr}
\begin{equation}
\label{eq:epsmccexperr}
\begin{aligned}
&\left(\int (f(\vx) - \hgeps(\hmW_1^T\vx))^2 \,\rho\,d\vx\right)^{1/2}\\
&\qquad\leq\;
C(1+M^{-1/2})\left(
\varepsilon(\lambda_1+\cdots+\lambda_n)^{1/2} + 
(\lambda_{n+1}+\cdots+\lambda_m)^{1/2}
\right),
\end{aligned}
\end{equation}
where $C$ is from Theorem \ref{thm:cexperr}, and $\varepsilon$ is the subspace error from \eqref{eq:suberr}. 
\end{theorem}

To summarize, we introduce four low-dimensional approximations for $f(\vx)$ that exploit the active subspace: (i) the conditional average $g$ \eqref{eq:cexp}, (ii) a Monte Carlo approximation $\hg$ \eqref{eq:mccexp} of the conditional average, (iii) the conditional average $\geps$ \eqref{eq:epscexp} constructed with estimated eigenvectors $\hmW_1$, and (iv) the Monte Carlo approximation $\hgeps$ \eqref{eq:epsmccexp} of the conditional average $\geps$. For each approximation, we bound the root-mean-squared error by $\mC$'s eigenvalues and the subspace error. Practical computations use $\hgeps$. 

\section{Approximating the Bayesian posterior}
\label{sec:postapprox}

We consider the following Bayesian inverse problem. Assume an additive noise model,
\begin{equation}
\label{eq:addnoise}
\vd \;=\; \vm(\vx) + \ve,
\end{equation}
where $\vd\in\mathbb{R}^d$ are the random data, $\vx\in\mathbb{R}^m$ are the random parameters, $\vm:\mathbb{R}^m\rightarrow\mathbb{R}^d$ is the deterministic map from parameters to observables (e.g., the observation operator composed with the forward model in PDE-based inverse problems), and $\ve\in\mathbb{R}^d$ is the random noise. We assume that $\vm$ is a differentiable function of $\vx$. For simplicity, we assume $\ve$ is a zero-mean Gaussian random vector with covariance $\sigma^2\mI$; correlated noise can be transformed to uncorrelated noise with standard linear transformations using the square root of the noise covariance matrix.
%If the given problem includes correlated noise, then we assume it has been properly decorrelated. 
%We assume a standard Gaussian prior density on the parameters,
%\begin{equation}
%\label{eq:prior}
%\prior(\vx) \;=\;
%\left(2\pi\right)^{\frac{-m}{2}}\, \exp\left(\frac{-\vx^T\vx}{2}\right).
%\end{equation}
%Again, a more general Gaussian prior can be transformed to a standard Gaussian by shifting according to the mean and linearly transforming with the square root of the prior covariance matrix. 
The Gaussian noise assumption leads to the likelihood,
\begin{equation}
\label{eq:likelihood}
\lik(\vd,\,\vx) \;=\; \exp\left(
\frac{-\|\vd-\vm(\vx)\|^2}{2\sigma^2}
\right),
\end{equation}
where $\|\cdot\|$ is the Euclidean norm. Denote the prior density on the parameters by $\prior(\vx)$, and let $\post(\vx)$ be the conditional density of $\vx$ given $\vd$. Bayes Theorem implies
\begin{equation}
\label{eq:bayes}
\post(\vx) \;=\; \cpost^{-1}\,\lik(\vd,\,\vx)\,\prior(\vx),
\qquad
\cpost \;=\; \int \lik(\vd,\,\vx)\,\prior(\vx)\,d\vx.
\end{equation}

\subsection{Identifying the active subspace}
\label{eq:identify}

To apply the active subspace machinery from section \ref{sec:as}, we must identify the scalar-valued function $f(\vx)$ and the density function $\rho(\vx)$ in \eqref{eq:C}. Similar to the likelihood-informed subspace of Cui, et al.~\cite{Cui2014}, we choose $f(\vx)$ to be the negative log-likelihood,
\begin{equation}
\label{eq:misfit}
f(\vx) \;=\; \frac{1}{2\sigma^2}\,\|\vd-\vm(\vx)\|^2.
\end{equation}
This function is often called the \emph{data misfit} function, or just \emph{misfit}, and it is closely related to an optimizer's objective function in deterministic inverse problems. For a given point in the parameter space, $f(\vx)$ measures how far the modeled observations are from the given data. If we use the misfit to define the active subspace, then the orthogonal complement---i.e., the inactive subspace---identifies directions along which the likelihood is relatively flat. Perturbing the parameters along the inactive subspace changes the likelihood relatively little, on average. We wish to exploit this structure, when present, to accelerate the MCMC. The gradient of the misfit is
\begin{equation}
\label{eq:misfitgrad}
\nabla f(\vx) \;=\; \frac{1}{\sigma^2}\,\nabla \vm(\vx)^T\,(\vd-\vm(\vx)),
\end{equation}
where $\nabla \vm \in\mathbb{R}^{d\times m}$ is the Jacobian of the parameter-to-observable map. The misfit and its gradient depend on the data $\vd$; this becomes important when we choose $\rho$. 

\subsection{Integrating against the prior}
\label{sec:prior}

Cui, et al.~\cite{Cui2014} average the misfit's prior-preconditioned Gauss-Newton Hessian with respect to the posterior density in \eqref{eq:bayes} to estimate the likelihood-informed subspace. Their subspace is then conditioned on the data. In contrast, we compute the averages defining $\mC$ in \eqref{eq:C} using $\rho=\prior$, which requires careful interpretation of the data $\vd$. In the model \eqref{eq:addnoise}, $\vd$ is a random variable whose mean depends on $\vx$; in other words, $\vd$ and $\vx$ are not independent. Therefore, we cannot integrate against the prior without $\vd$ changing as $\vx$ varies. However, if we treat the realization $\vd$ as a fixed and constant vector, then we can integrate $\nabla f\nabla f^T$, which depends on $\vd$, against a density function equal to $\prior$ without issue. The integrals defining $\mC$ from \eqref{eq:C} and the approximation $g$ from \eqref{eq:cexp} are well-defined. We sacrifice the probabilistic interpretation of the data $\vd$. Hence, we also sacrifice the interpretation of $\mC$ and all derived quantities as random variables conditioned on $\vd$; instead, they are functions of the fixed vector $\vd$. Additionally, if we are given a new set of data, then we must recompute the active subspace. But we gain the practical advantage of estimating $\mC$ with simple Monte Carlo as in section \ref{sec:compas}, since we can sample independently from the density $\prior$. With this choice of $\rho$, the Monte Carlo-based eigenvector estimates $\hmW$ from \eqref{eq:mc} are computable, and approximations to the misfit $f(\vx)$, namely, $\hg$ from \eqref{eq:mccexp}, $\geps$ from \eqref{eq:epscexp}, and $\hgeps$ from \eqref{eq:epsmccexp}, are also well-defined.

Consider the approximation $f(\vx)\approx g(\mW_1^T\vx)$ from \eqref{eq:cexp}. The prior can be factored as
\begin{equation}
\label{eq:sepprior}
\prior(\vx) \;=\; \prior(\vy,\vz) \;=\; \priory(\vy)\,\priorz(\vz|\vy),
\end{equation}
where $\priory(\vy)$ is the marginal density of $\vy$, and $\priorz(\vz|\vy)$ is the conditional density of $\vz$ given $\vy$. If $\prior(\vx)$ is a standard Gaussian on $\mathbb{R}^m$, then $\prior(\vy)$ is a standard Gaussian on $\mathbb{R}^n$, and $\prior(\vz|\vy)$ is a standard Gaussian on $\mathbb{R}^{m-n}$ that is independent of $\vy$. We construct an approximate posterior $\pi(\vx)$ as
\begin{equation}
\label{eq:approxpost}
\begin{aligned}
\post(\vx) &\approx \pi(\vx)\\
&= \cpi^{-1}\,\exp(-g(\mW_1^T\vx))\,\prior(\vx)\\
&= \cpi^{-1}\,\exp(-g(\vy))\,\priory(\vy)\,\priorz(\vz|\vy),
\end{aligned}
\end{equation}
where
\begin{equation}
\label{eq:cpi}
\cpi \;=\; \int \exp(-g(\mW_1^T\vx))\,\prior(\vx)\,d\vx.
\end{equation}
Since $\exp(\cdot)>0$, $\cpi>0$. Also, $g\geq 0$ implies $\exp(-g)\leq 1$, so $\cpi\leq 1$. The approximation in \eqref{eq:approxpost} suggests a strategy for MCMC that runs the Markov chain only on the active variables $\vy$ while sampling independently from the prior $\priorz(\vz|\vy)$ on the inactive variables; we explore this strategy in section \ref{sec:mcmc}. Before the computational exploration, we study the approximation properties of $\pi$ and similar constructions using $\hg$, $\geps$, and $\hgeps$ in place of $g$. We use the Hellinger distance~\cite{Gibbs2002} to quantify the approximation errors. 

\begin{theorem}
\label{thm:postapprox}
Let $\pi$ be defined as in \eqref{eq:approxpost}, and define the approximate posteriors $\hpi$, $\pieps$, and $\hpieps$ using $\hg$, $\geps$, and $\hgeps$, respectively, in place of $g$---all constructed with $\rho=\prior$. Define the constant $L$ as
\begin{equation}
\label{eq:L}
L^2 \;=\; \frac{1}{8}
\left[
\left(
\int \exp(-f)\,\prior\,d\vx
\right)
\left(
\exp\left(-\int f\,\prior\,d\vx \right)
\right)
\right]^{-1/2}.
\end{equation}
Then the Hellinger distances between the approximate posteriors and the true posterior $\post$ are bounded as follows.
\begin{align}
H(\post,\pi) &\leq\; 
L\,C\,
(\lambda_{n+1}+\cdots+\lambda_m)^{1/2},
\label{eq:hprior1}\\
H(\post,\hpi) &\leq\;
L\,C\,
\left(1+M^{-1/2}\right)(\lambda_{n+1}+\cdots+\lambda_m)^{1/2},
\label{eq:hprior2}\\
H(\post,\pieps) &\leq\; 
L\,C\,
\left(
\varepsilon\,(\lambda_1+\cdots+\lambda_n)^{1/2} +
(\lambda_{n+1}+\cdots+\lambda_m)^{1/2}
\right),
\label{eq:hprior3}\\
H(\post,\hpieps) &\leq\; 
L\,C\,
\left(1+M^{-1/2}\right)
\left(
\varepsilon\,(\lambda_1+\cdots+\lambda_n)^{1/2} +
(\lambda_{n+1}+\cdots+\lambda_m)^{1/2}
\right),
\label{eq:hprior4}
\end{align}
where $C$ is the Poincar\'{e} constant associated with $\prior$, $M$ is from \eqref{eq:mccexp} and \eqref{eq:epsmccexp}, and $\varepsilon$ is from \eqref{eq:suberr}. 
\end{theorem}

The bound \eqref{eq:hprior1} is an improved and extended version of Theorem 4.9 in~\cite{asm2015}. The proof of Theorem \ref{thm:postapprox} is in Appendix \ref{sec:proofb}. If the eigenvalues $\lambda_{n+1},\dots,\lambda_m$ are small, and if the error $\varepsilon$ in the numerically estimated active subspace is small, then the Hellinger distances between the posterior and its approximations are small. The Hellinger distance is a useful metric, because it provides an upper bound on the posterior mean and covariance; see~\cite[Lemma 6.37]{Stuart2010}. 

\subsection{Linear forward model}
\label{sec:linear}

Consider the case where the forward model is linear in the parameters,
\begin{equation}
\vm(\vx) \;=\; \mM\vx, \qquad \mM\in\mathbb{R}^{d\times m}. 
\end{equation}
Assume $\mM$ has rank $r$. If the prior and measurement noise are Gaussian, then the posterior is a Gaussian density whose mean and covariance have closed-form expressions in terms of the data, the prior covariance, and the noise covariance. A Gaussian density is completely characterized by its mean and covariance, so this case is not interesting to study with MCMC. Nevertheless, we can examine how the posterior approximation $\pi$ from \eqref{eq:approxpost} compares to the true posterior. The gradient of the misfit \eqref{eq:misfit} with the linear forward model is
\begin{equation}
\nabla f(\vx) \;=\; \frac{1}{\sigma^2}\,\mM^T(\mM\vx-\vd).
\end{equation}
Consider the case where $\prior$ is a standard Gaussian, i.e., 
\begin{equation}
\prior(\vx) \;=\;
\left(2\pi\right)^{\frac{-m}{2}}\, \exp\left(\frac{-\vx^T\vx}{2}\right).
\end{equation}
Using $\rho=\prior$, the matrix $\mC$ from \eqref{eq:C} is
\begin{equation}
\label{eq:Clin}
\mC \;=\; \frac{1}{\sigma^4}\,\mM^T(\mM\mM^T + \vd\vd^T)\mM. 
\end{equation}
If $\mW_1$ are the first $n<m$ eigenvectors of $\mC$ from \eqref{eq:Clin}, then the conditional average $g(\vy)$ from \eqref{eq:cexp} is
\begin{equation}
\begin{aligned}
g(\vy) &= \frac{1}{2\sigma^2} \int \left\|
\mM(\mW_1\vy + \mW_2\vz) - \vd
\right\|^2\,(2\pi)^{\frac{-(m-n)}{2}}\,\exp\left(
\frac{-\vz^T\vz}{2}
\right)\,d\vz\\
&= \frac{1}{2\sigma^2} \left(\left\|
\mM\mW_1\vy - \vd
\right\|^2
+ \gamma^2
\right),
\end{aligned}
\end{equation}
where 
\begin{equation}
\gamma^2 \;=\; \int \vz^T\mW_2^T\mM^T\mM\mW_2\vz
\,(2\pi)^{\frac{-(m-n)}{2}}\,\exp\left(
\frac{-\vz^T\vz}{2}
\right)\,d\vz
\end{equation}
is independent of $\vy$. The posterior approximation $\pi$ is
\begin{equation}
\label{eq:pilin}
\pi(\vx) \;=\; \cpi^{-1}\,\exp\left(\frac{-\gamma^2}{2\sigma^2}\right)
\,\exp\left(
\frac{-\|\mM\mW_1\mW_1^T\vx - \vd\|^2}{2\sigma^2}
\right)\,\prior(\vx).
\end{equation}
Note that $\exp(-\gamma^2/2\sigma^2)$ is independent of $\vx$. Using standard manipulations as in~\cite[Chapter 8]{Calvetti2007}, we can write down $\pi$'s mean $\mu$ and covariance matrix $\Gamma$,
\begin{equation}
\label{eq:meancovapprox}
\begin{aligned}
\mu &= 
\Mtil^T(\Mtil\Mtil^T + \sigma^2\mI)^{-1}\vd,\\
\Gamma &= 
\mI - \Mtil^T(\Mtil\Mtil^T + \sigma^2\mI)^{-1}\Mtil,
\end{aligned}
\end{equation}
where
\begin{equation}
\Mtil \;=\; \mM\mW_1\mW_1^T.
\end{equation}
Since $\rank(\mM)$ is $r$, $\rank(\mC)\leq r$ for $\mC$ from \eqref{eq:Clin}. Therefore, $\mC$'s eigenvalues $\lambda_{r+1},\dots,\lambda_m$ are zero. Applying Theorem \ref{thm:postapprox}, the Hellinger distance between $\pi$ and $\post$ is zero when the number $n$ of active variables is greater than or equal to $r$. In that case, the mean and covariance approximations in \eqref{eq:meancovapprox} are exact. 

Flath, et al.~\cite{Flath2011} and Bui-Thanh, et al.~\cite{Bui2012} construct an approximation to the posterior covariance with a rank-$k$ update of the prior covariance. They derive the rank-$k$ update from the $k$ dominant eigenpairs of the negative log-likelihood's Hessian, preconditioned by the prior. Spantini, et al.~\cite{spantini2015optimal} recently showed that such an update is optimal in a general class of matrix norms that implies an optimal posterior approximation in the Hellinger distance. The approximation $\Gamma$ in \eqref{eq:meancovapprox} is indeed a rank-$n$ update to the prior covariance $\mI$. However, it is difficult to compare this update to the Hessian-based constructions when $n<r$ because the eigenvectors $\mW_1$ depend on the fixed data vector $\vd$, while the eigenpairs of the negative log-likelihood's Hessian do not. It is possible to construct simple cases where the rank-$n$ update in \eqref{eq:meancovapprox} is a poor approximation for a particular fixed $\vd$. However, such a deficiency reveals little about the inverse problem with nonlinear forward models---where MCMC methods are most appropriate. For nonlinear forward models, it may be that data-dependent approximations (i.e., approximations that depend on $\vd$) are preferable to data-independent approximations for a particular realization of the data. 
%Such explorations are beyond the scope of the present paper. 

\section{MCMC with the active subspace}
\label{sec:mcmc}

Recall that choosing $\rho=\prior$ in \eqref{eq:C} and \eqref{eq:cexp} requires us to interpret the data $\vd$ as a fixed and constant vector; we sacrifice the probabilistic interpretation of $\vd$ and any quantities that depend on $\vd$. However, using $\rho=\prior$ allows us to use simple Monte Carlo to estimate the integrals, which can be done in parallel. Theorem \ref{thm:postapprox} contains error estimates for all likelihood approximations: the conditional expectation $g$, its Monte Carlo estimate $\hg$, and the analogous approximations using the estimated eigenvectors, $\geps$ and $\hgeps$. We restrict attention to the computable approximation $\hgeps(\hmW_1^T\vx)$ from \eqref{eq:epsmccexp} of the misfit $f(\vx)$ from \eqref{eq:misfit}. Recall that $\hgeps$ is a Monte Carlo estimate of the conditional expectation of $f$ given the estimated active variables $\hvy=\hmW_1^T\vx$. 

We first compute the eigenpair estimates $\hmW$ and $\hLambda$ with Monte Carlo as in section \ref{sec:compas} using the misfit's gradient \eqref{eq:misfitgrad}, where the samples $\vx_j$ are drawn independently according to the prior. One can view this step as preprocessing before running any MCMC. This preprocessing step checks for exploitable, low-dimensional structure indicated by (i) a gap in the estimated eigenvalues and (ii) small eigenvalues following the gap. A gap in the spectrum indicates that the Monte Carlo procedure in section \ref{sec:compas} can accurately estimate the active subspace; see \eqref{eq:suberrbnd}. In~\cite[Section 4]{constantine2015computing}, we describe a practical bootstrap procedure that can aid in the assessing the quality of the estimated active subspace. If the eigenvalues following the gap are small, then Theorem \ref{thm:postapprox} gives confidence that the approximate posterior is close to the true posterior. If these two conditions are not satisfied, then the problem may not be a good candidate for the proposed active subspace-accelerated MCMC, since there is no evidence that an exploitable active subspace exists based on the misfit function $f(\vx)$ and the prior density $\rho=\prior$. 

Assuming we have identified an active subspace, we propose an MCMC method that exploits the low-dimensional structure. Algorithm \ref{alg:mcmc} is an active subspace-exploiting variant of the Metropolis-Hastings method, sometimes called the \emph{random walk} method, outlined by Kaipio and Sommersalo~\cite[Chapter 3]{kaipio2005statistical}. The essential idea is to run the Markov chain only on the $n$ active variables instead of all $m$ variables. Therefore, we expect the chain to mix faster than MCMC on all $m$ variables. We do not compare to more sophisticated sampling schemes, because they can be adapted to exploit the active subspace in the same way. 

\begin{algorithm}{Markov chain Monte Carlo with the Active Subspace}{\label{alg:mcmc}}

\noindent Pick an initial value $\hvy_1$, and compute $\hgeps(\hvy_1)$. Set $k=1$. 
\begin{enumerate}
\item Draw $\hvy'\in\mathbb{R}^n$ from a symmetric proposal density centered at $\hvy_k$.
\item Compute $\hgeps(\hvy')$ as in \eqref{eq:epsmccexp}, where $f$ is the misfit function \eqref{eq:misfit}, and $\hvz_i$ are drawn independently according to $\rho(\hvz|\hvy)$, which is a standard Gaussian on $\mathbb{R}^{m-n}$.
\item Compute the acceptance ratio
\begin{equation}
\label{eq:mcmcaccept} 
\gamma(\hvy_k,\hvy') \;=\; \text{minimum}\,\left(1,\, 
\frac{
\exp(-\hgeps(\hvy'))\priory(\hvy')
}{
\exp(-\hgeps(\hvy_k))\priory(\hvy_k)
} \right).
\end{equation}
\item Draw $t$ uniformly from $[0,1]$.
\item If $\gamma(\hvy_k,\hvy')\geq t$, set $\hvy_{k+1}=\hvy'$. Otherwise, set $\hvy_{k+1}=\hvy_k$.
\item Increment $k$ and repeat. 
\end{enumerate}
\end{algorithm}

Step 2 computes $\hgeps(\hvy')$ from \eqref{eq:epsmccexp}. Each sample in \eqref{eq:epsmccexp} requires an independent evaluation of the misfit $f(\vx)$, and hence the parameter-to-observable map $\vm(\vx)$. Therefore, each step of the Markov chain in the $n$ active variables uses $M$ forward model evaluations---compared to one forward model evaluation for each step in the standard Metropolis-Hastings. We expect that in many problems (such as the example in section \ref{sec:experiment}), the dimension reduction enabled by the active subspace is far more valuable---in terms of forward model evaluations needed to sufficiently sample the space---than the penalty of a factor of $M$ increase. A factor $M$ increase is much smaller than the exponential growth of the parameter space with dimension.

How large should $M$ be to compute $\hgeps$ in Step 2? The error estimate \eqref{eq:hprior4} in Theorem \ref{thm:postapprox} provides some guidance. Note that $M$ enters the error bound through the multiplicative term $1+M^{-1/2}$. Thus, the effect of $M$ is bounded between 1 (for very large $M$) and 2 (for $M=1$). In other words, $M$ has relatively little effect on the error estimate; the eigenvalues matter much more. If the eigenvalues $\lambda_{n+1},\dots,\lambda_m$ are small, and if the gap $\lambda_n-\lambda_{n+1}$ is large---implying a small subspace error $\varepsilon$---then $M$ can be surprisingly small with little effect on the error. However, it is difficult to provide a universal numerical condition (e.g., a tolerance) for the eigenvalues that directly leads to a useful $M$. In the experiment in section \ref{sec:experiment}, we perform a preliminary computational experiment that justifies $M=10$; we recommend such an experiment in practice. 
%It's worth noting that an accurate estimate of $f(\vx)$ from $\hgeps(\hmW_1^T\vx)$ is sufficient but not necessary. We only need a good approximation of the acceptance ratio \eqref{eq:mcmcaccept}---and we only need that approximation for the majority of the steps, since many steps may be discarded during the MCMC. Analyzing the cost in this way is beyond the scope of this paper, but it is worth pursuing. For example, if we could show that $M=1$ was sufficient for Step 2 inside Algorithm 1, then the resulting algorithm would mix like a chain in $n<m$ dimensions, and each step would cost the same as a standard MCMC. 

Algorithm \ref{alg:mcmc} generates the set $\{\hvy_k\}$. These samples must be transformed to the space of the original parameters $\vx$ for inference. For each $\hvy_k$, draw independent realizations of $\hvz_{k,\ell}$ from the conditional density $\priorz(\hvz|\hvy_k)$ with $\ell = 1,\dots,P$. Then construct
\begin{equation}
\label{eq:reconstruct}
\vx_{k,\ell} \;=\; \hmW_1\hvy_k + \hmW_2\hvz_{k,\ell}.
\end{equation} 
To be sure, constructing $\{\vx_{k,\ell}\}$ from $\{\hvy_k\}$ requires no forward model evaluations; it only requires that one be able to draw independent samples from the conditional density $\priorz(\hvz|\hvy)$. If $\prior(\vx)$ is a standard Gaussian as in our numerical examples, then $\priorz(\hvz|\hvy)$ is a standard Gaussian, so drawing independent samples is straightforward. By a derivation similar to \eqref{eq:approxpost}, the set $\{\vx_{k,\ell}\}$ contains correlated samples from the approximate posterior $\hpieps$ defined as
\begin{equation}
\hpieps(\vx) \;=\; \chpieps^{-1}\,\exp(-\hgeps(\hmW_1^T\vx))\,\prior(\vx),
\end{equation}
where
\begin{equation}
\chpieps \;=\; \int \exp(-\hgeps(\hmW_1^T\vx))\,\prior(\vx)\,d\vx.
\end{equation}
Recall that Theorem \ref{thm:postapprox} bounds the Hellinger distance between $\hpieps$ and the true posterior in \eqref{eq:hprior4}. We expect the correlation in the set of samples $\{\vx_{k,\ell}\}$ to be much smaller than a set of samples drawn with MCMC directly on the parameters $\vx$, since the $\vx_{k,\ell}$'s from \eqref{eq:reconstruct} contain many independently sampled components, $\hvz_{k,\ell}$. The problem in section \ref{sec:pdebayes} shows an example of such behavior. 

\section{Numerical experiments}
\label{sec:experiment}

The following experiments used Matlab 2015b and Enthought Canopy Python 2.7.9 on a 2013 MacBook Air with 8GB of RAM for the computations with the quadratic model from section \ref{sec:2parm} and most postprocessing and plotting. The PDE-based experiment from section \ref{sec:pdebayes} ran on two processors from one node of Colorado School of Mines' Mio cluster (\url{inside.mines.edu/mio}) using the same Python distribution. The scripts and data to produce the figures for the following numerical experiments can be found at \url{bitbucket.org/paulcon/accelerating-mcmc-with-active-subspaces}. 

\subsection{Two-parameter model}
\label{sec:2parm}

\begin{figure}[ht]
\centering
\subfloat[Eigenvalues, $\varepsilon=0.01$]{\label{fig:eval0}
\includegraphics[width=0.3\linewidth]{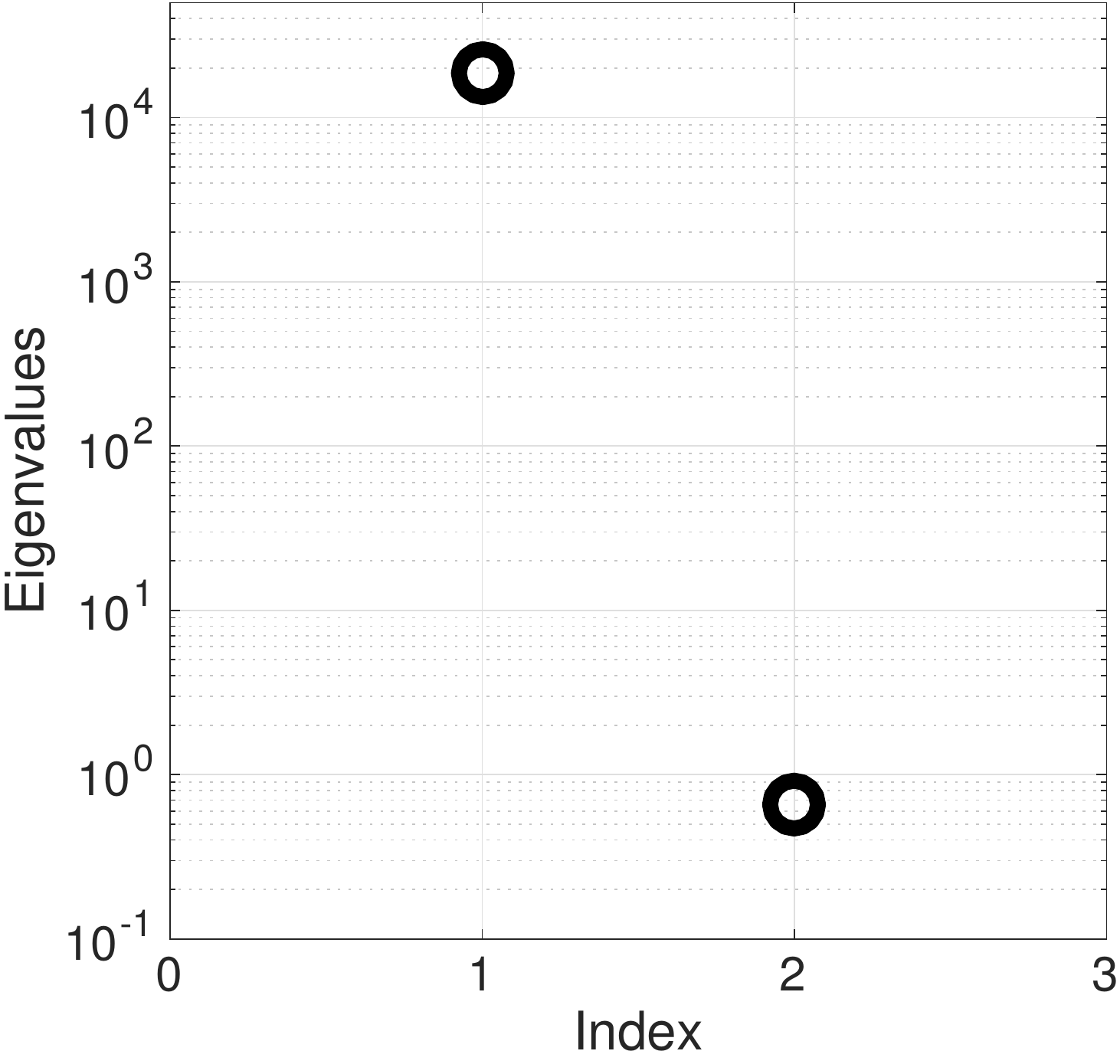}
}\;
\subfloat[Posterior contours, $\varepsilon=0.01$]{\label{fig:post0}
\includegraphics[width=0.3\linewidth]{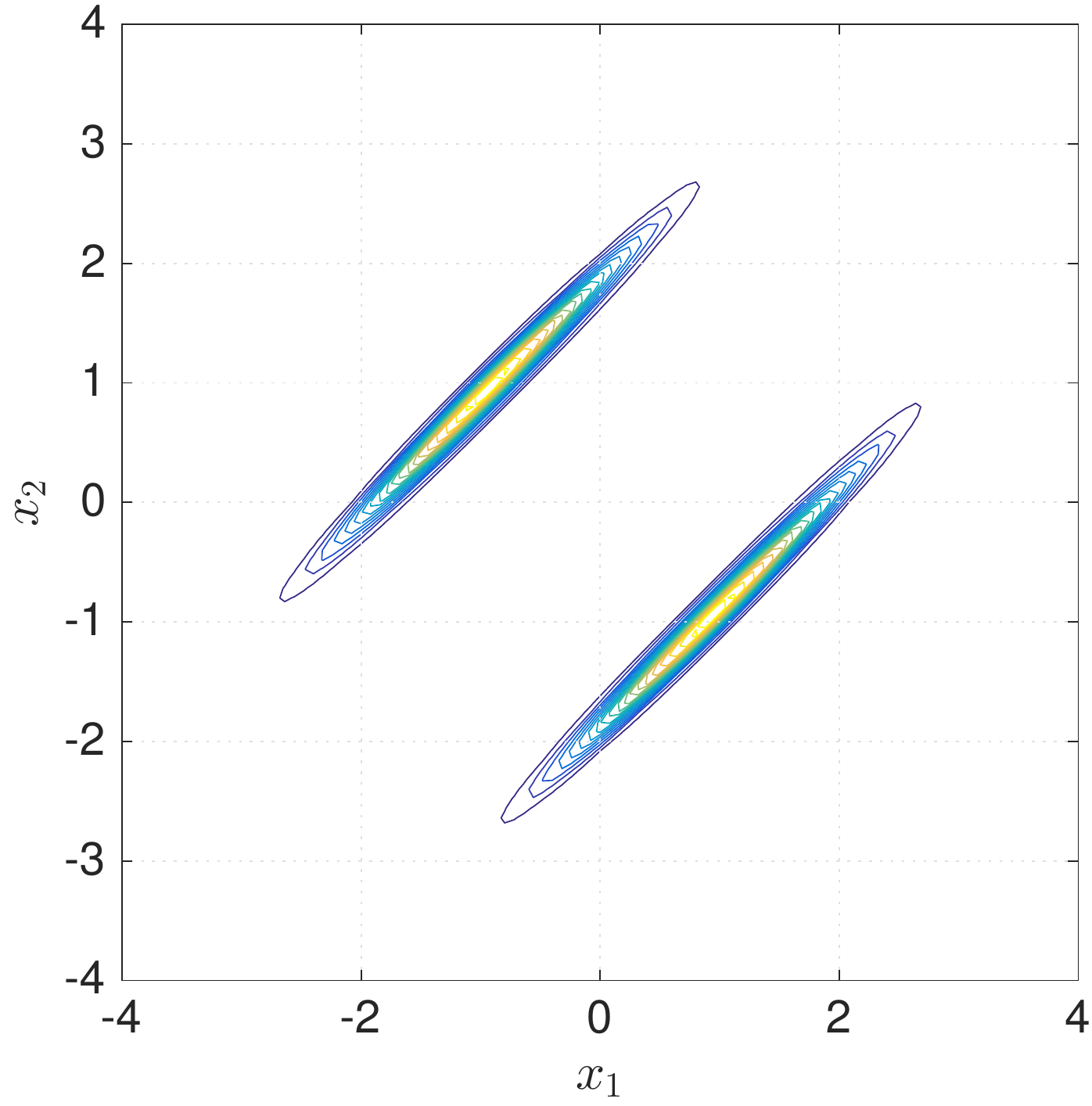}
}\\
\subfloat[Eigenvalues, $\varepsilon=0.95$]{\label{fig:eval1}
\includegraphics[width=0.3\linewidth]{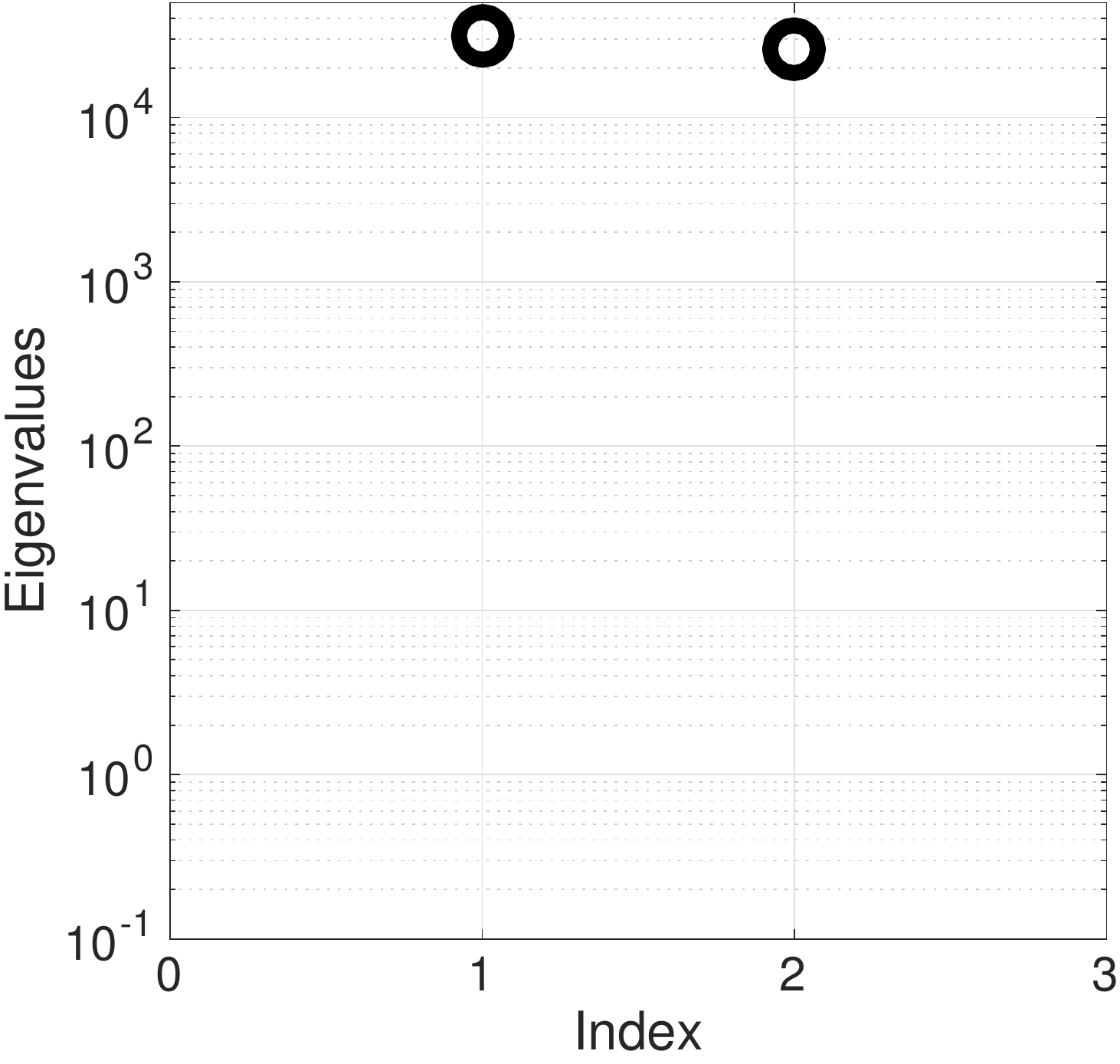}
}\;
\subfloat[Posterior contours, $\varepsilon=0.95$]{\label{fig:post1}
\includegraphics[width=0.3\linewidth]{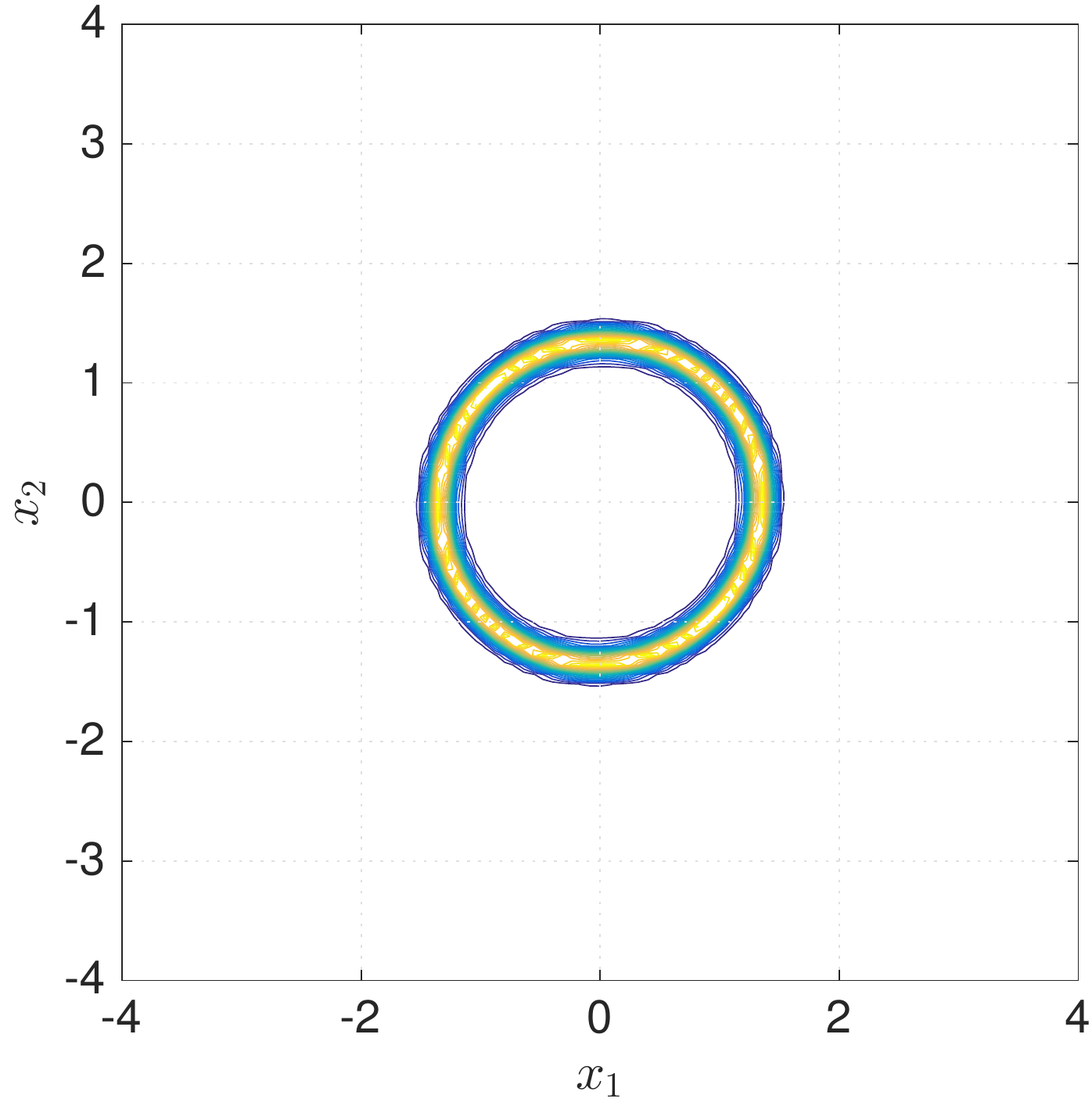}
}
\caption{Two cases of $\mA$ in the forward map \eqref{eq:quadfwd}. The first uses $\varepsilon=0.01$ in \eqref{eq:amat} (top row) and the second uses $\varepsilon=0.95$ (bottom row). The leftmost figures show the two eigenvalues from $\mC$ in \eqref{eq:C} computed with a high order Gauss-Hermite quadrature rule. The rightmost figures show the contours of the posterior density conditioned on $\vd=0.9$. (Colors are visible in the electronic version.)}
\label{fig:smallex}
\end{figure}

We first demonstrate the approach on a simple two-parameter quadratic forward model with a scalar output; in this case $\vd$ has one component. Let $\vx=[x_1,x_2]^T$, and we assume a standard Gaussian prior. Define the parameter-to-observable map $\vm(\vx)$ as
\begin{equation}
\label{eq:quadfwd}
\vm(\vx) \;=\; \frac{1}{2}\,\vx^T\mA\vx,
\end{equation}
where 
\begin{equation}
\label{eq:amat}
\mA = \mQ\,\bmat{1 & \\ & \varepsilon}\,\mQ^T,\qquad
\mQ = \frac{1}{2}\bmat{\sqrt{2} &\sqrt{2}\\ -\sqrt{2} & \sqrt{2}}.
\end{equation}
This map goes into the likelihood \eqref{eq:likelihood} with noise parameter $\sigma^2=0.1$, and we set the fixed data $\vd=0.9$. The parameter $\varepsilon$ in \eqref{eq:amat} controls how active the one-dimensional active subspace is. Figure \ref{fig:smallex} shows two cases: (i) $\varepsilon=0.01$ and (ii) $\varepsilon=0.95$. We compute the elements of $\mC$---each a two-dimensional integral---using a tensor product Gauss-Hermite quadrature rule with 50 points in each dimension. The high order integration rule eliminates the finite sampling errors present in Monte Carlo estimates. The eigenvalues of $\mC$ with $\varepsilon=0.01$ are shown in Figure \ref{fig:eval0}. The gap between the eigenvalues suggests an active one-dimensional active subspace. This is confirmed in Figure \ref{fig:post0}, which plots the skewed contours of the posterior density conditioned on $\vd=0.9$. The comparable figures for $\varepsilon=0.95$ are shown in Figures \ref{fig:eval1} and \ref{fig:post1}, respectively. In this case, both eigenvalues have the same order of magnitude and the smallest eigenvalue is large, so they do not suggest an exploitable active subspace. This is confirmed by the posterior contours that vary significantly in all directions. 

\begin{figure}[ht]
\centering
\subfloat[MCMC posterior $\vx$]{\label{fig:mcx}
\includegraphics[width=0.3\linewidth]{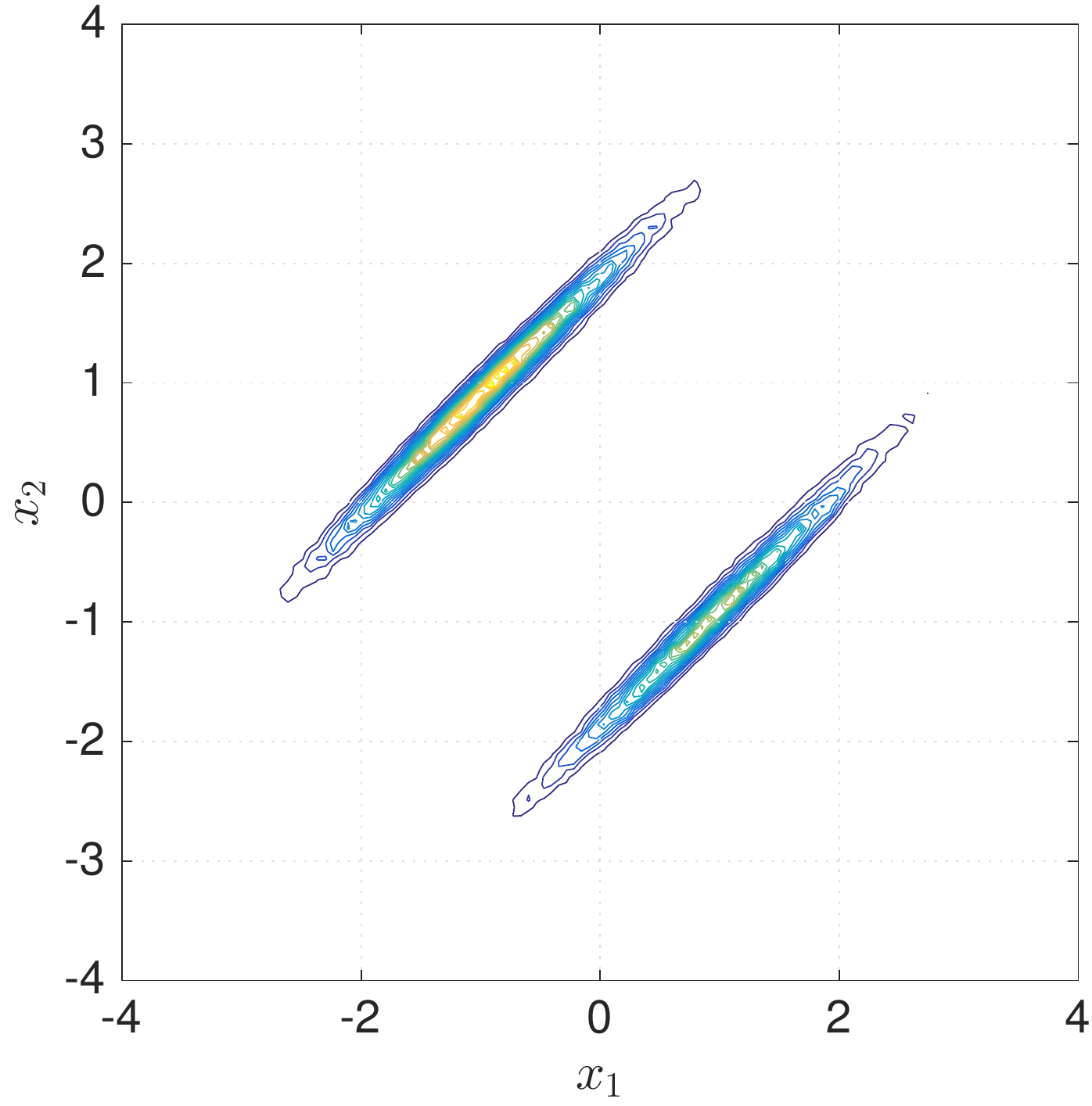}
}\;
\subfloat[MCMC-AS posterior $\vx$]{\label{fig:mcasx}
\includegraphics[width=0.3\linewidth]{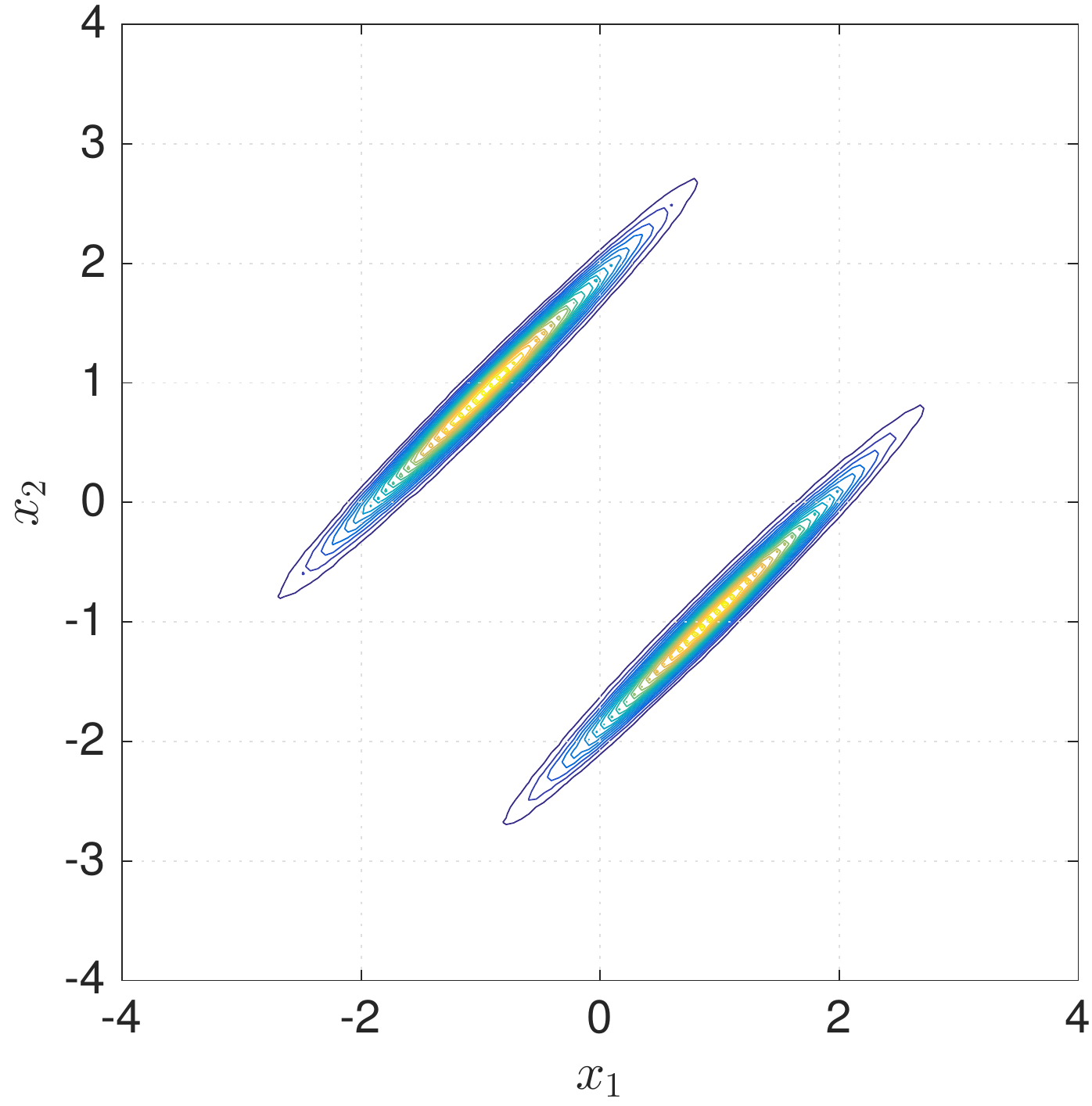}
}\;
\subfloat[MCMC-AS posterior $\vy$]{\label{fig:mcasy}
\includegraphics[width=0.3\linewidth]{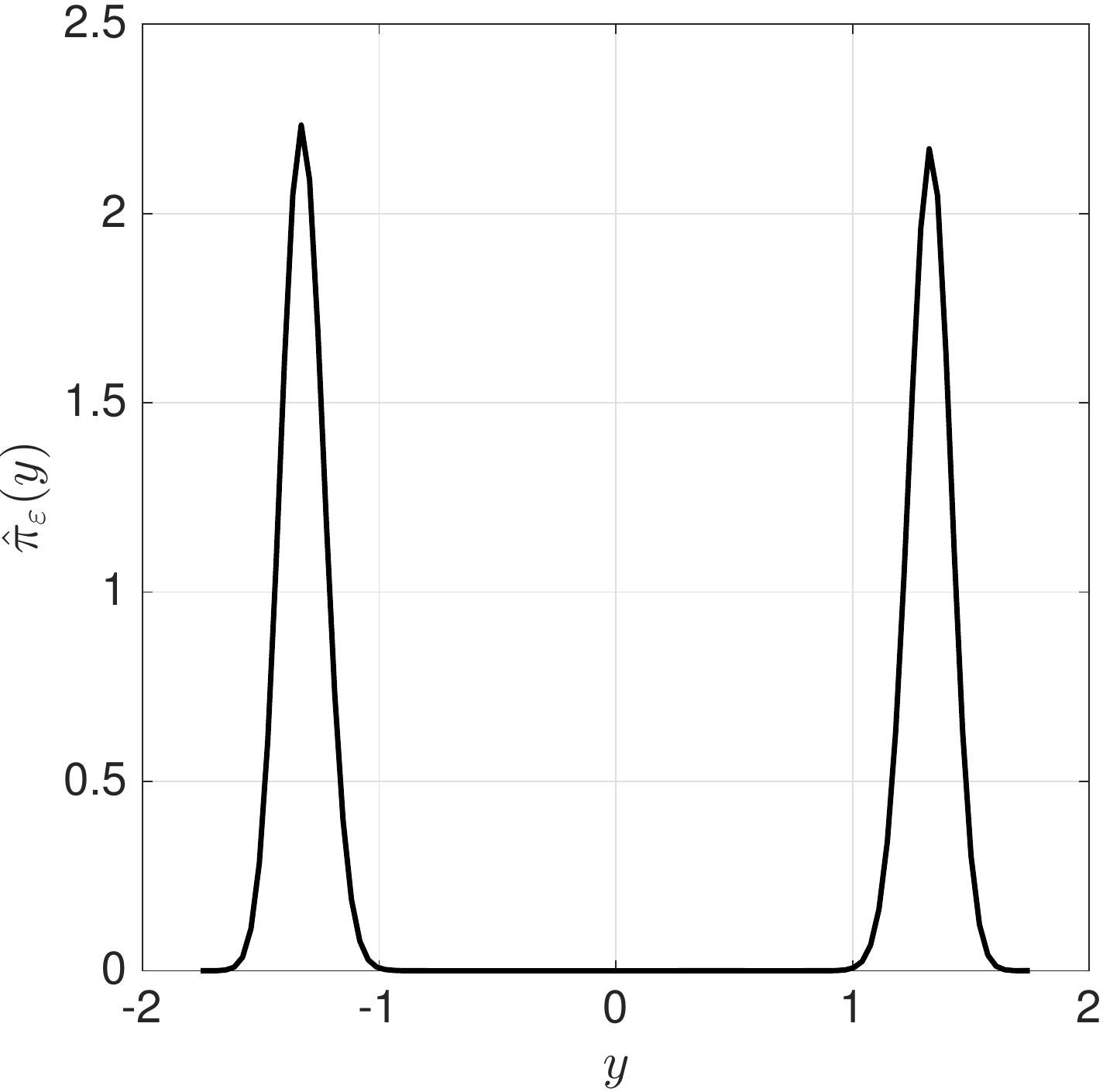}
}\\
\subfloat[MCMC trace $\vx$]{\label{fig:trx}
\includegraphics[width=0.3\linewidth]{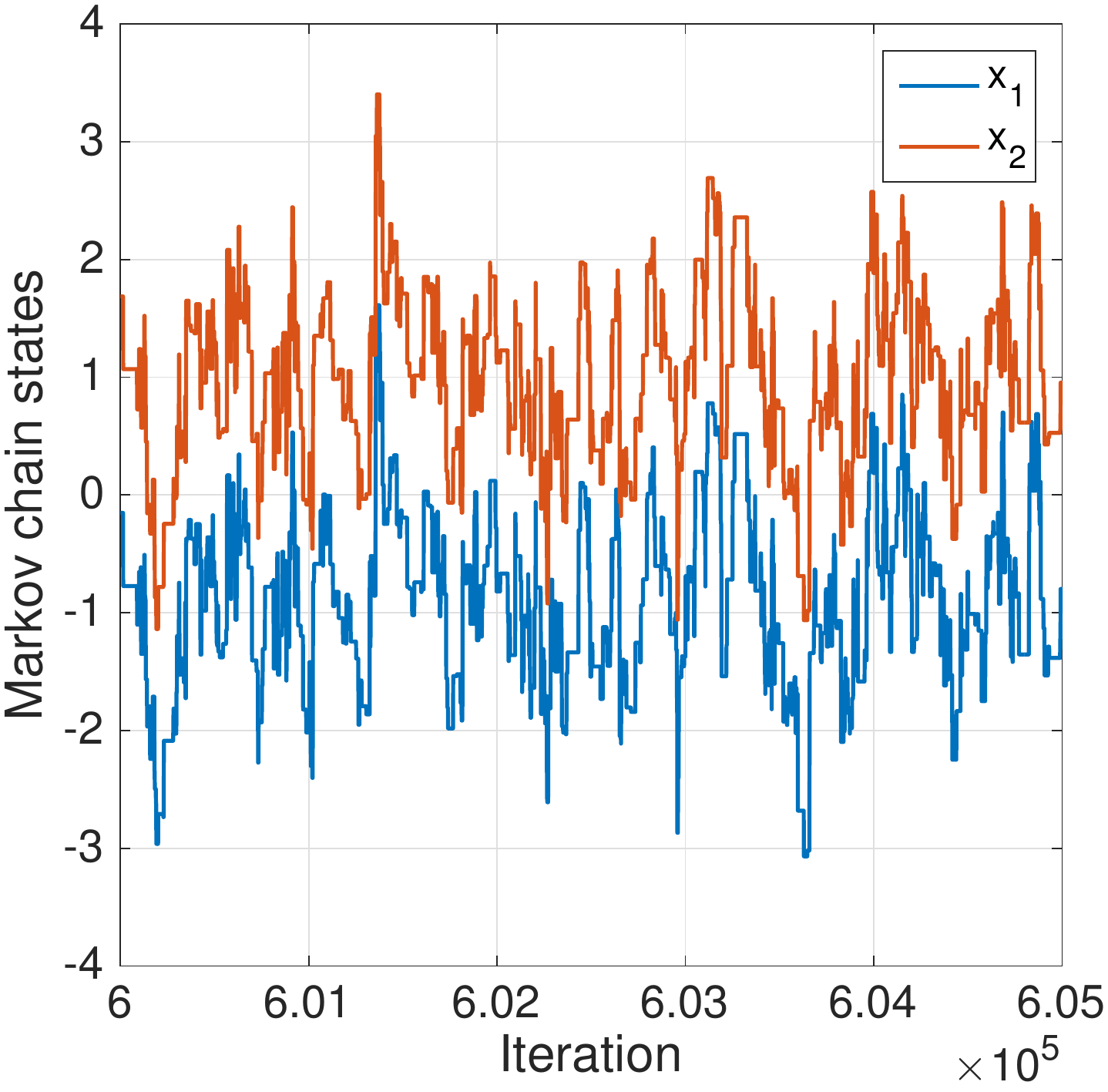}
}\:
\subfloat[MCMC-AS trace $\vx$]{\label{fig:trasx}
\includegraphics[width=0.3\linewidth]{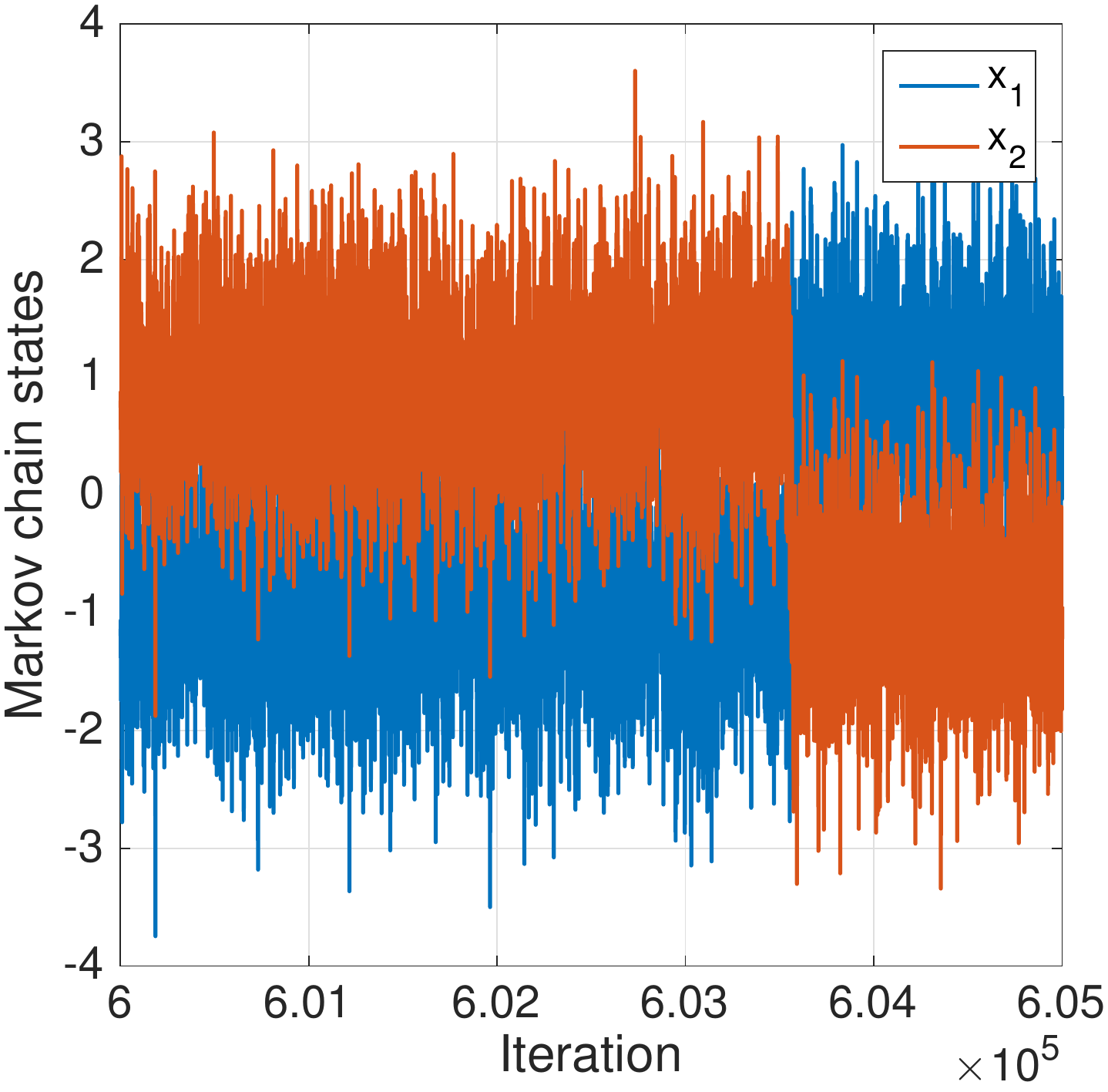}
}\;
\subfloat[MCMC-AS trace $\vy$]{\label{fig:trasy}
\includegraphics[width=0.3\linewidth]{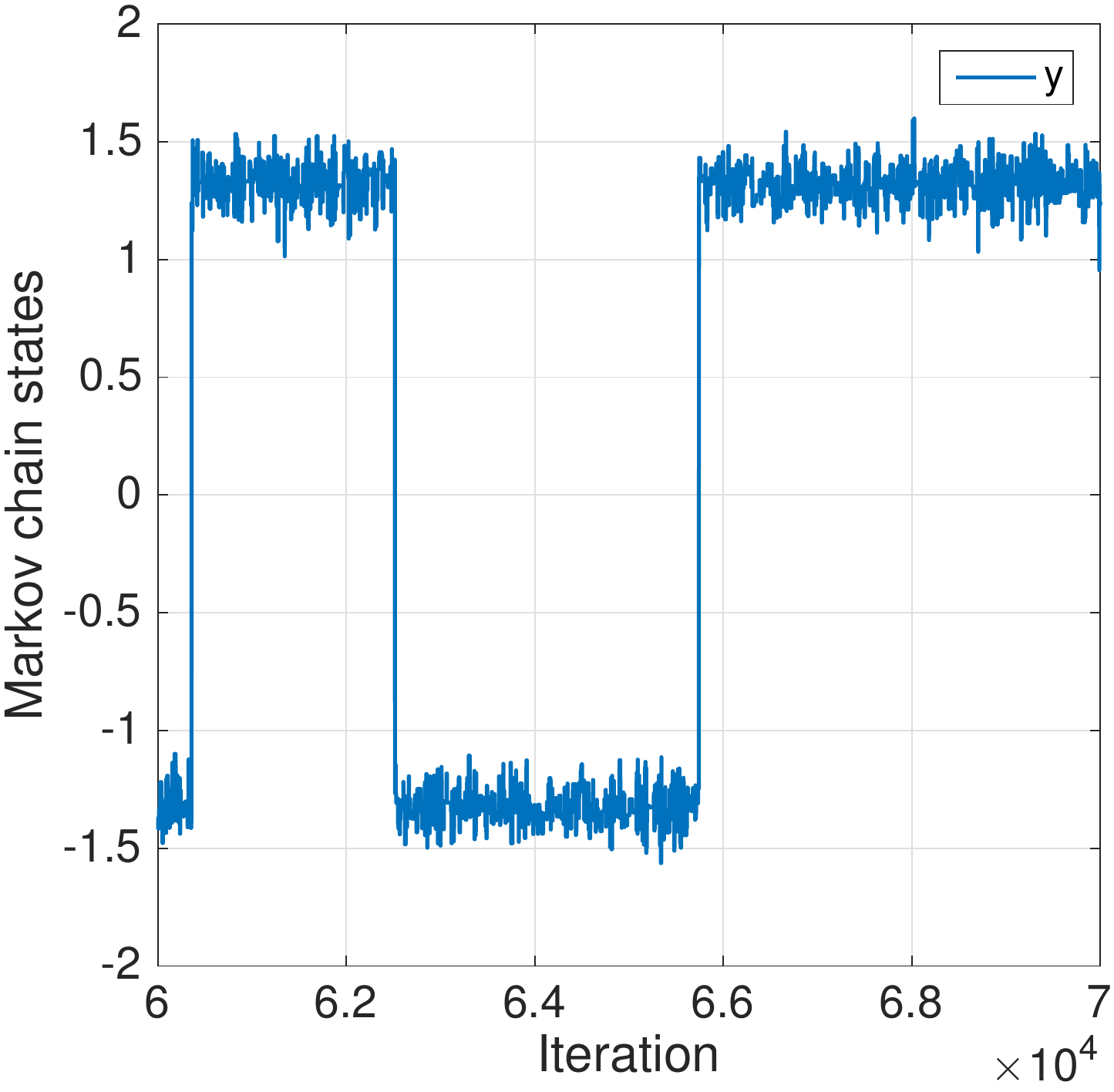}
}
\caption{Comparison of standard Metropolis-Hastings MCMC to the active subspace accelerated MCMC. Figure \ref{fig:mcx} shows density estimate of the posterior contours for the standard MCMC, and Figure \ref{fig:mcasx} shows the same contours for the samples from Algorithm \ref{alg:mcmc} using the reconstruction \eqref{eq:reconstruct}. The plots directly beneath them show trace plots for steps 600000 to 605000; notice that the active subspace-accelerated chain mixes better. The rightmost plots show the univariate posterior on the active variable $\hvy$ (above) and its trace plot for 10000 representative steps (below). (Colors are visible in the electronic version.)}
\label{fig:smallexmc}
\end{figure}

Next we compare standard Metropolis-Hastings MCMC on the two-dimensional parameter space to an active subspace-accelerated MCMC; we use the $\varepsilon=0.01$ case of the forward model \eqref{eq:quadfwd} with data $\vd=0.9$. The active subspace-accelerated method (i) runs Metropolis-Hastings on the estimated active variables $\hvy$ and (ii) samples independently from the prior on the inactive variables $\hvz$ for each sample of $\hvy$. In this case, $\hvy$ and $\hvz$ each have one component. The bivariate standard Gaussian prior on $\vx$ implies that the conditional density of $\hvz$ given $\hvy$ is a univariate standard Gaussian prior. 

We run the standard MCMC for $10^6$ steps using a two-dimensional standard Gaussian proposal density with variance $\sigma^2=0.5$. The acceptance rate was 12\%, which is not surprising since the true posterior in Figure \ref{fig:post0} shows two small, separated regions of large posterior probability. Figure \ref{fig:mcx} shows the contours of a bivariate Gaussian kernel density estimate~\cite{kde2d} using all $10^6$ steps. Figure \ref{fig:trx} shows the trace plots for steps 600000 to 605000. 

We run the active subspace-accelerated MCMC for $10^5$ steps. At each step, we evaluate the misfit approximation $\hgeps(\hvy)$ with a 10-point Gauss-Hermite quadrature rule, which is appropriate because (i) $\hvz$ has one component and (ii) the conditional density of $\hvz$ given $\hvy$ is a standard Gaussian. Since each approximate misfit uses 10 forward model evaluations, the total number of forward model evaluations is the same for each MCMC method. The proposal density for Metropolis-Hastings on $\hvy$ is a univariate standard Gaussian with variance $\sigma^2=0.5$, and the acceptance rate was 13\%. We draw $P=10$ independent samples from the Gaussian prior on $\hvz$ for each step in the chain as in \eqref{eq:reconstruct}; recall that drawing these samples needs no forward model evaluations. Figure \ref{fig:mcasx} shows the contours of the approximate posterior density constructed with the same bivariate Gaussian kernel density estimation from samples of the active subspace-accelerated MCMC; these contours compare favorably to the true posterior contours in Figure \ref{fig:post0}. The trace plots for the $\vx$ components are shown in Figure \ref{fig:trasx} for steps 600000 to 605000. Notice the difference between the trace plots in Figure \ref{fig:trx} and those in Figure \ref{fig:trasx}. Since $\hvz$ are drawn independently according to their Gaussian prior, the samples from the active subspace-accelerated MCMC are much less correlated than the samples from MCMC on the two-dimensional $\vx$ space. Figure \ref{fig:mcasy} shows a univariate Gaussian kernel density estimate of the posterior on $\hvy$ computed with the $10^5$ samples from the Metropolis-Hastings; Figure \ref{fig:trasy} shows the trace plot for steps 60001 to 70000. 

In this small example where the misfit admits a one-dimensional active subspace, the active subspace-accelerated MCMC performs very well. For the same number of forward model evaluations, the active subspace enables the Metropolis-Hastings to run on only the active variable $\hvy$, while the inactive variable $\hvz$ is sampled independently according to its prior. This procedure produces samples of $\vx$ with very little correlation. 

\subsection{PDE model}
\label{sec:pdebayes}

\begin{figure}[ht]
\centering
\subfloat{
\includegraphics[width=0.3\linewidth]{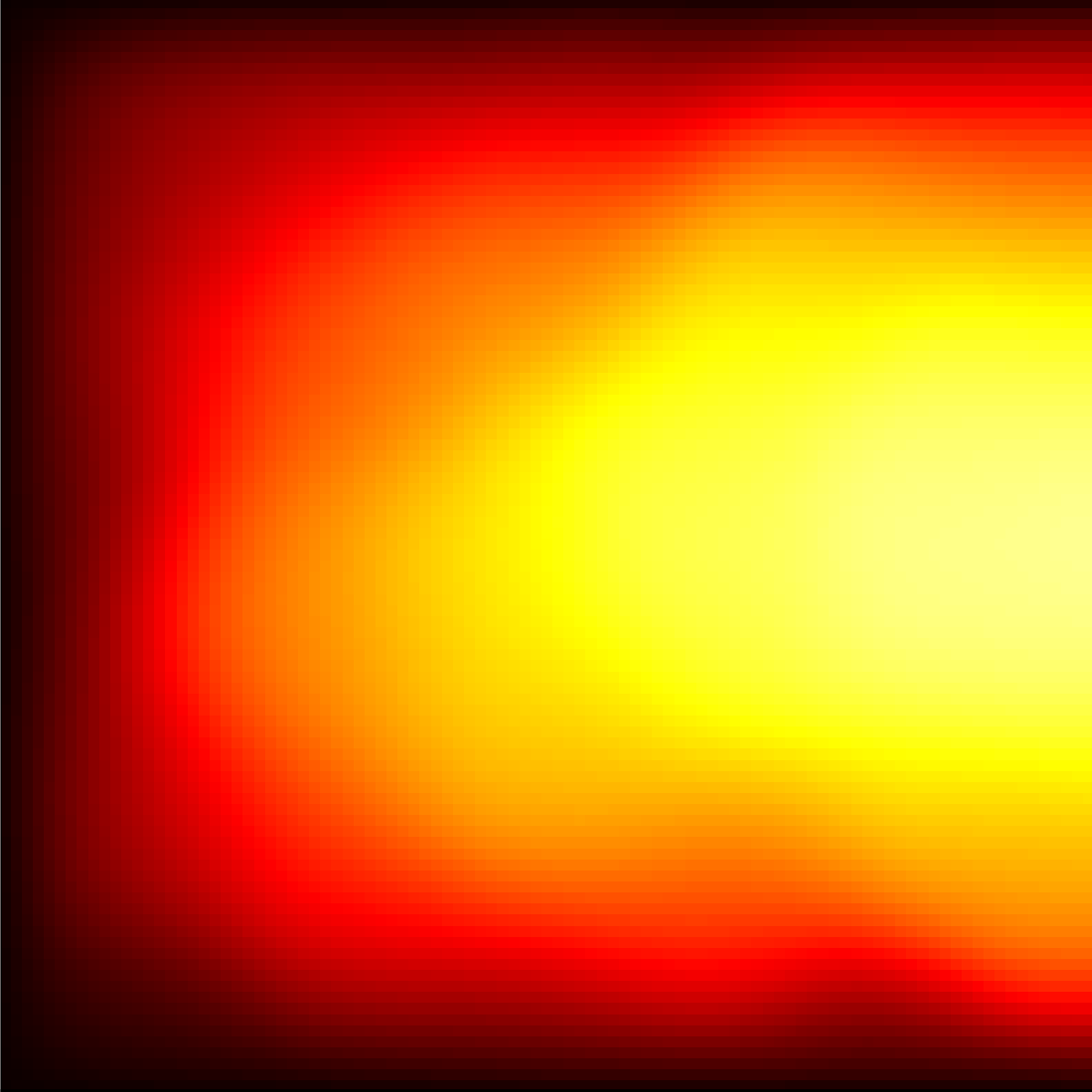}
}\;
\subfloat{
\includegraphics[width=0.3\linewidth]{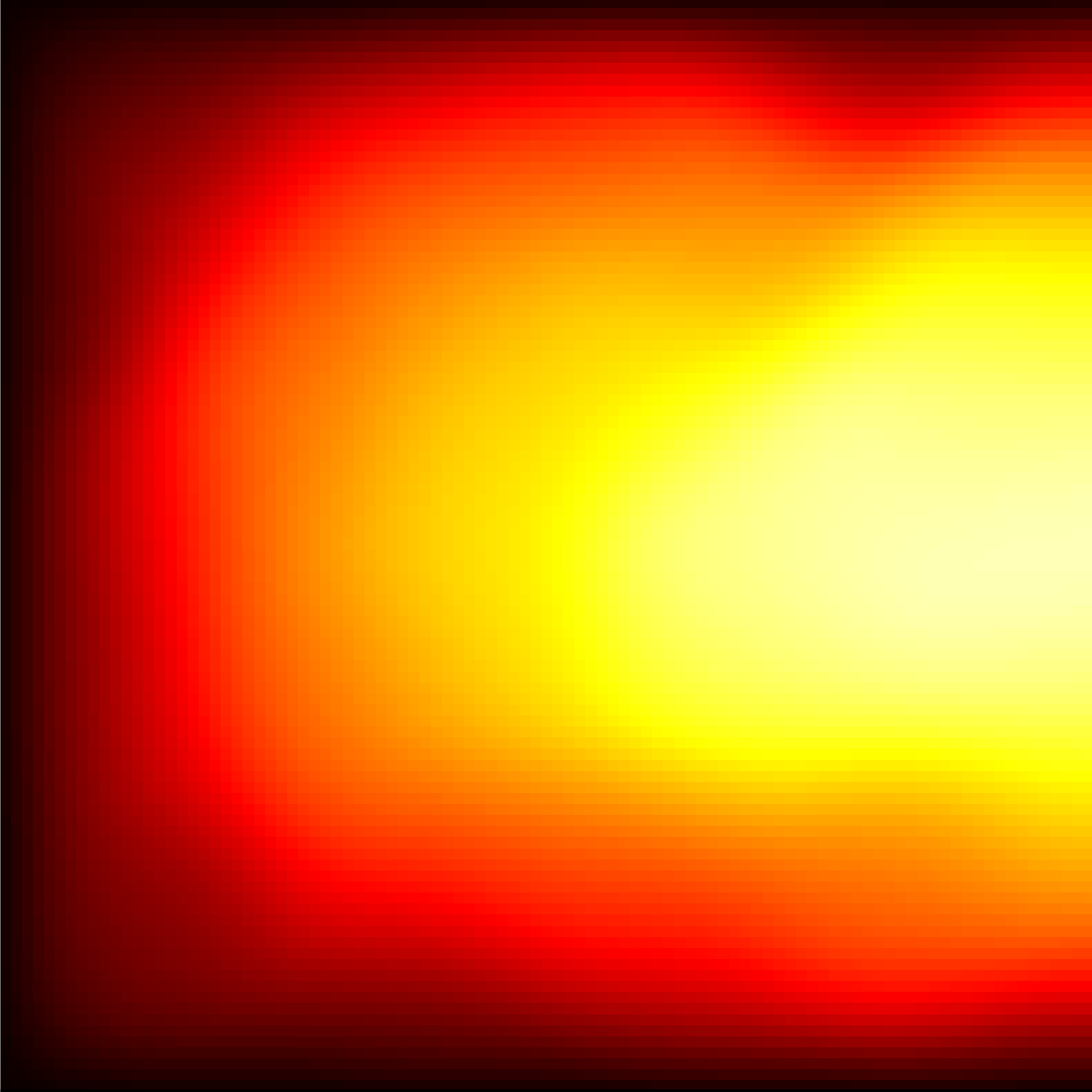}
}\;
\subfloat{
\includegraphics[width=0.3\linewidth]{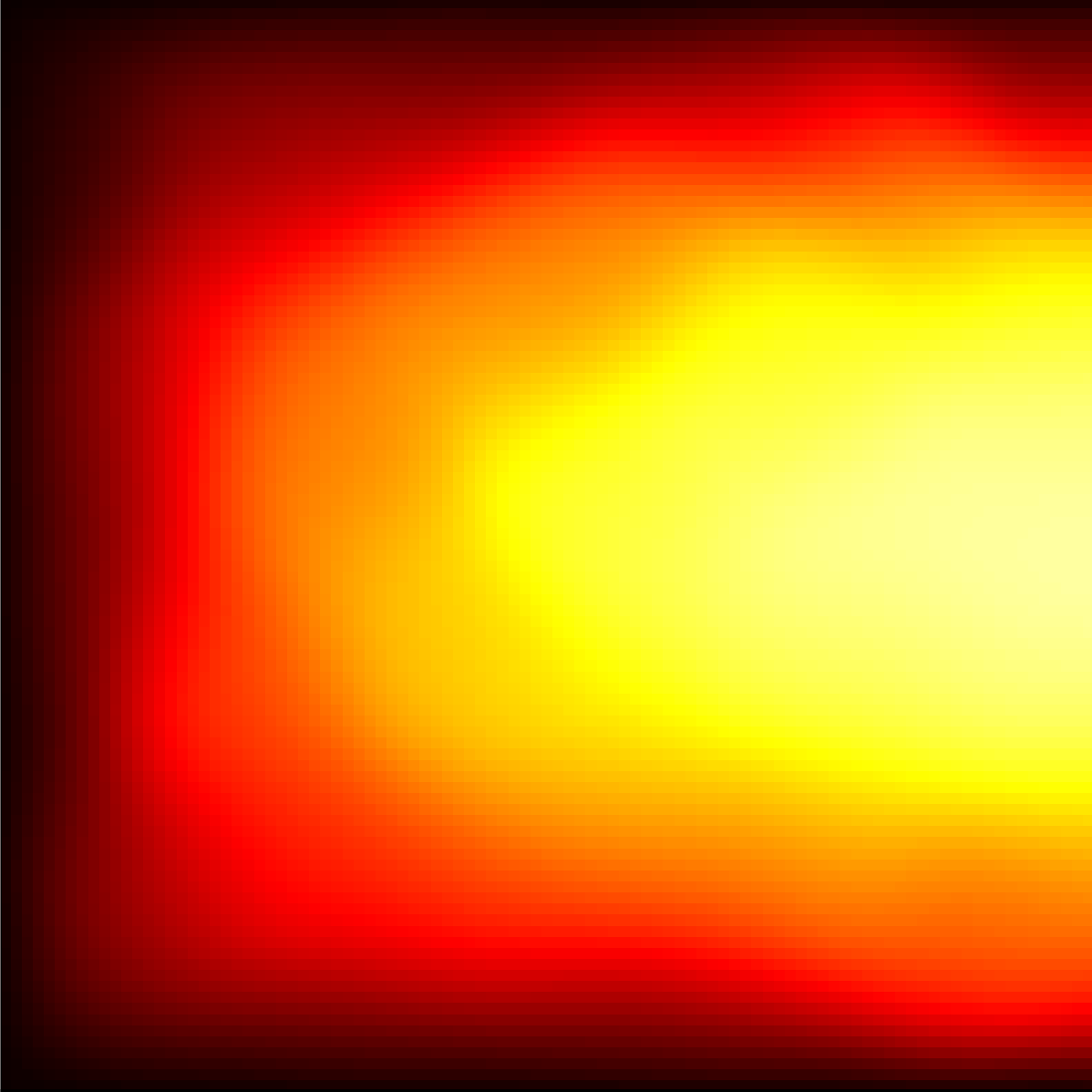}
}\\
\subfloat{
\includegraphics[width=0.3\linewidth]{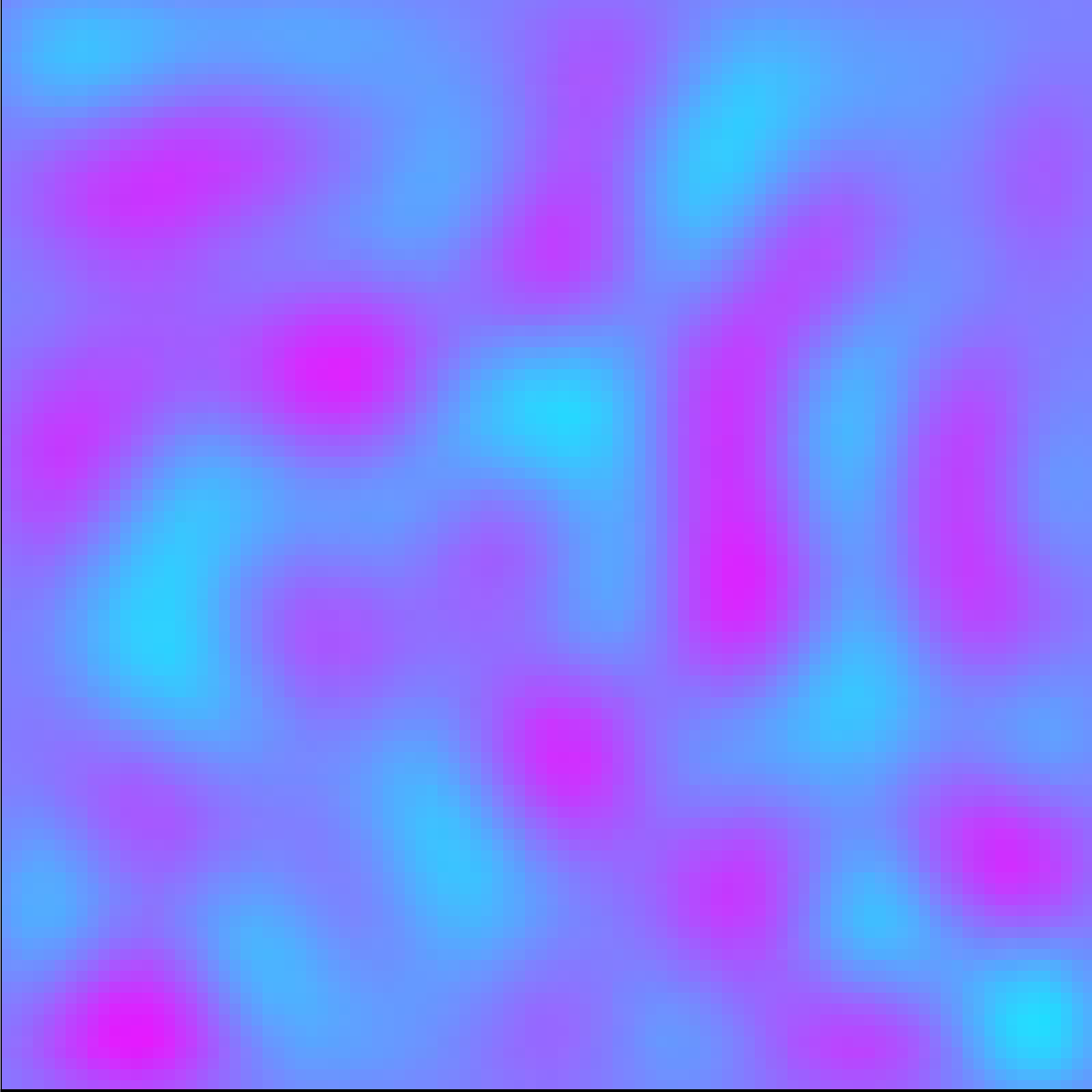}
}\:
\subfloat{
\includegraphics[width=0.3\linewidth]{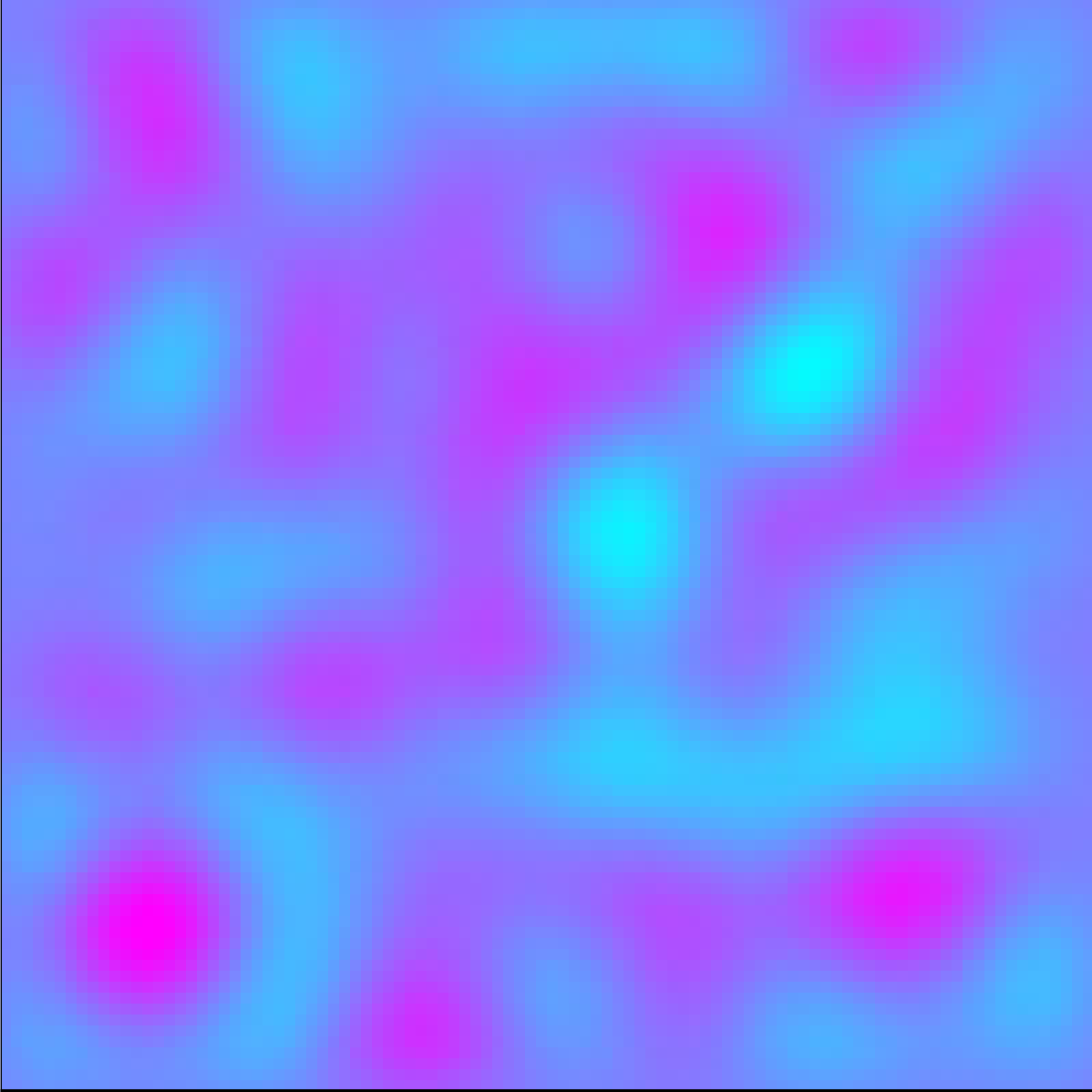}
}\;
\subfloat{
\includegraphics[width=0.3\linewidth]{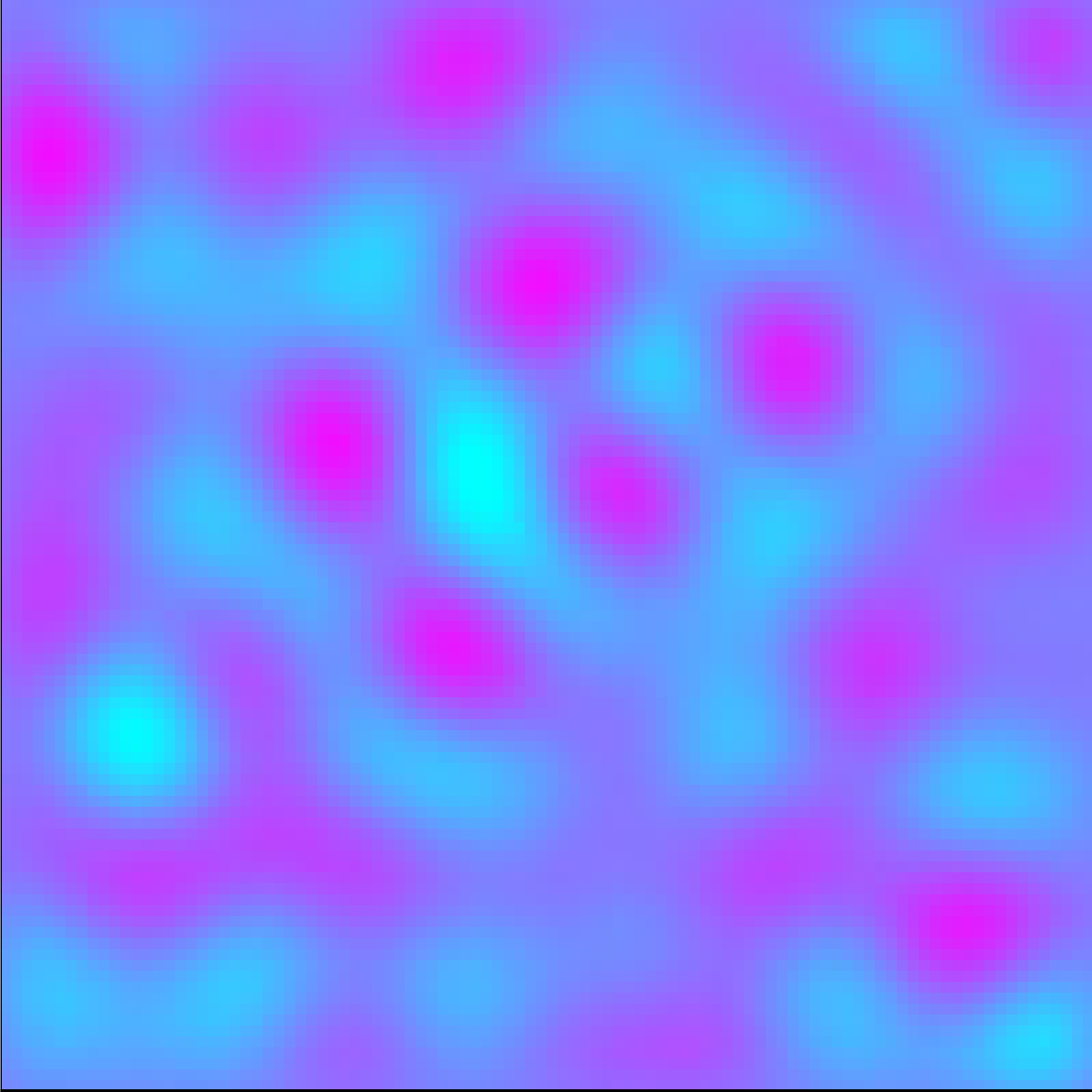}
}
\caption{The bottom row shows three realizations of the coefficients $\log(a)$ with parameters $\vx$ drawn according to the prior. The top row shows the corresponding PDE solutions. (Colors are visible in the electronic version.)}
\label{fig:1}
\end{figure} 

Next we apply the active subspace-accelerated MCMC on a synthetic inverse problem derived from a PDE forward model. The PDE solution $u=u(\vs,\vx)$ satisfies
\begin{equation}
\label{eq:pde}
-\nabla_{\vs}\cdot(a\,\nabla_{\vs}u) \;=\; 1, \qquad \vs\in[0,1]^2.
\end{equation}
The coordinates $\vs$ represent the solution's spatial variables, and $\nabla_{\vs}$ denotes differentiation with respect to $\vs$. The boundary conditions are $u=0$ on the bottom ($s_2=0$), left ($s_1=0$), and top ($s_2=1$) of the domain, and $\vn\cdot(a\,\nabla_{\vs}u)=0$ on the right ($s_1=1$) of the domain, where $\vn$ is the unit normal vector. The log of the coefficients $a = a(\vs,\vx)$ is a truncated Karhunen-Loeve series,
\begin{equation}
\label{eq:kl}
\log(a(\vs,\vx)) \;=\; \sum_{i=1}^m \sqrt{\sigma_i}\,\phi_i(\vs)\,x_i,
\end{equation}
with $m=100$ terms. The pairs $(\sigma_i,\phi_i(\vs))$ are the eigenpairs of an exponential two-point correlation function,
\begin{equation}
\label{eq:corr}
\sC(\vs_1,\vs_2) \;=\; \exp\left(-\frac{\|\vs_1-\vs_2\|_1}{\beta}\right),
\end{equation}
where $\|\cdot\|_1$ is the 1-norm on $\mathbb{R}^2$, and $\beta=0.02$ is the correlation length. Note that the correlation length 0.02 is short compared to the size of the spatial domain, and the 1-norm generates rough coefficient fields. The Karhunen-Loeve eigenvalues $\sigma_i$ are shown in Figure \ref{fig:01}; their slow decay reflects the short correlation length. Given a point $\vx=[x_1,\dots,x_m]^T$, the forward model evaluates the Karhunen-Loeve expansion \eqref{eq:kl} and solves the PDE \eqref{eq:pde} with second-order finite differences---discretized using 100 points in each spatial dimension. Figure \ref{fig:1} shows three realizations of the log-coefficients with parameters drawn according to a standard Gaussian (bottom row) and their corresponding PDE solutions (top row). 

\subsection{Parameters and data}

\begin{figure}[ht]
\centering
\subfloat[Forward model realizations and data]{
\label{fig:00}
\includegraphics[width=0.45\linewidth]{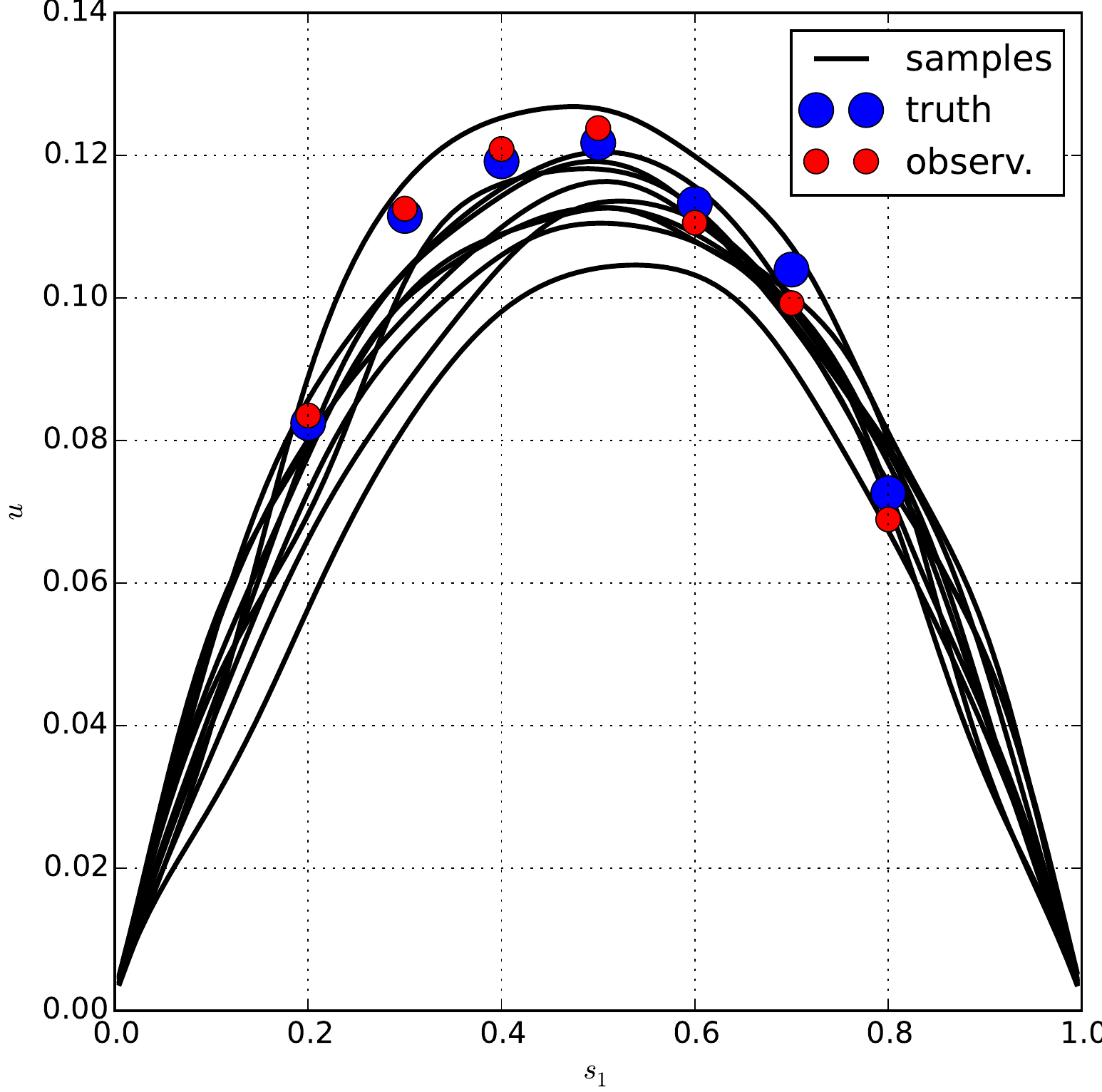}
}\;
\subfloat[Karhunen-Loeve eigenvalues]{
\label{fig:01}
\includegraphics[width=0.45\linewidth]{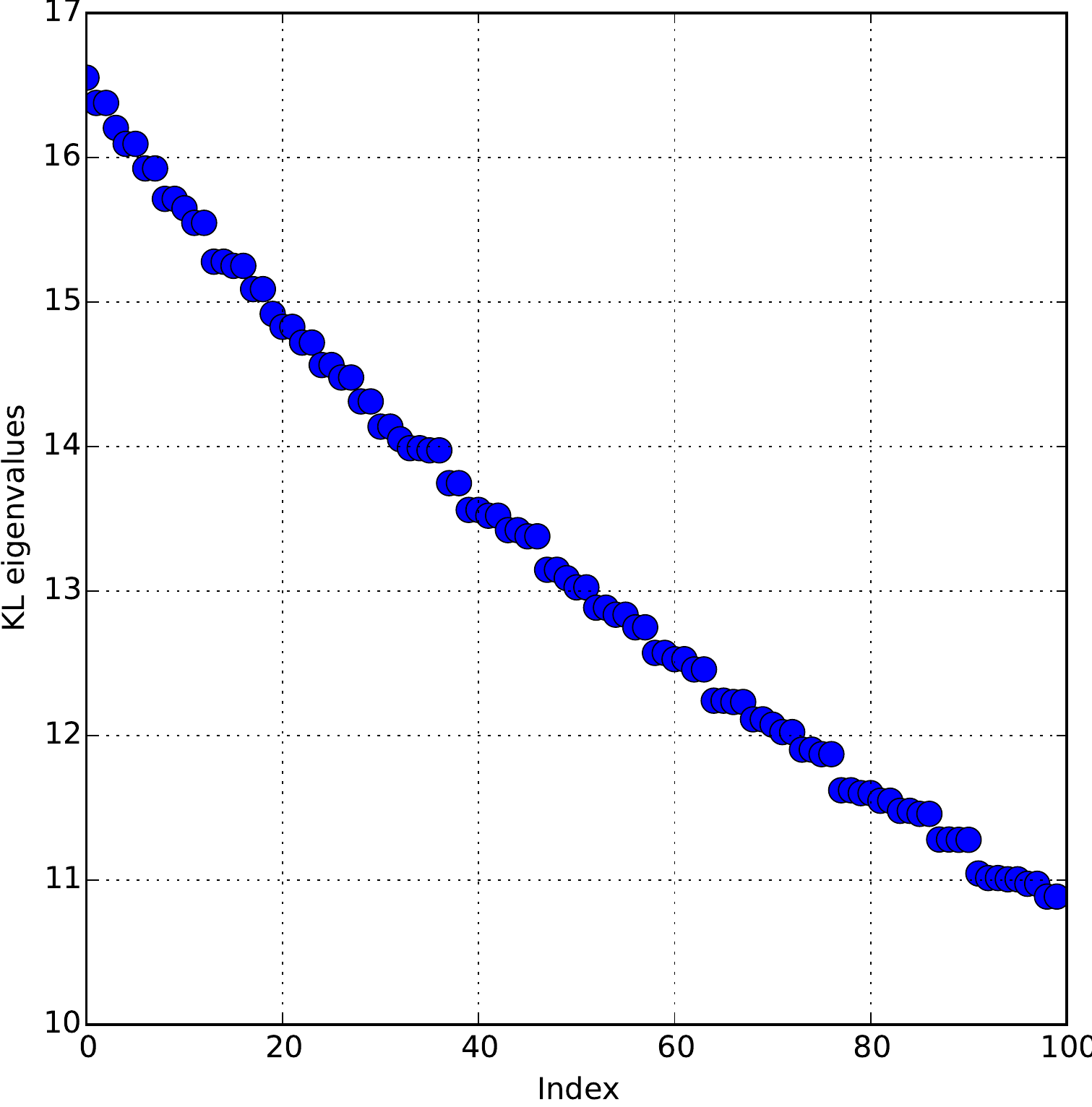}
}
\caption{On the left, Figure \ref{fig:00} shows ten realizations of the PDE solution with the parameters drawn according to the prior. The blue dots are the observations corresponding to $\vx=\xtrue$, and the red dots show the observations perturbed by 1\% noise to generate synthetic data for the inverse problem. Figure \ref{fig:01} shows the eigenvalues of the Karhunen-Loeve expansion in \eqref{eq:kl}. The slow decay corresponds to the short correlation length. (Colors are visible in the electronic version.)}
\label{fig:0}
\end{figure} 

The coefficients $\vx$ in \eqref{eq:kl} are parameters in the statistical inverse problem. We use a standard Gaussian prior on $\vx$, which implies a Gaussian random field prior on $\log(a)$ with a low-rank ($m$-term) correlation function approximating \eqref{eq:corr}. This setup is similar to the \emph{a priori} dimension reduction proposed by Marzouk and Najm~\cite{Marzouk2009}, which, for this case, reduces the dimension of the parameter space from 10000 to 100 via the truncated Karhunen-Loeve series.

The forward model's observations come from the PDE solution $u$ evaluated at seven points on the right boundary, $s_2\in\{0.2,0.3,0.4,0.5,0.6,0.7,0.8\}$. The parameter-to-observable map $\vm(\vx)$ from \eqref{eq:addnoise} takes the Karhunen-Loeve coefficients $\vx$ and returns the PDE solution at those seven points on the boundary. To generate the data for the synthetic inverse problem, we draw $\xtrue$ according to the prior and evaluate \eqref{eq:addnoise} with $\vx=\xtrue$. The noise is $\sigma^2=0.0001\,\|\vm(\xtrue)\|^2_2$, which is roughly 1\%. In Figure \ref{fig:00}, the black lines are ten realizations of the PDE solution $u$ on the right boundary ($s_1=1$) with parameters $\vx$ drawn according to the prior. The blue dots are the observations from the solution evaluated with $\vx=\xtrue$, and the red dots are the observations perturbed by the noise that constitute the synthetic data for the inverse problem. 

\subsection{Estimating the active subspace}

\begin{figure}[ht]
\centering
\subfloat[Active subspace eigenvalues]{
\label{fig:20}
\includegraphics[width=0.45\linewidth]{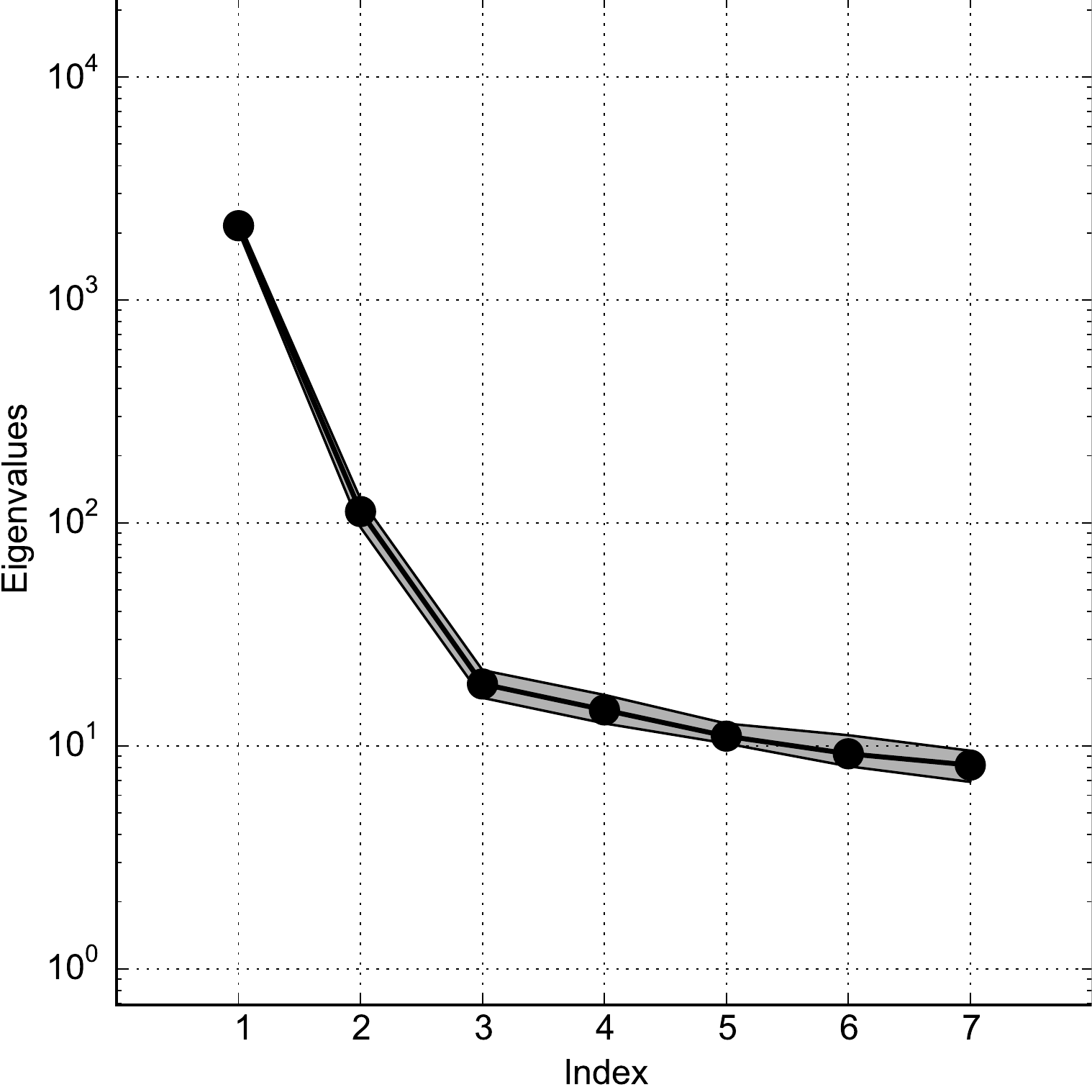}
}\;
\subfloat[Active subspace errors]{
\label{fig:21}
\includegraphics[width=0.45\linewidth]{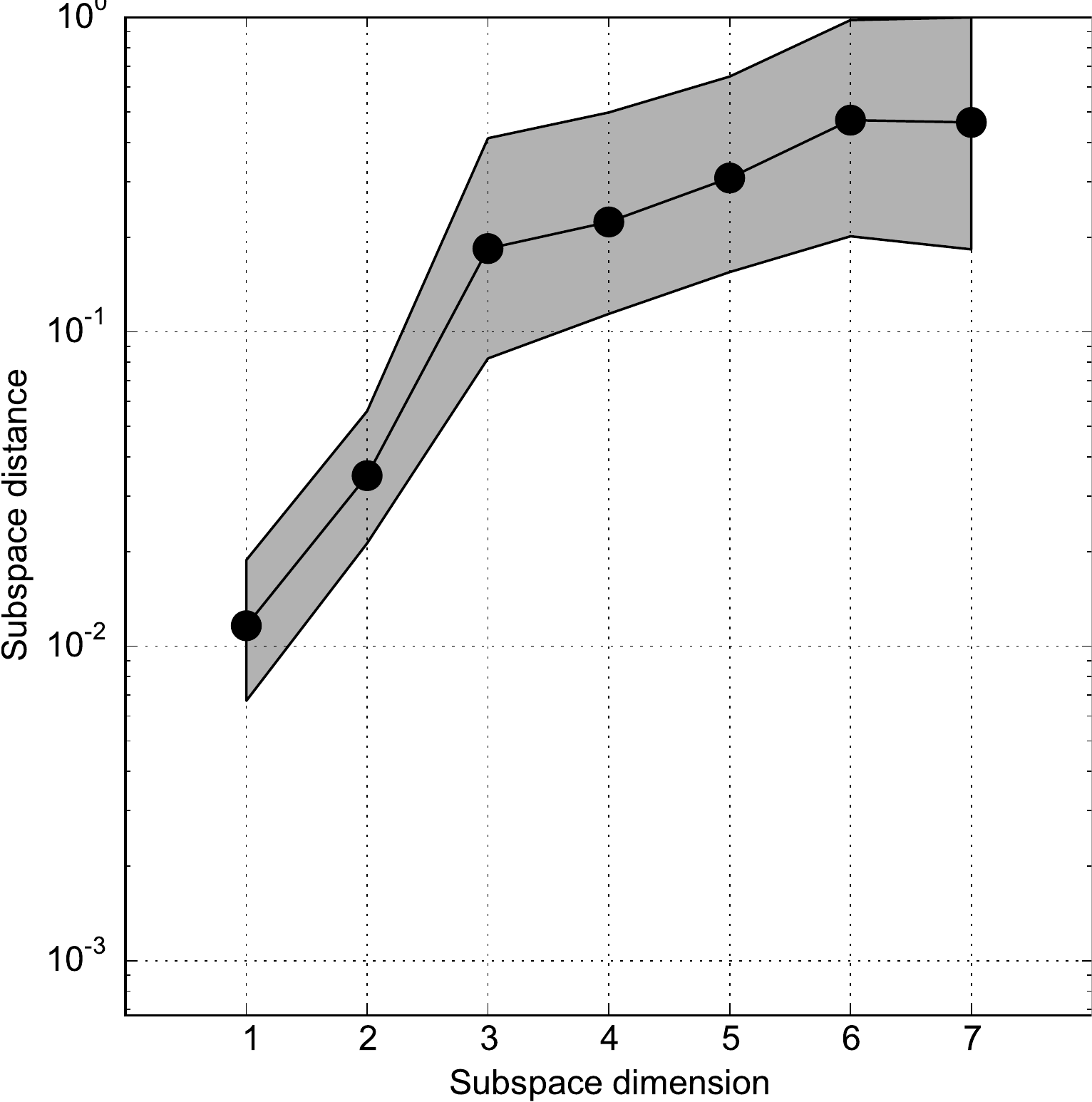}
}
\caption{The left figure shows the first seven of the $m=100$ eigenvalues $\hlambda_i$ from the active subspace analysis; see \eqref{eq:C}. The large gaps between eigenvalues $\hlambda_1$, $\hlambda_2$ and eigenvalues $\hlambda_2$, $\hlambda_3$ indicate a one-dimensional and a two-dimensional active subspace. The right figure shows the estimated subspace error $\varepsilon$ from \eqref{eq:suberr} using the bootstrap as described in~\cite{constantine2015computing}.}
\label{fig:2}
\end{figure} 

The Python code that evaluates the forward model is built using Wang's automatic differentiation package \texttt{numpad}~\cite{numpad}. This allows us to compute the misfit's gradient \eqref{eq:misfitgrad} efficiently. With the prior $\prior(\vx)$, misfit $f(\vx)$, and gradient $\nabla f(\vx)$ defined, we can estimate the active subspace. We use $N=1000$ gradient samples in the computation \eqref{eq:mc}. Figure \ref{fig:20} shows the first seven of the $m=100$ eigenvalues and their bootstrap ranges. The gaps between the first two pairs of eigenvalues suggest that the Monte Carlo approximation \eqref{eq:mc} can accurately estimate a one- or two-dimensional active subspace. Figure \ref{fig:21} shows bootstrap estimates of the subspace error \eqref{eq:suberr} as described in~\cite{constantine2015computing}. The relatively small error estimates for the one- and two-dimensional subspaces are consistent with the eigenvalue gaps. We choose to study a two-dimensional active subspace, since (i) the gap between eigenvalues $\hlambda_2$ and $\hlambda_3$ and the estimated subspace error suggest a good approximation of the subspace, and (ii) the error estimate in Theorem \ref{thm:postapprox} then includes the last 98 eigenvalues instead of the last 99 for a one-dimensional subspace. In other words, the theoretical error bound for the approximate posterior is smaller for a two-dimensional active subspace than for a one-dimensional active subspace. 

\subsection{Spatial sensitivity}

\begin{figure}[ht]
\centering
\subfloat[First eigenvector]{
\label{fig:30}
\includegraphics[width=0.45\linewidth]{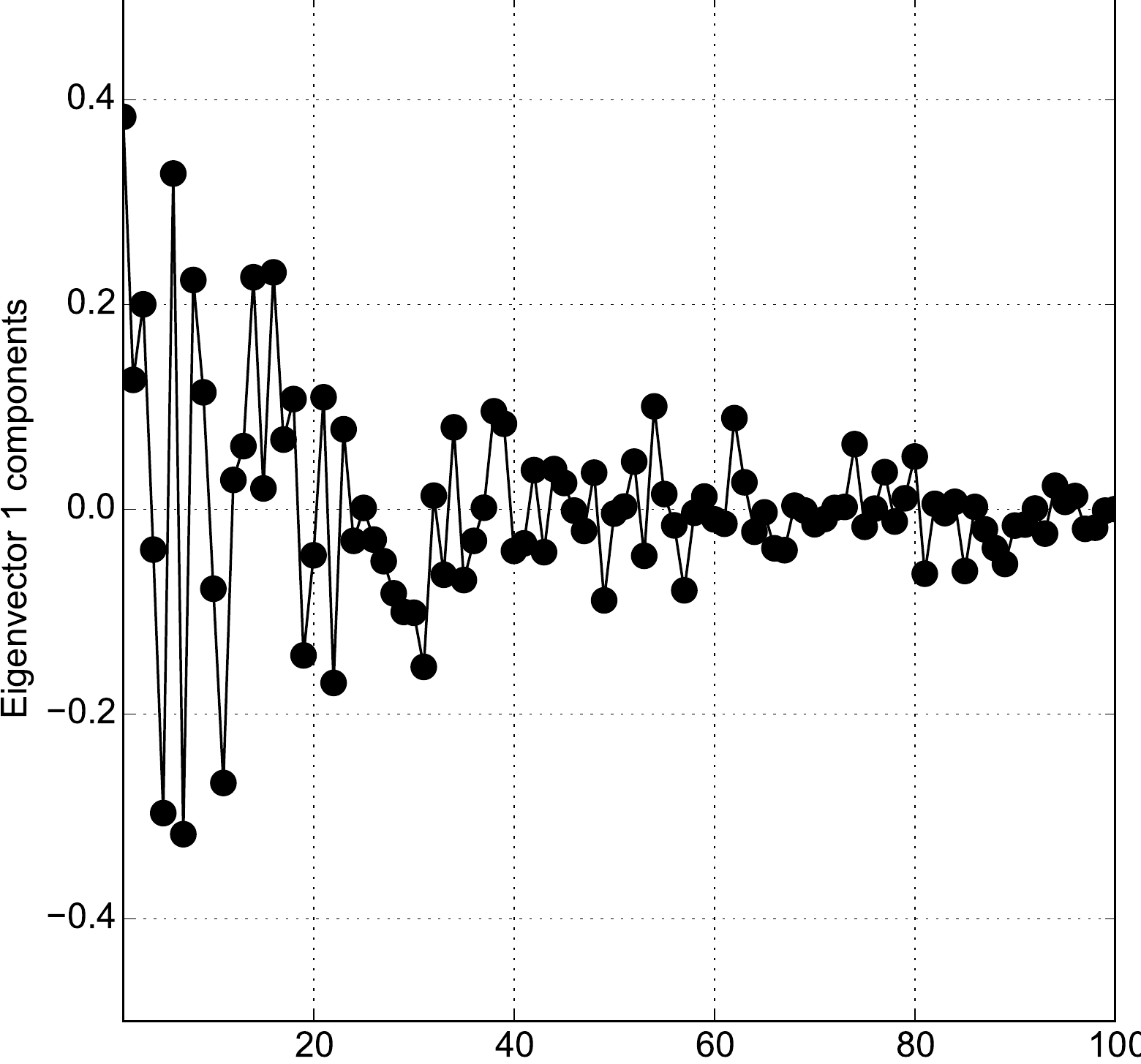}
}\;
\subfloat[Second eigenvector]{
\label{fig:31}
\includegraphics[width=0.45\linewidth]{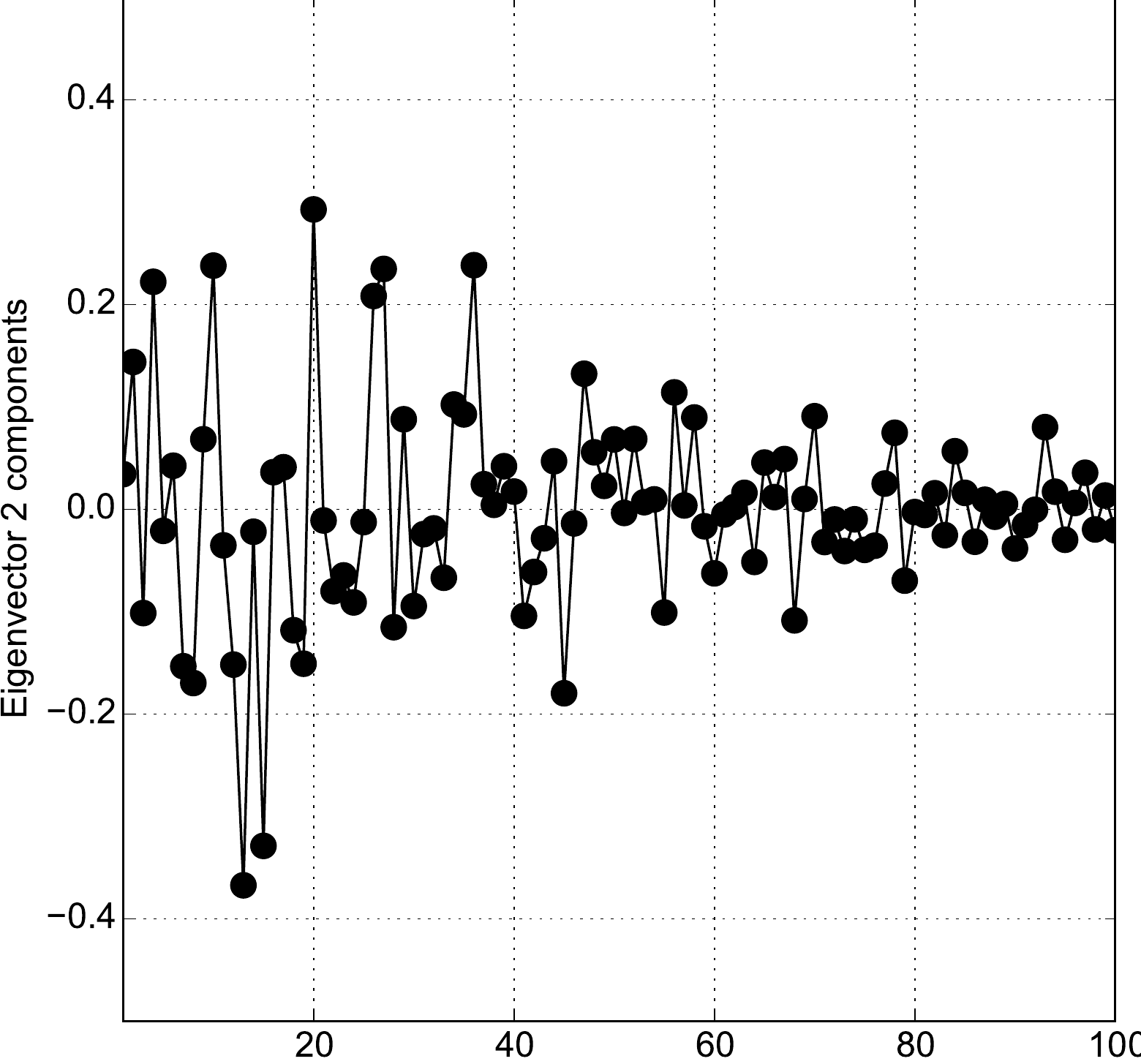}
}\\
\subfloat[Coefficient field with first eigenvector]{
\label{fig:32}
\includegraphics[width=0.45\linewidth]{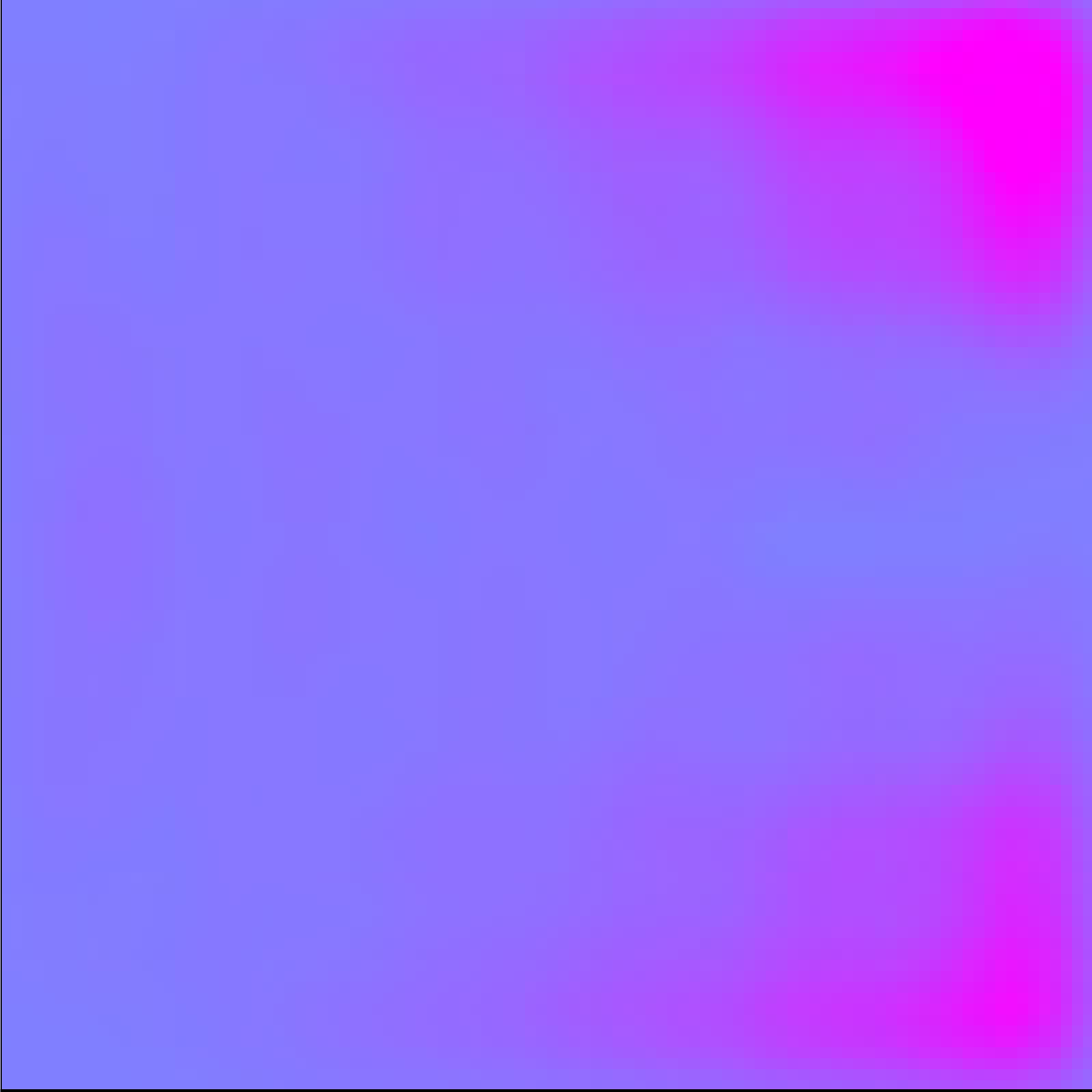}
}\;
\subfloat[Coefficient field with second eigenvector]{
\label{fig:33}
\includegraphics[width=0.45\linewidth]{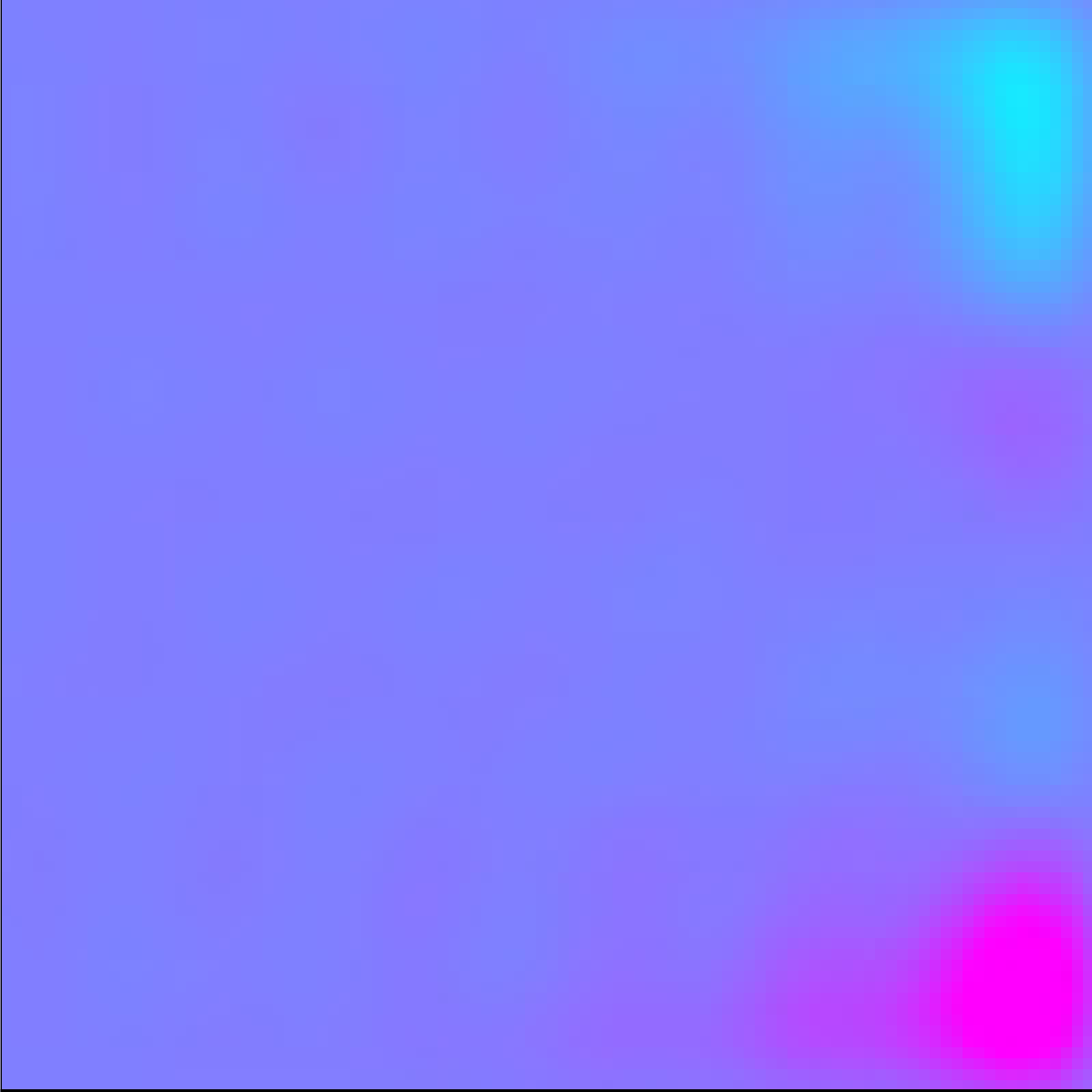}
}
\caption{The top two figures show the first and second eigenvector from $\hmC$ in \eqref{eq:mc}. The bottom two figures show the corresponding log-coefficient fields when these eigenvectors are plugged into $\vx$ in \eqref{eq:kl}. These fields can be treated as quantifying the spatial sensitivity of the misfit. (Colors are visible in the electronic version.)}
\label{fig:3}
\end{figure} 

Figures \ref{fig:30} and \ref{fig:31} show the first and second eigenvector components from $\hmC$ in \eqref{eq:mc}. Recall the interpretation of these eigenvectors: perturbing $\vx$ along some linear combination of the first two eigenvectors changes the misfit more, on average, than perturbing $\vx$ in a direction orthogonal to the span of the first two eigenvectors. These eigenvectors quantify global sensitivity of the misfit to the parameters. If we use the eigenvector components as $\vx$ in \eqref{eq:kl}, we get a spatially varying perturbation of the log-coefficients $\log(a)$ that changes the misfit the most, on average. The spatial perturbation corresponding to the first eigenvector is in Figure \ref{fig:32}, and the spatial perturbation corresponding to the second eigenvector is in Figure \ref{fig:33}. The sensitivity is larger near the boundary where the observations are taken.

\subsection{Applying the active subspace-accelerated MCMC}

\begin{figure}[ht]
\centering
\includegraphics[width=0.45\linewidth]{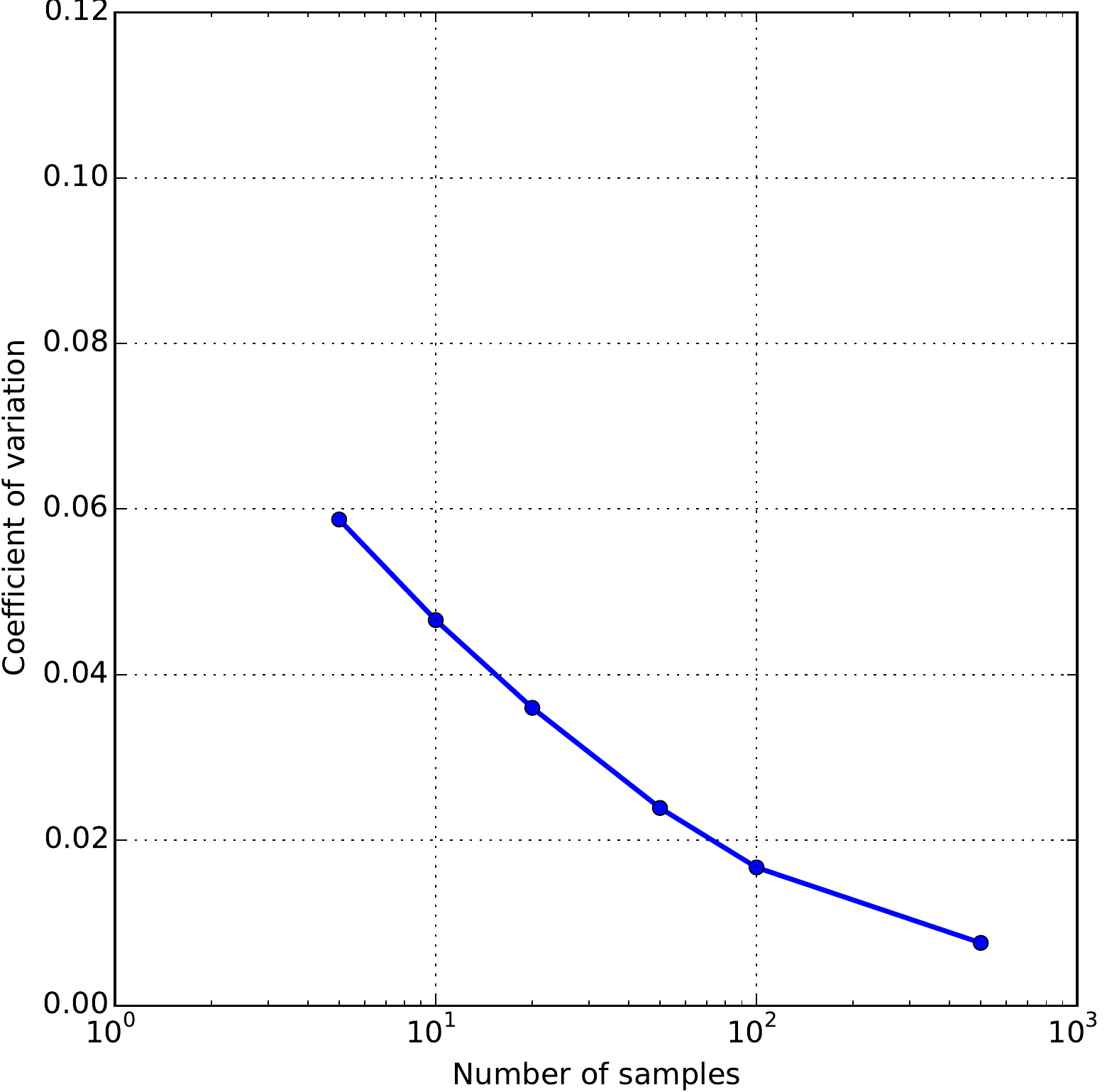}
\caption{Average coefficient of variation as a function of the number $M$ of samples for 100 Monte Carlo estimates of the conditional expectation $\hgeps(\hvy)$ from \eqref{eq:epsmccexp}. We use $M=10$ in Step 2 of Algorithm \ref{alg:mcmc}, which gives one-to-two digits of accuracy. (Colors are visible in the electronic version.)}
\label{fig:mcexp}
\end{figure}

%\begin{table}[htb]
%\begin{center}
%\begin{tabular}{c|c|c|c|c|c}
%Proposal dim. & Proposal var. & \# steps & \# forward models & %Accept rate & Eff.~sample size\\
%\hline
%2 & 0.1 & 50k & 500k & 80\% & 682 \\
%2 & 0.3 & 50k & 500k & 62\% & 4,030\\
%100 & 0.1 & 500k & 500k & 60\% & 766 
%\end{tabular}
%\caption{The first two rows show cases of Algorithm \ref{alg:mcmc} that exploits the active subspace. The last row is a standard Metropolis-Hastings MCMC on all 100 dimensions, i.e., the vanilla case. The first column is the dimension of the proposal density. The second column is the proposal variance. The third column is the number of steps in the MCMC, and the fourth column shows the number of forward model evaluations. The fifth column shows the acceptance rate, and the last column is the effective sample size computed from the first component of $\hvy$ samples in the first two rows and the 10th component of $\vx$ samples in the last row. The effect of the low-dimensional MCMC is apparent in the differing effective sample sizes.}
%\label{tab:rates}
%\end{center}
%\end{table}

\begin{table}[htb]
\begin{center}
\begin{tabular}{l|l|l|l}
 & Vanilla & AS2, 0.1 & AS2, 0.3 \\
\hline
Proposal dimension & 100 & 2 & 2 \\
Proposal variance & 0.1 & 0.3 & 0.1 \\
Number of steps & 500k & 50k & 50k \\
Number of forward models & 500k & 500k & 500k \\
Acceptance rate & 60\% & 80\% & 62\% \\
Min.~eff.~sample size $\hvy$ & N/A & 198 & 1053\\
Min.~eff.~sample size $\vx$ & 604 & 81986 & 47281
\end{tabular}
\caption{The first column shows statistics for standard Metropolis-Hastings MCMC on all 100 dimensions, i.e., the vanilla case. The second two rows show cases of Algorithm \ref{alg:mcmc} that exploits the active subspace. The first row is the dimension of the proposal density. The second row is the proposal variance. The third row is the number of steps in the MCMC, and the fourth row shows the number of forward model evaluations. The fifth row shows the acceptance rate. The sixth row shows the smallest effective sample size over the two active variables for chains operating in the active subspace. The last row shows the smallest effective sample size over 100 components of the chain in the original 100-dimensional space.}
\label{tab:rates}
\end{center}
\end{table}

\begin{figure}[ht]
\centering
\subfloat[Autocorrelation $x_{10}$]{
\label{fig:40}
\includegraphics[width=0.45\linewidth]{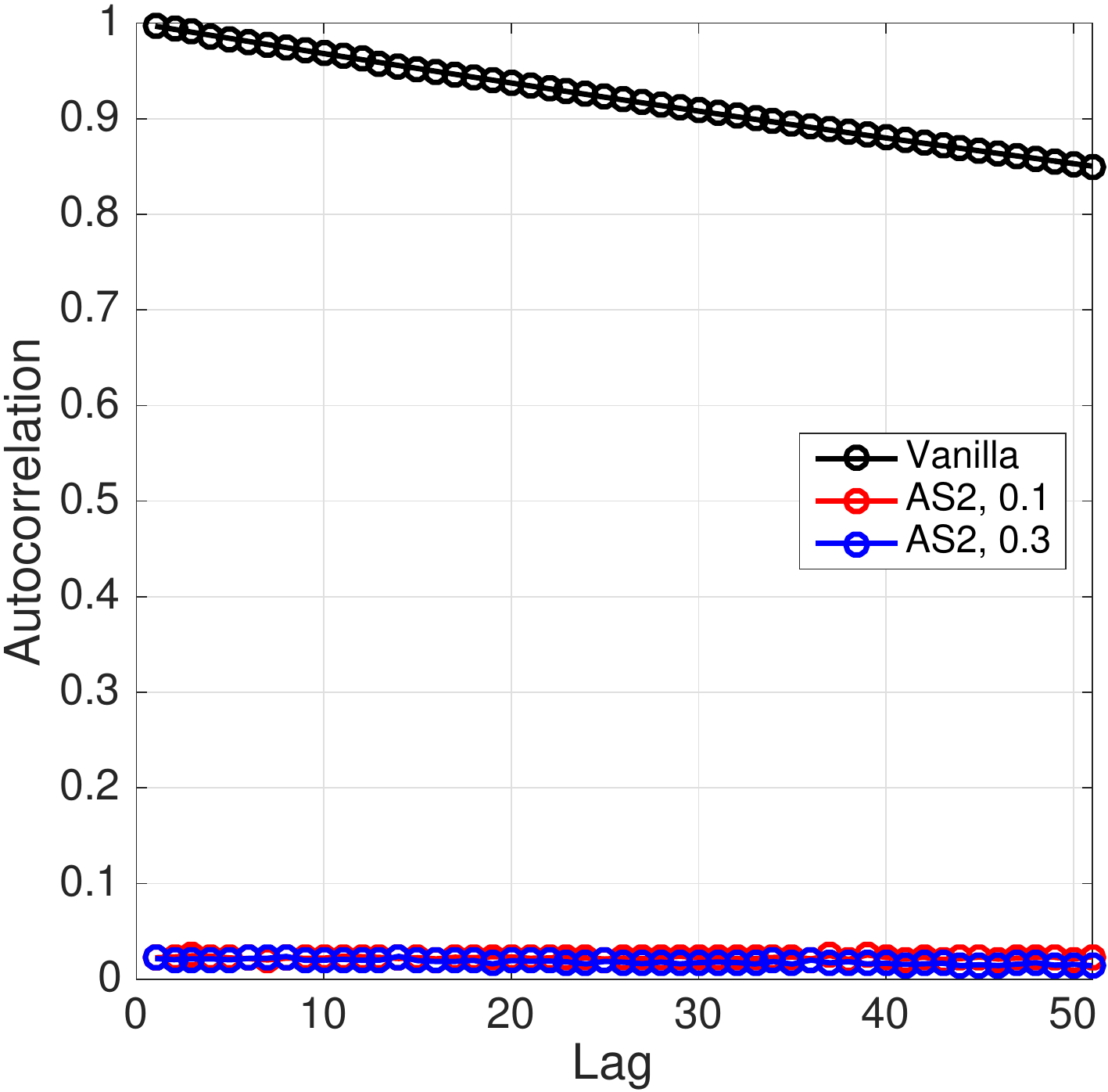}
}\;
\subfloat[Trace plot $x_{10}$, Vanilla]{
\label{fig:41}
\includegraphics[width=0.45\linewidth]{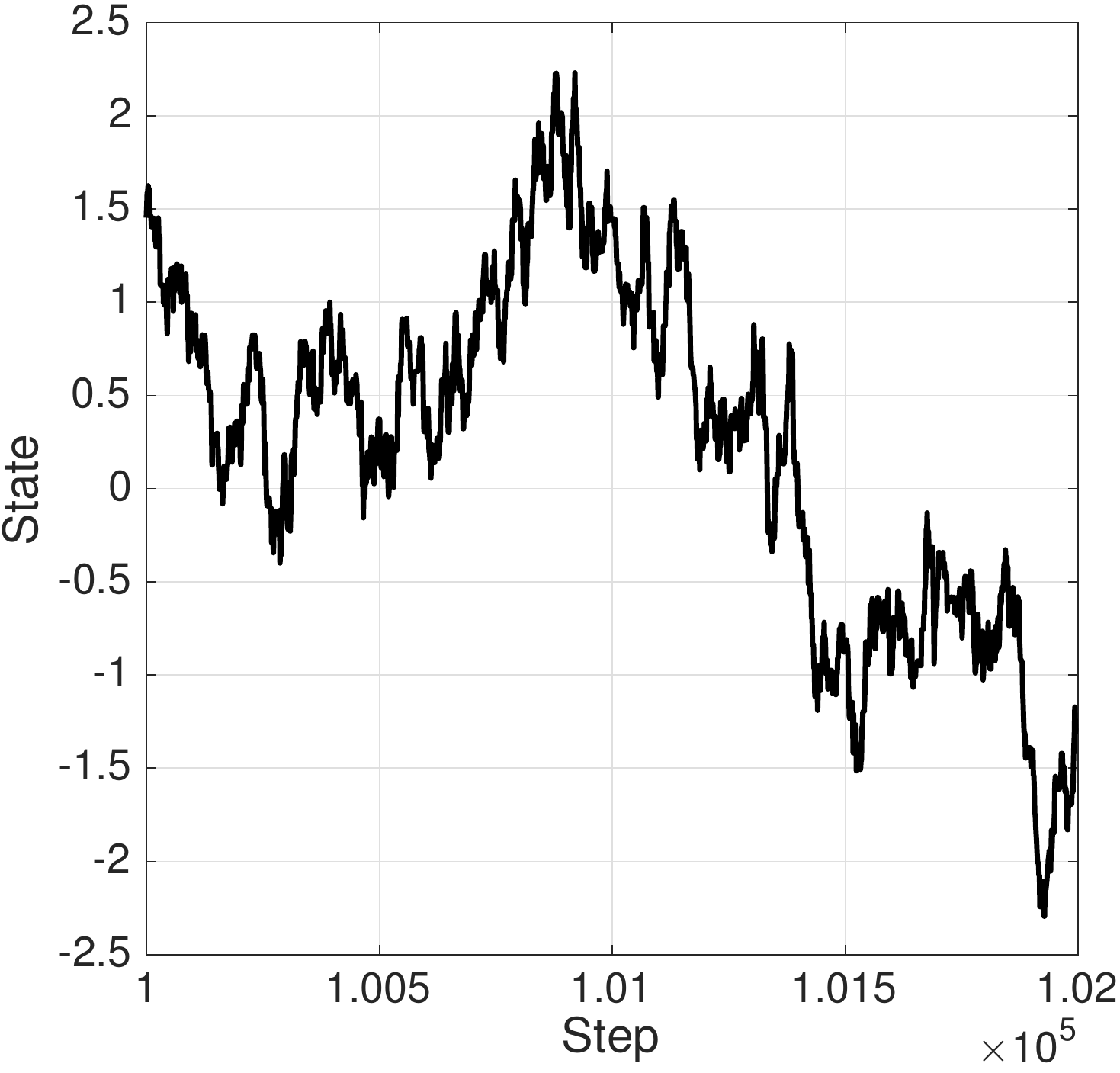}
}\\
\subfloat[Trace plot $x_{10}$, AS2, 0.1]{
\label{fig:42}
\includegraphics[width=0.45\linewidth]{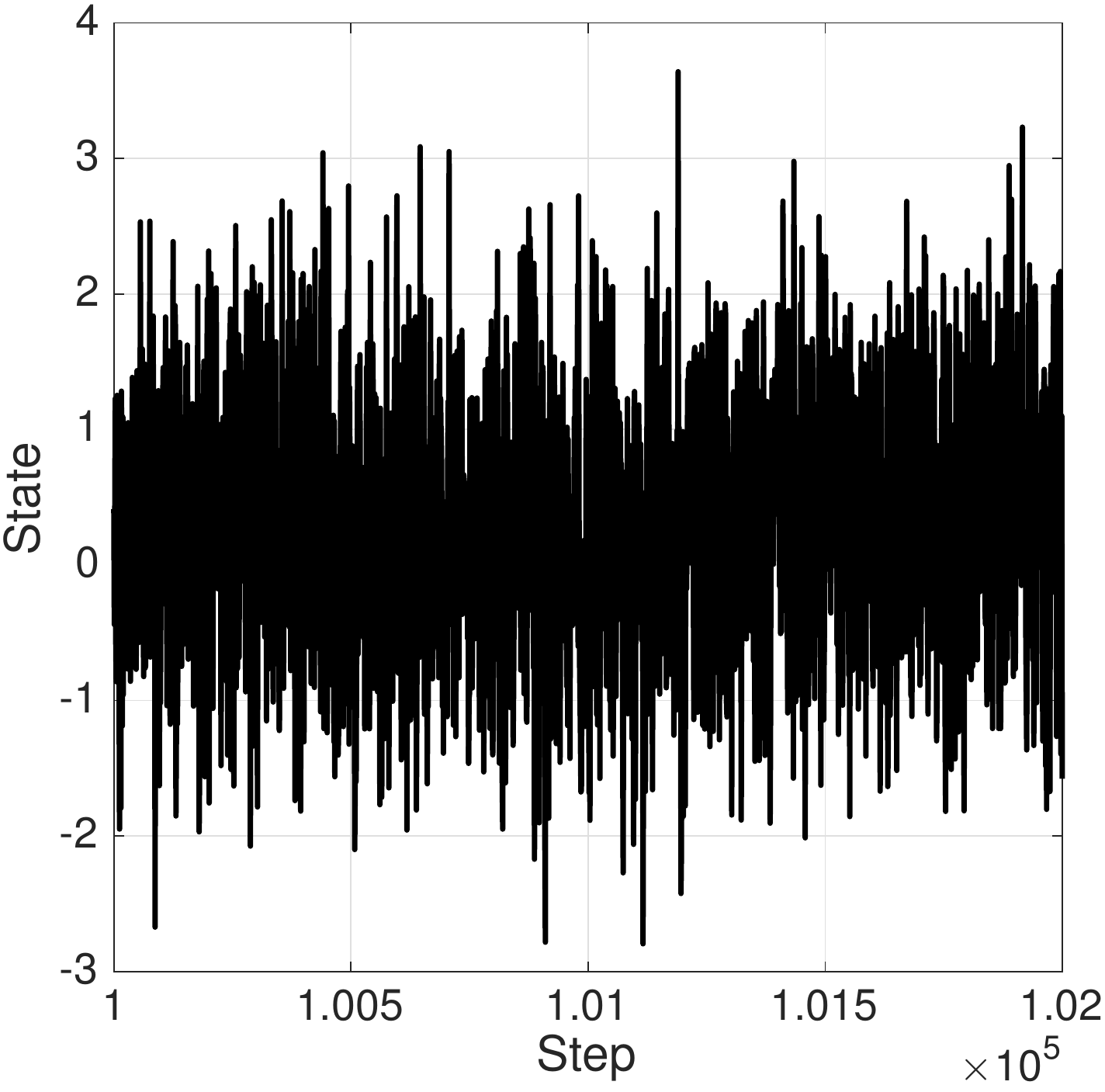}
}\;
\subfloat[Trace plot $x_{10}$, AS2, 0.3]{
\label{fig:43}
\includegraphics[width=0.45\linewidth]{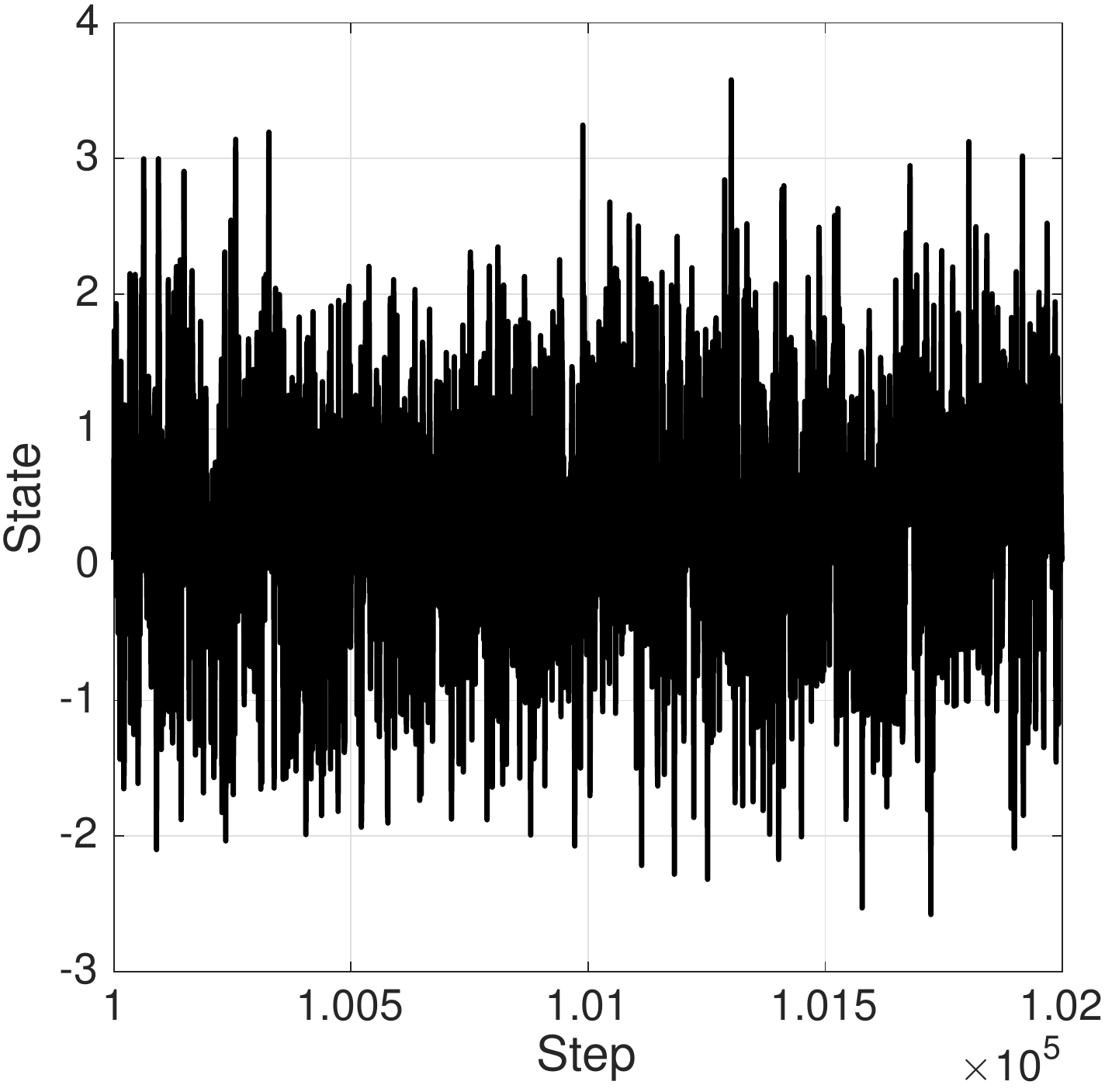}
}
\caption{The top left figure shows the autocorrelation function for the 10th component of $\vx$ in three cases: (i) the Metropolis-Hastings MCMC in 100 dimensions (Vanilla), (ii) Algorithm \ref{alg:mcmc} with a two-dimensional subspace and proposal variance 0.1 (AS2, 0.1), and (iii) Algorithm \ref{alg:mcmc} with a two-dimensional subspace and proposal variance 0.3 (AS2, 0.3). The remaining figures show trace plots of $x_{10}$ for all three cases. The active subspace-accelerated MCMC produces Markov chain iterates with much smaller correlation. (Colors are visible in the electronic version.)}
\label{fig:4}
\end{figure}

With the eigenvectors $\hmW$ from \eqref{eq:mc} and the choice of $n=2$ (i.e., a two-dimensional active subspace)---justified by the eigenvalues (Figure \ref{fig:20}) and subspace error estimates (Figure \ref{fig:21})---we can apply the active subspace-accelerated MCMC in Algorithm \ref{alg:mcmc}. To choose the number $M$ of Monte Carlo samples in Step 2, we perform the following experiment.
\begin{enumerate}
\item Draw 100 $\vx$'s independently at random according to the Gaussian prior.
\item For each $\vx$, let $\hvy=\hmW_1^T\vx$, and compute the coefficient of variation for the Monte Carlo estimates,
\begin{equation}
\frac{
\left(\frac{1}{M-1}\sum_{i=1}^M (f(\hmW_1\hvy + \hmW_2\hvz_i) - \hgeps(\hvy))^2\right)^{1/2}
}{
\sqrt{M}\,\hgeps(\hvy)
},
\end{equation}
where $\hgeps(\hvy)$ is from \eqref{eq:epsmccexp}, for $M=$1, 5, 10, 20, 50, 100, and 500. 
\item For each $M$, average the coefficients of variation over all $\vx$'s. 
\end{enumerate}
Figure \ref{fig:mcexp} shows the average coefficients of variation as a function of $M$. We choose $M=10$ for Step 2 in Algorithm \ref{alg:mcmc}, which is sufficient for one-to-two digits of accuracy from the Monte Carlo estimates $\hgeps(\hvy)$.

We use a standard Gaussian proposal density in Step 1 of Algorithm \ref{alg:mcmc}. We compare results from two variants of the active subspace-accelerated MCMC: (i) a two-dimensional active subspace and a proposal variance of 0.1, and (ii) a two-dimensional active subspace with a proposal variance of 0.3. Each case ran 50k steps of the Markov chain, which used 500k forward model evaluations and 120 CPU hours. We discard 10k steps as a burn-in. For each sample $\hvy_k$ from the MCMC, we use $P=10$ independent samples of the 98-dimensional inactive variables $\hvz$ drawn according to their standard Gaussian prior as in \eqref{eq:reconstruct} to construct a chain on the 100-dimensional parameter space.

We compare the results to a standard Metropolis-Hastings MCMC in all $m=100$ dimensions with a Gaussian proposal density with variance $0.1$. We refer to this chain as the \emph{vanilla} case. The 100-dimensional Markov chain took 500k steps, which corresponds to 500k forward model evaluations, and we discarded 100k steps as a burn-in. We do not compare the active subspace-accelerated method to other MCMC variants, because any of those variants can exploit the active subspace the same way Algorithm \ref{alg:mcmc} does. 

Table \ref{tab:rates} displays several characteristics of the Markov chains for the three cases. The first two rows show the proposal density dimension and variance. The two-dimensional chains on the active variables $\hvy$ use 50k steps, while the 100-dimensional chain uses 500k steps. Each evaluation of the approximate likelihood $\hgeps(\hvy)$ in Step 2 of Algorithm \ref{alg:mcmc} uses $M=10$ forward model evaluations, so the total number of forward model evaluations is the same (500k) across all chains. The fifth row shows the acceptance rates for the chains. Note that the two-dimensional chain can use a larger proposal variance while maintaining a comparable acceptance rate to the vanilla case. We compute the effective sample size of the chain's components as
\begin{equation}
\frac{N_{\text{steps}}}{1 + 2 \sum_{k=1}^{2000} \rho_k},
\end{equation}
where $N_{\text{steps}}$ is the number of steps in the chain, and $\rho_k$ is the autocorrelation with lag $k$. The sixth row of Table \ref{tab:rates} shows the minimum effective sample size over the two components of the two-dimensional chain in the coordinates of the active subspace. The last row shows the minimum effective sample size over all 100 components of the chains on the full 100-dimensional parameter space. The much larger effective sample sizes for the active subspace-accelerated chains is due to the independent sampling according to the prior on the inactive variables. 

Figure \ref{fig:40} shows the autocorrelation function for the 10th component of $\vx$ for all three chains on the 100-dimensional parameter space: the standard MCMC in 100 dimensions (vanilla), the active subspace-accelerated MCMC in 2 dimensions with proposal variance 0.1 (AS2, 0.1), and the active subspace-accelerated MCMC in 2 dimensions with proposal variance 0.3 (AS2, 0.3). The other components of $\vx$ had similar autocorrelation functions. The slow decay in the vanilla case is due to the Markov chain operating in all 100 dimensions, while the active subspace accelerated cases run the Markov chain in only two dimensions and draw the remaining components independently according to the prior. Thus, the iterates in $\vx$ appear uncorrelated. The remaining subplots in Figure \ref{fig:4} show trace plots of $\vx$'s 10th component for all three cases. The active subspace-accelerated method mixes much better in the space of $\vx$ due to the independent sampling of $\hvz$; this faster mixing justifies the term \emph{accelerated} in the title.

\begin{figure}[ht]
\centering
\subfloat[Mean, AS2, 0.1]{
\label{fig:50}
\includegraphics[width=0.45\linewidth]{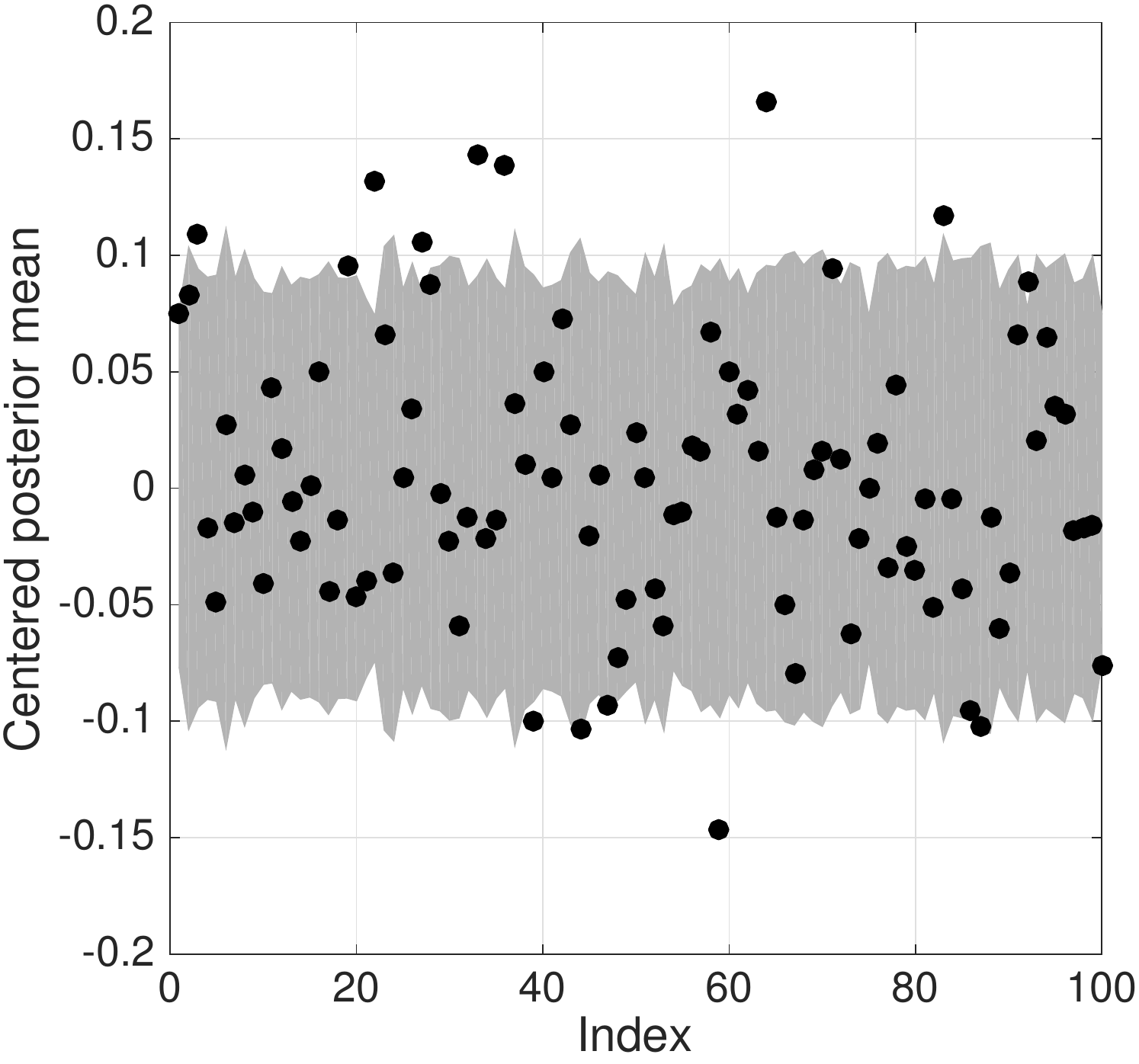}
}\;
\subfloat[Mean, AS2, 0.3]{
\label{fig:51}
\includegraphics[width=0.45\linewidth]{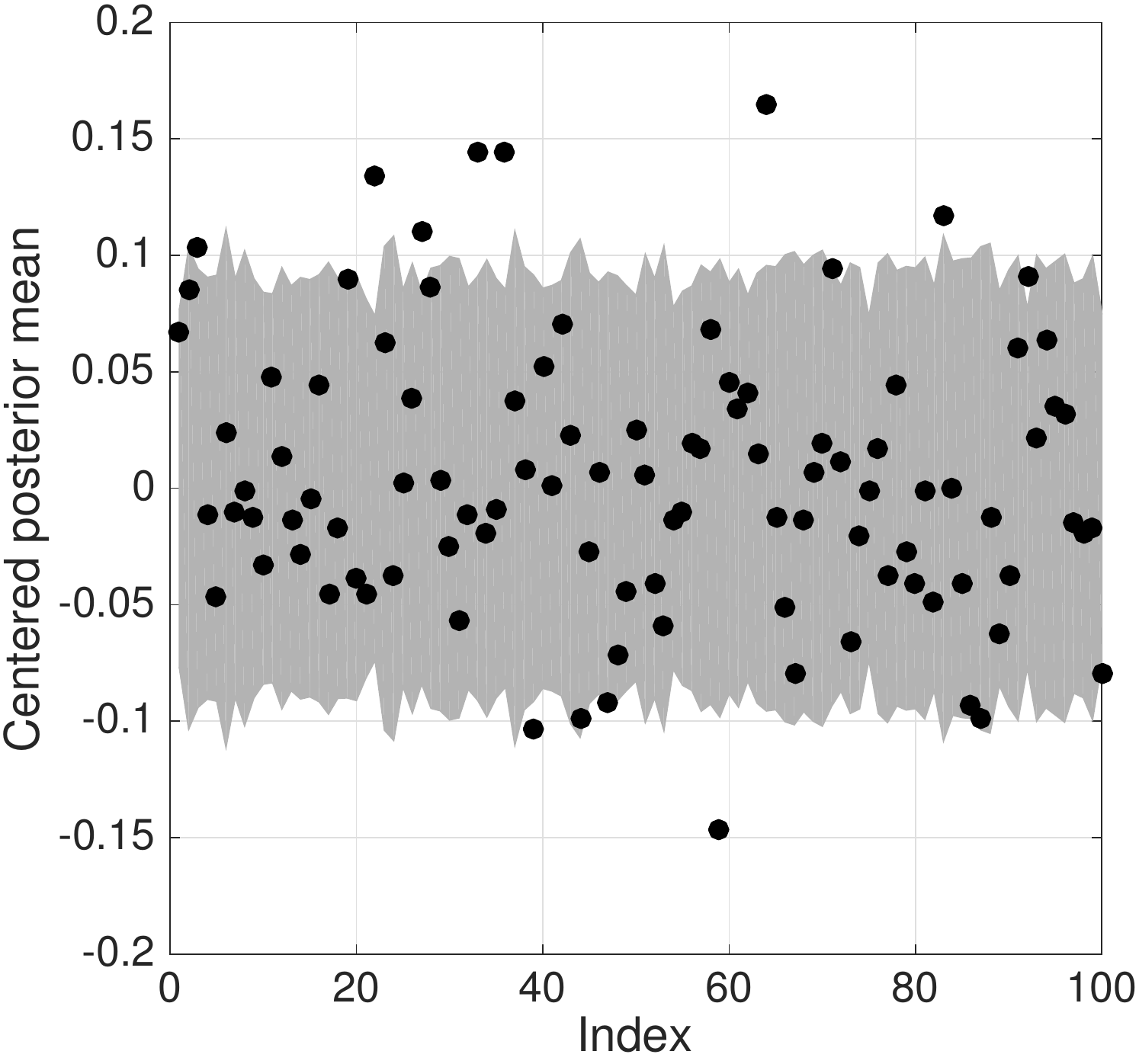}
}\\
\subfloat[Variance, AS2, 0.1]{
\label{fig:52}
\includegraphics[width=0.45\linewidth]{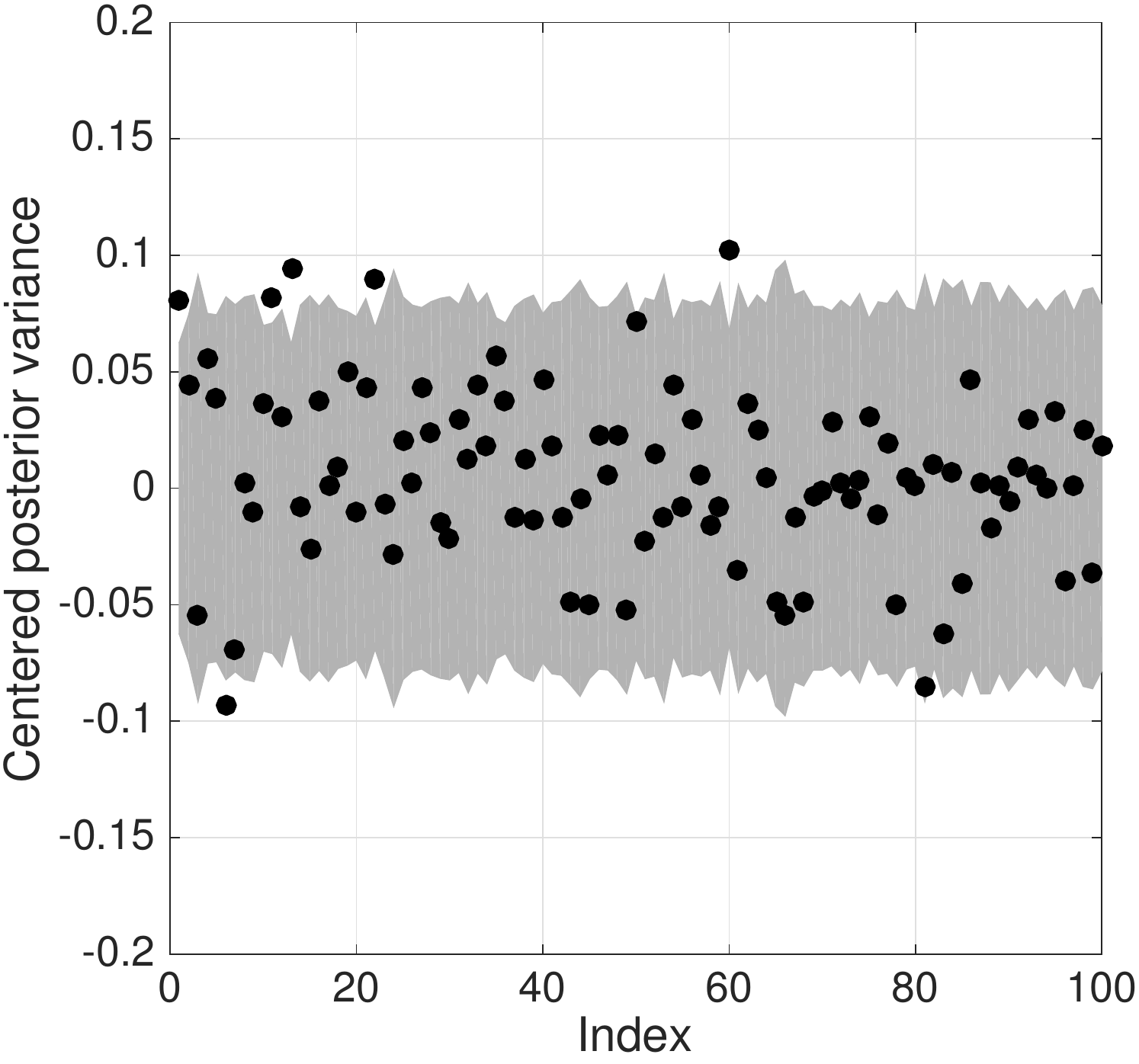}
}\;
\subfloat[Variance, AS2, 0.3]{
\label{fig:53}
\includegraphics[width=0.45\linewidth]{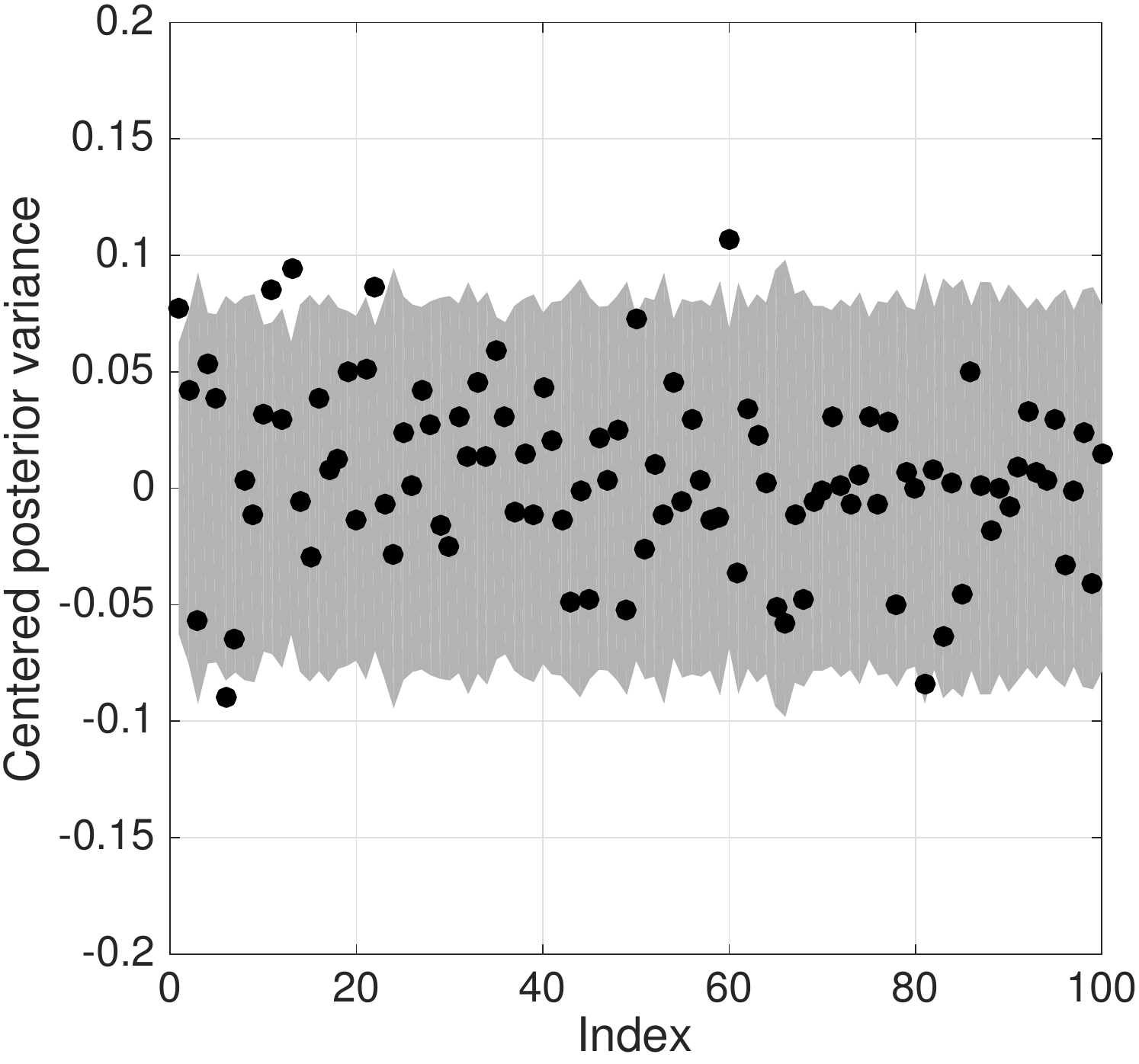}
}
\caption{These figures show the sample moments (black dots) of the active subspace-accelerated chains shifted by the sample moments from the vanilla chain; the top row shows sample means and the bottom row shows sample variances. The gray regions are the asymptotically valid 99\% confidence intervals, shifted by the sample moments, computed with consistent batch means as in~\cite{Flegal2008}. The left column shows the accelerated chain with proposal variance 0.1, and the right column shows the chain with proposal variance 0.3.}
\label{fig:5}
\end{figure}

Figure \ref{fig:4} and Table \ref{tab:rates} suggest that the active subspace-accelerated MCMC mixes much faster than the vanilla MCMC. But are its iterates producing correlated samples from a density close to the true posterior? MCMC convergence metrics can be difficult to interpret for a single chain, so comparing results from different chains is especially challenging. The reader should treat the following results as qualitative, since any quality metrics for the vanilla MCMC are computed from a 100-dimensional chain with high autocorrelation and low effective sample size; see Table \ref{tab:rates}. We perform the following test to check the mean and variance of the iterates produced by the active subspace-accelerated MCMC. We first compute asymptotically valid 99\% confidence intervals on the posterior mean and variance from the vanilla MCMC using \emph{consistent batch means} as in section 3.1 of Flegal, et al.~\cite{Flegal2008} with parameter $\theta=2/3$. Let $\hat{\mu}\in\mathbb{R}^{100}$ be the sample posterior mean, and let $\hat{\sigma}^2\in\mathbb{R}^{100}$ be the sample posterior variance---both computed from the vanilla MCMC. Denote the confidence intervals 
\begin{equation}
\hat{\mu}_{\ell} \leq \hat{\mu} \leq \hat{\mu}_u,\qquad
\hat{\sigma}_{\ell}^2 \leq \hat{\sigma}^2 \leq \hat{\sigma}_u^2,
\end{equation}
where the inequalities are interpreted component-wise. Denote the sample mean and variance from the active subspace-accelerated chains as $\hat{\mu}_{\text{as}}$ and $\hat{\sigma}^2_{\text{as}}$, respectively. Figure \ref{fig:5} compares the shifted moments, $\hat{\mu}_{\text{as}}-\hat{\mu}$ and $\hat{\sigma}^2_{\text{as}}-\hat{\sigma}^2$, to the the shifted confidence intervals,
\begin{equation}
[\hat{\mu}_{\ell}-\hat{\mu},\,\hat{\mu}_u-\hat{\mu}], \qquad
[\hat{\sigma}^2_{\ell}-\hat{\sigma}^2,\,\hat{\sigma}^2_u-\hat{\sigma}^2].
\end{equation}
Shifting by the sample mean and variance allows easier visual comparison. Figures \ref{fig:50} and \ref{fig:52} show (i) the shifted mean and variance (black dots), respectively, for the active subspace-accelerated chain with proposal variance 0.1 and (ii) the shifted confidence intervals (gray region). Figures \ref{fig:51} and \ref{fig:53} show the same shifted moments for the accelerated chain with proposal variance 0.3. We also computed the consistent batch means-based confidence intervals for the accelerated chain's moments, but the intervals were very small (within the marker size), so we do not include them in the plots. The small confidence intervals for the accelerated chain's moments are consistent with the observed rapid mixing. The bulk of the black dots fall within the confidence region, which indicates general agreement in the sample moments. The dots that fall outside the confidence region may indicate bias introduced by the active subspace-based dimension reduction. 

\begin{figure}[ht]
\centering
\subfloat{
\includegraphics[width=0.3\linewidth]{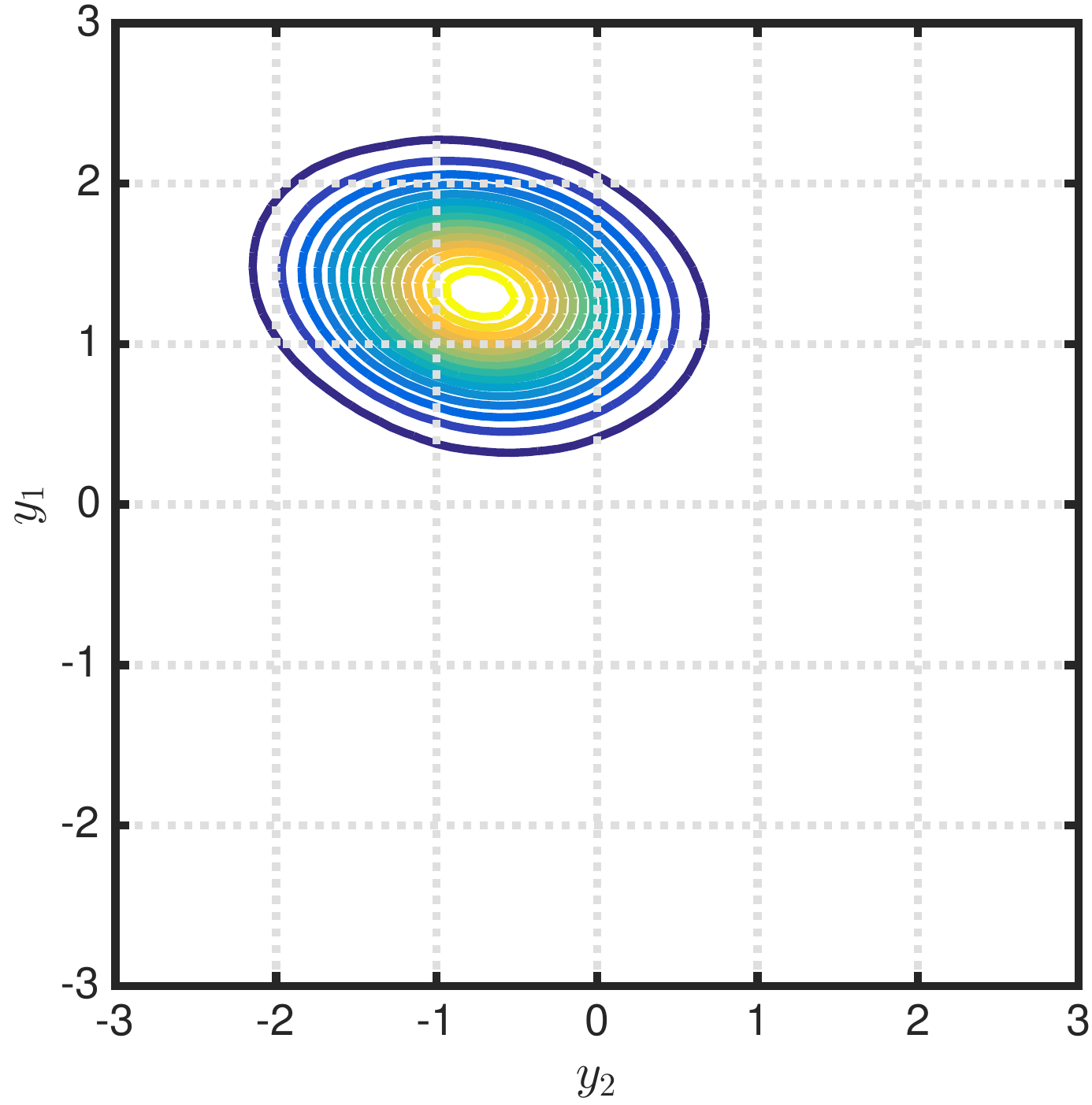}
}\;
\subfloat{
\includegraphics[width=0.3\linewidth]{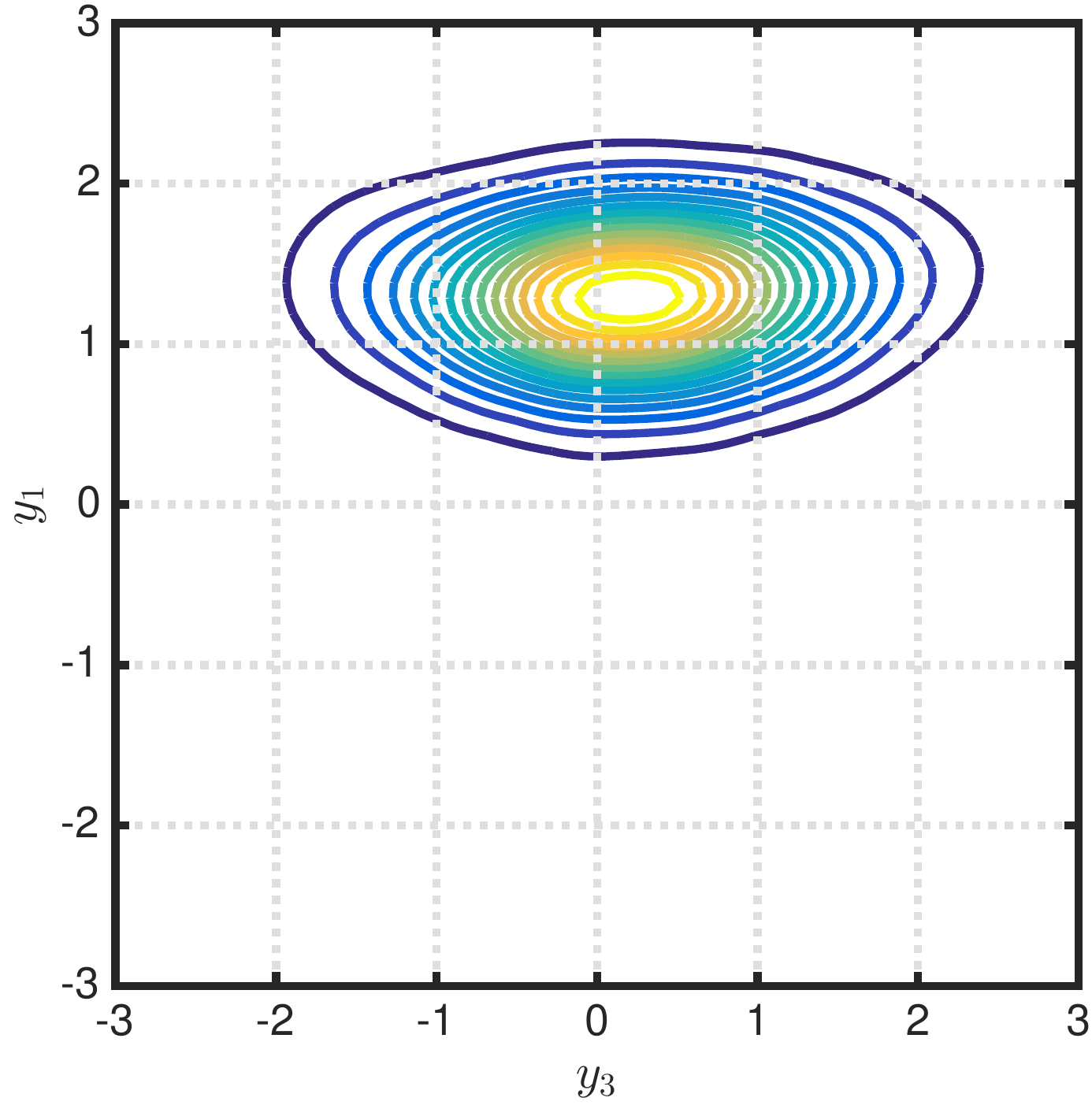}
}\;
\subfloat{
\includegraphics[width=0.3\linewidth]{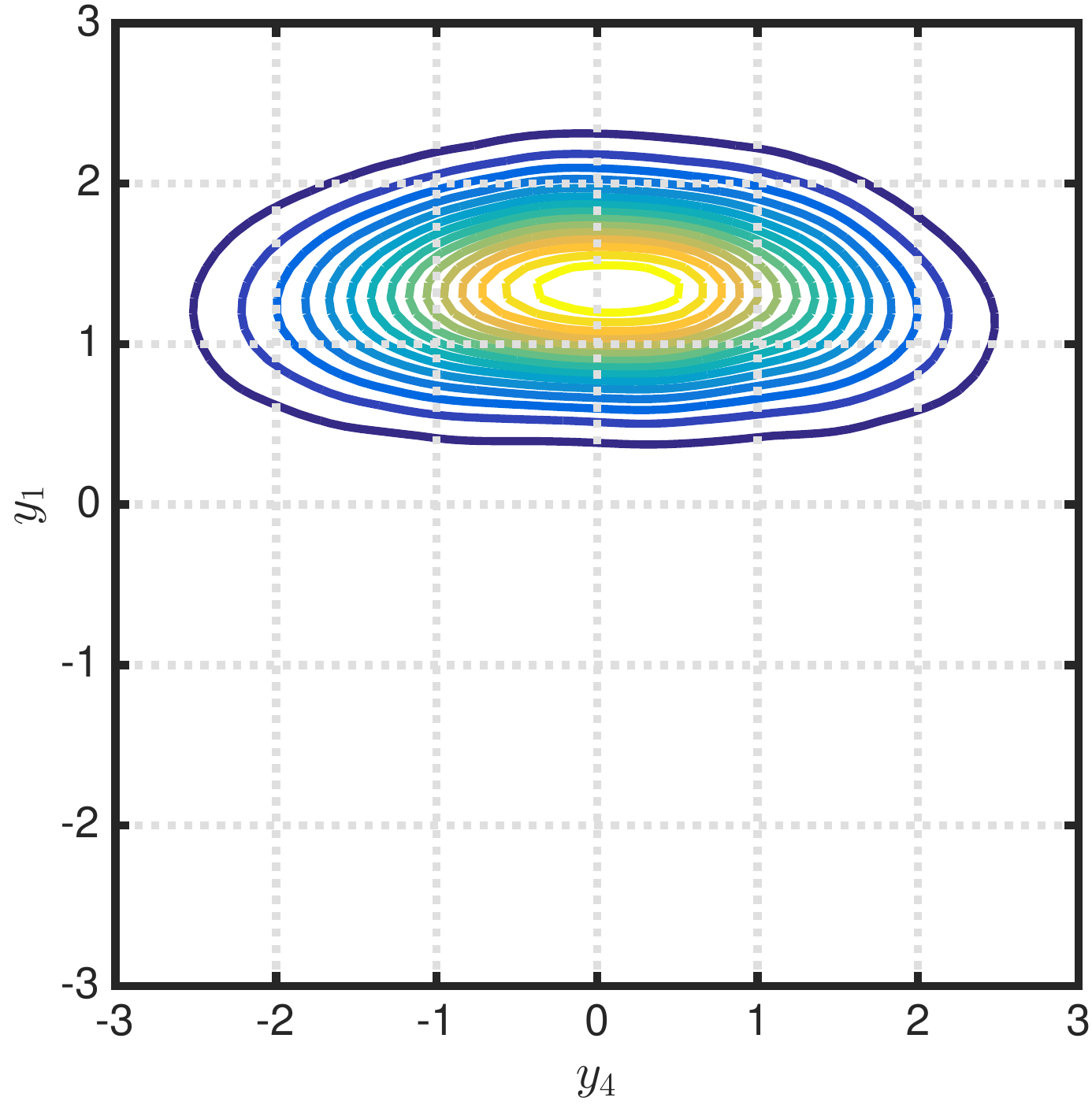}
}\\
\subfloat{
\includegraphics[width=0.3\linewidth]{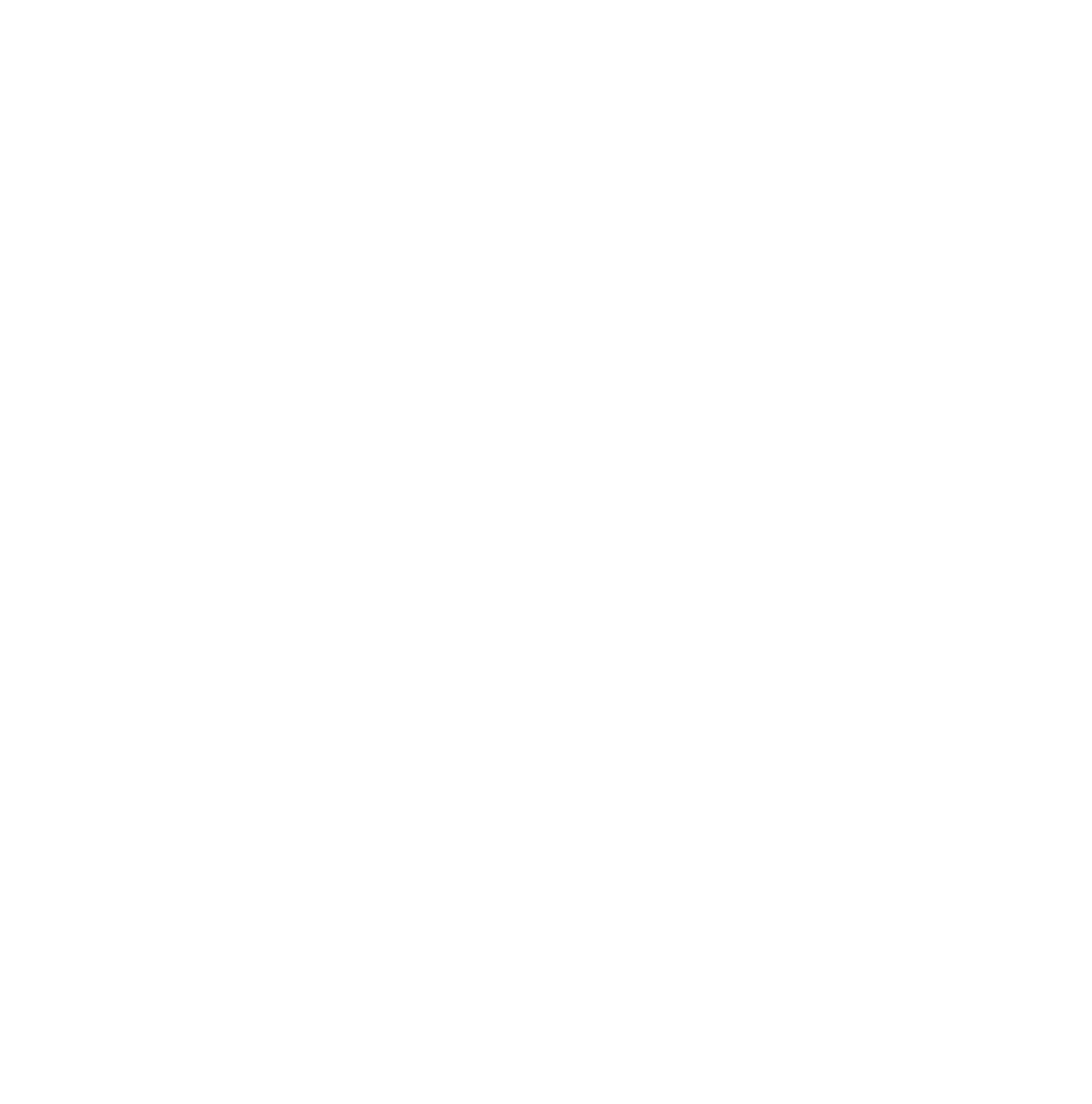}
}\:
\subfloat{
\includegraphics[width=0.3\linewidth]{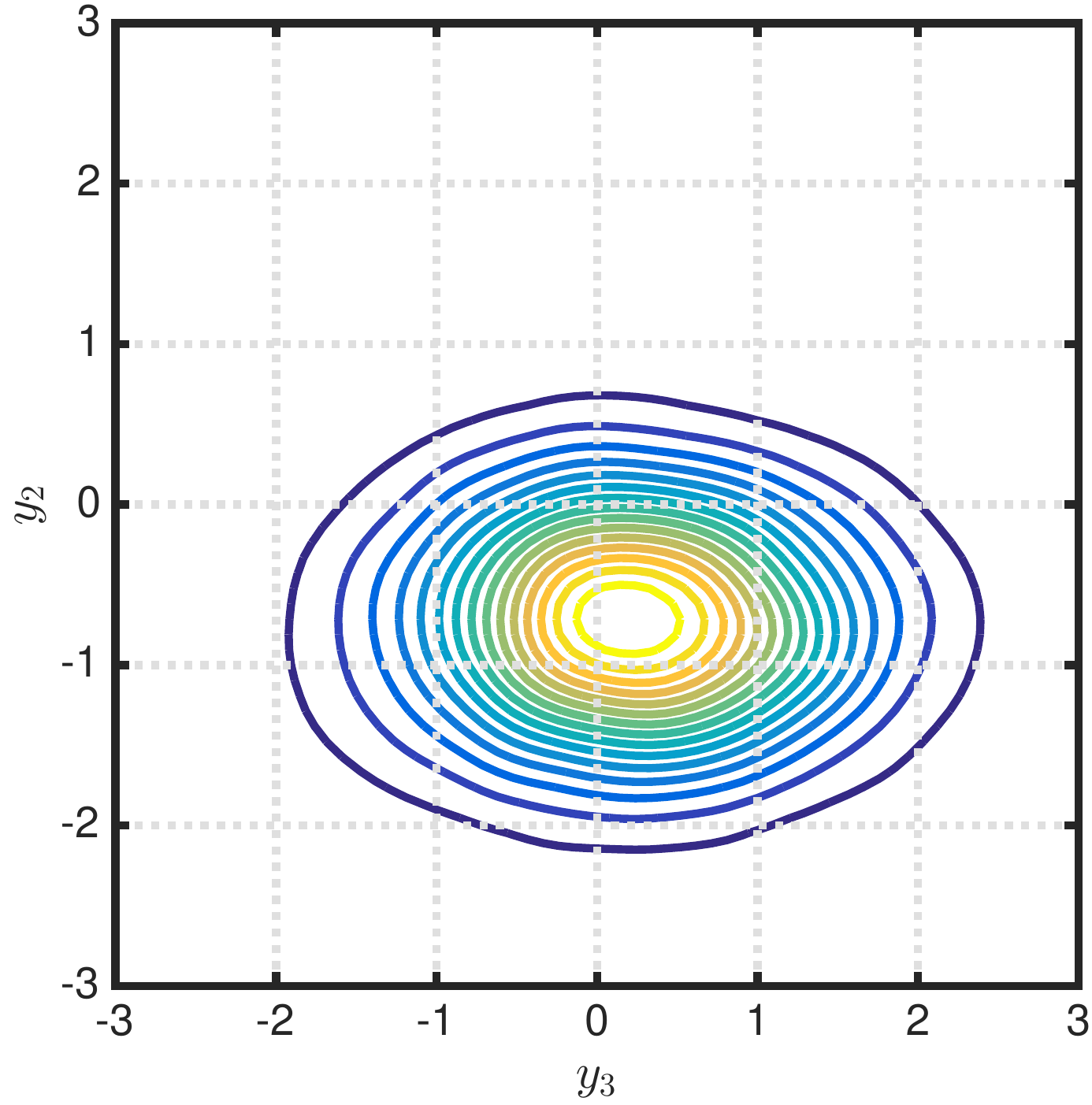}
}\;
\subfloat{
\includegraphics[width=0.3\linewidth]{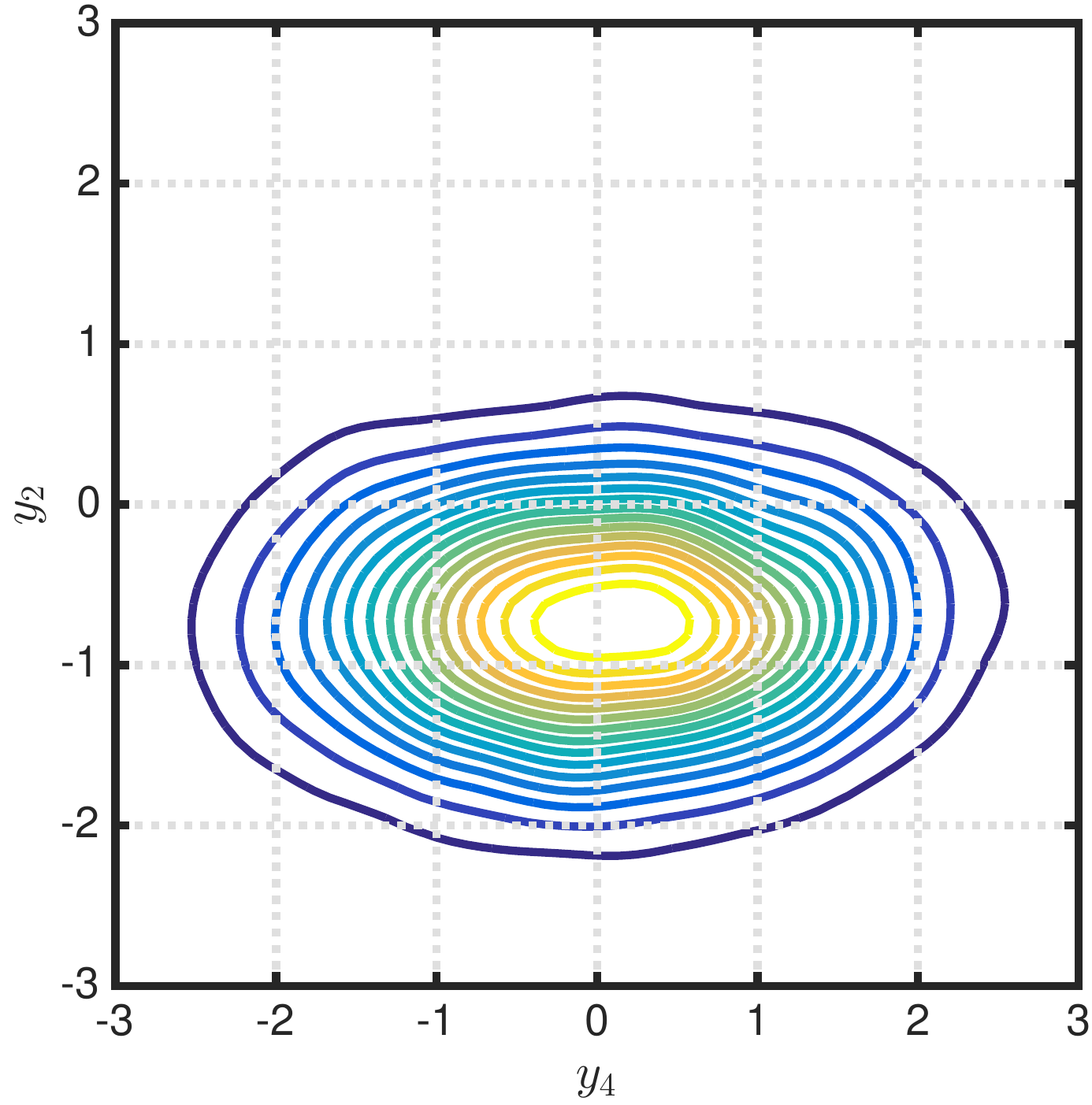}
}\\
\subfloat{
\includegraphics[width=0.3\linewidth]{figs/blank}
}\:
\subfloat{
\includegraphics[width=0.3\linewidth]{figs/blank}
}\;
\subfloat{
\includegraphics[width=0.3\linewidth]{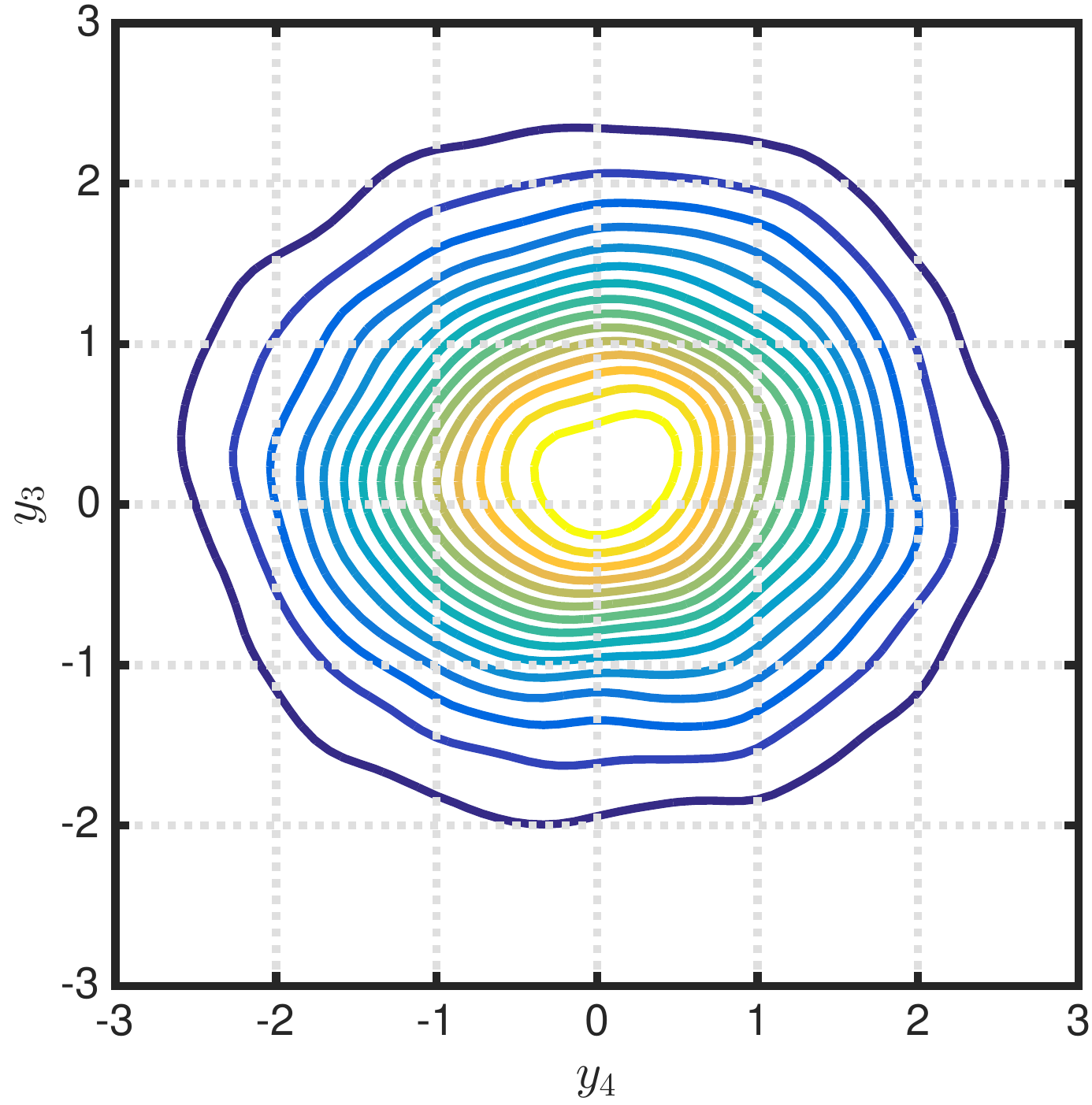}
}
\caption{Bivariate posterior marginals for the first four components from the transformed variables $\hvy_k=\hmW^T\vx_k$, where the set $\{\vx_k\}$ contains samples from the vanilla MCMC in 100 dimensions with proposal variance 0.1. (Colors are visible in the electronic version.)}
\label{fig:6}
\end{figure} 

We perform one final qualitative check for the two-dimensional structure using the iterates from the 100-dimensional vanilla MCMC. With the eigenvectors $\hmW$ from \eqref{eq:mc}, we can transform all iterates $\vx_k$ from the vanilla chain to the space of active variables, $\hvy_k = \hmW^T\vx_k$. In this transformation, $\hvy_k$ has 100 components, but they are ordered like the eigenvalues of $\hmC$. In other words, the first two components of $\hvy_k$ are comparable to the active variables. Figure \ref{fig:6} shows contours of the kernel density estimates of the bivariate marginals using the first four components from the samples $\hvy_k$; the axis labels indicate which components are shown. Note that the active subspace was not employed to generate these samples; they are from the vanilla MCMC. Nevertheless, the density estimate for first two active variable components, $\hat{y}_1$ and $\hat{y}_2$, departs from the standard Gaussian prior. In other words, the data updates this density from the prior. In contrast, the density estimate for components $\hat{y}_3$ and $\hat{y}_4$ more closely resembles the standard Gaussian prior; this pattern continues for the remaining $\hvy$ components. Post-processing the vanilla MCMC iterates with the eigenvectors $\hmW$ validates the two-dimensional structure suggested by the eigenvalues in Figure \ref{fig:20}, providing further confirmation that the data informs only a two-dimensional subspace in the 100-dimensional parameter space. 

\section{Summary and conclusions}
\label{sec:conclusion}

We have shown how to use the active subspace---derived from the scalar-valued data misfit and defined with respect to the prior---to accelerate MCMC for high-dimensional Bayesian inverse problems with nonlinear parameter-to-observable maps. One can estimate the active subspace as a preprocessing step before beginning the MCMC. Since the active subspace is defined with respect to the prior, its components can be estimated in parallel. If the low-dimensional structure is apparent, then the MCMC can be run on only the active variables, which correspond to the subspace informed by the data. The inactive variables are constrained by the prior; they can be sampled independently as a post-processing step to reconstruct a chain on the original parameter space. We have bounded the error in the Hellinger distance between the true posterior and its approximation with the active subspace. The bound is in terms of the eigenvalues and the error the numerically estimated active subspace basis. We demonstrated the approach on (i) a two-dimensional example with a quadratic forward model and one-dimensional active subspace and (ii) a 100-dimensional example with a PDE-based forward model and a two-dimensional active subspace. 

The proposed approach has several limitations that practitioners should consider when evaluating its appropriateness for their own Bayesian inverse problems. The eigenvalues of $\mC$ derived from the misfit's gradient provide evidence of exploitable low-dimensional structure. However, these eigenvalues may be misleading for a suitably irregular parameter-to-observable map. The eigenvectors associated with small eigenvalues are reasonable directions to ignore in the parameter space \emph{on average}. A parameter-to-observable map with large, localized variability in a high-dimensional parameter space may yield an estimated active subspace that misses important variation in the misfit. (Such irregular models cause problems for most inversion methods.) Moreover, since the active subspace is derived from the misfit's gradient, a misfit with large gradients but small variation (such as an oscillating function) may have an active subspace that favors directions of oscillation---even when other directions are more apprpropriate for exploring the range of misfit values. We have not encountered this scenario in practice, but it is possible to construct such functions; see~\cite[Chapter 1]{asm2015}. For a linear forward model, when the given data is in the tails of the forward model's likely outputs (i.e., propagated according to the prior), the active subspace may be unduly influenced by the data vector. This insight is not directly relevant for MCMC, since MCMC is not appropriate for linear forward models. However, it suggests that deep understanding of how the data affect the active subspace may be needed for some nonlinear forward models; we intend to pursue such analysis in future studies. 

The practitioner should also be aware of scaling limitations of the proposed approach as the problem's dimensions increase. Without additional acceleration (e.g., surrogate forward models or alternative structure-exploiting techniques), random walk MCMC is practically limited to inference on a handful of parameters in expensive nonlinear simulation models. If the active subspace for a given problem has dimension greater than 5-to-10, then the practitioner should consider a more sophisticated approach than random walk MCMC for efficient inference. For PDE-based inverse problems, the holy grail is to infer a spatially varying parameter field, where the number of parameters is the number of spatial discretization cells---which may be in the millions for modern computational models. In the example from section \ref{sec:pdebayes}, we first assumed we could reduce the number of parameters from 10000 to 100 using known correlation structure in the parameters. Without this reduction, we may have needed thousands of realizations of the 10000-component gradient vector to estimate the eigenpairs of $\mC$. Such computations are large enough to stress workstation-sized computers. A more spatially refined computation---or one that originates from a PDE in three spatial dimensions---would cause more stress. Moreover, there is no guarantee that the misfit's active subspace dimension would be small enough to permit efficient MCMC. Therefore, for practical PDE-based inverse problems, we expect that estimating and exploiting an active subspace in the misfit would be one tool in a comprehensive toolbox for Bayesian inference.

There are other variants of active subspace-accelerated MCMC that are worth studying that we did not explore in this paper. For instance, one could define the active subspace using the posterior instead of the prior as the integration measure; this is similar to the likelihood-informed subspace~\cite{Cui2014}. Using the posterior would allow the proper interpretation of the subspace components as conditional random variables, but it would be more computationally expensive. Alternatively, one could use the joint density of the data and parameters in place of the prior, which would produce a data independent subspace. Another idea is to use the likelihood directly as the differentiable, scalar-valued function, instead of its negative log (i.e., the misfit). Lastly, there may be a way to combine subspaces from the state covariance, the matrix $\mC$ from \eqref{eq:C}, and the average Hessian as in~\cite{Cui2014} to produce a more robust dimension reduction. A thorough quantitative comparison to LIS is beyond the scope of the current manuscript, but such a comparison would be worthwhile.

\section*{Acknowledgments}
We thank Luis Tenorio and Aaron Porter at Colorado School of Mines and Youssef Marzouk and Tiangang Cui at MIT for their helpful comments. We also acknowledge the support of the J.~Tinsley Oden Faculty Fellowship Research Program for the first author's 2015 summer visit to the Institute for Computational Engineering and Sciences at UT Austin. This material is based upon work supported by (i) the U.S. Department of Energy Office of Science, Office of Advanced Scientific Computing Research, Applied Mathematics program under Award Number DE-SC-0011077 and (ii) the Defense Advanced Research Projects Agency's Enabling Quantification of Uncertainty in Physical Systems.

%\pagebreak
\appendix

\section{Proof of Theorem \ref{thm:postapprox}}
\label{sec:proofb}

First, we carefully work out the derivation for \eqref{eq:hprior1}, which uses the approximation $g$ from \eqref{eq:cexp}. The following quantities depend on $\vx$:
\begin{equation}
\post=\post(\vx),\quad
\pi=\pi(\vx),\quad
\prior=\prior(\vx),\quad
f=f(\vx),\quad
g=g(\mW_1^T\vx).
\end{equation}
In the next derivation, we omit the explicit dependence on $\vx$ to keep the notation clean.
\begin{align}
H^2(\post,\pi) 
&= \frac{1}{2}\int \left(
\left(\post\right)^{\frac{1}{2}}
- \left(\pi\right)^{\frac{1}{2}}
\right)^2\,d\vx\label{eq:helldef}\\
&= \frac{1}{2}\int \left(
\left(
\frac{\exp(-f)\,\prior}{\cpost}
\right)^{\frac{1}{2}}
- \left(
\frac{\exp(-g)\,\prior}{c_\pi}
\right)^{\frac{1}{2}}
\right)^2\,d\vx\label{eq:hp0}\\
&= \frac{1}{2}\int \left(
\left(
\frac{\exp(-f)}{\cpost}
\right)^{\frac{1}{2}}
- \left(
\frac{\exp(-g)}{c_\pi}
\right)^{\frac{1}{2}}
\right)^2\,\prior\,d\vx\label{eq:hp1}\\
&= \frac{1}{2(\cpost\,c_\pi)^{\frac{1}{2}}}
\left[
\int \left(
\left(\exp(-f)\right)^{\frac{1}{2}}
- \left(\exp(-g)\right)^{\frac{1}{2}}
\right)^2\,\prior\,d\vx
- \left(\cpost^{\frac{1}{2}} - \cpi^{\frac{1}{2}}\right)^2
\right]\label{eq:hp3}\\
&\leq \frac{1}{2(\cpost\,c_\pi)^{\frac{1}{2}}}
\int \left(
\left(\exp(-f)\right)^{\frac{1}{2}}
- \left(\exp(-g)\right)^{\frac{1}{2}}
\right)^2\,\prior\,d\vx\label{eq:hp4}\\
&= \frac{1}{2(\cpost\,c_\pi)^{\frac{1}{2}}}
\int \left(
\exp\left(\frac{-f}{2}\right)
- \exp\left(\frac{-g}{2}\right)
\right)^2\,\prior\,d\vx\\
&\leq \frac{1}{2(\cpost\,c_\pi)^{\frac{1}{2}}}
\int \left(\frac{1}{2}\,(f - g)\right)^2
\,\prior\,d\vx\label{eq:hp5}\\
&\leq \frac{C}{8(\cpost\, c_\pi)^{\frac{1}{2}}}
(\lambda_{n+1}+\cdots+\lambda_m)\label{eq:hp6}
\end{align}
Line \eqref{eq:helldef} is the definition of the squared Hellinger distance. Line \eqref{eq:hp0} plugs in the definitions of the posterior $\post$ and approximate posterior $\pi$ in terms of the misfit and its approximation. Line \eqref{eq:hp1} factors out the prior. Line \eqref{eq:hp3} is verified by inspection using the definitions of $\cpost$ \eqref{eq:bayes} and $\cpi$ \eqref{eq:cpi}. Line \eqref{eq:hp4} follows since the omitted squared term is positive. Line \eqref{eq:hp5} follows from the mean value theorem and the fact that $|\exp(-x)|\leq 1$ for $x\geq 0$. The last line follows from Theorem \ref{thm:cexperr}.

The constant $\cpi$ is bounded below using Jensen's inequality,
\begin{equation}
\begin{aligned}
\cpi &= \int \exp(-g)\,\prior\,d\vx\\
&\geq \exp\left(-\int g\,\prior\,d\vx\right)\\
&= \exp\left(-\int f\,\prior\,d\vx\right),
\end{aligned}
\end{equation}
where the last line follows from the construction of $g$ in \eqref{eq:cexp}. Then the constant from \eqref{eq:hp6} can be bounded
\begin{equation}
\label{eq:cpibnd}
\frac{1}{8}(\cpost\, \cpi)^{\frac{-1}{2}}
\;\leq\; \frac{1}{8}\left[
\cpost\,\exp\left(-\int f\,\prior\,d\vx\right)
\right]^{\frac{-1}{2}}.
\end{equation}
Recalling the definition of $\cpost$ from \eqref{eq:bayes} with the definition of the misfit,
\begin{equation}
\cpost \;=\; \int \exp(-f)\,\prior\,d\vx,
\end{equation}
which completes the proof of \eqref{eq:hprior1}. 

Equations \eqref{eq:helldef} through \eqref{eq:hp5} are identical if $\hpi$ replaces $\pi$, $\hg$ replaces $g$, and $\chpi$ replaces $\cpi$. The constant $\chpi$ is bounded as
\begin{equation}
\begin{aligned}
\chpi &= \int \exp(-\hg)\,\prior\,d\vx\\
&\geq \exp\left(-\int \hg\,\prior\,d\vx\right)\\
&= \exp\left(-\int f\,\prior\,d\vx\right),
\end{aligned}
\end{equation}
where the last line follows from the unbiasedness of the Monte Carlo approximation \eqref{eq:mccexp}. Then \eqref{eq:cpibnd} holds with $\chpi$ replacing $\cpi$. Finally, using Theorem \ref{thm:mccexperr} recovers \eqref{eq:hprior2}. Lines \eqref{eq:hprior3} and \eqref{eq:hprior4} use identical reasoning with the estimated eigenvectors $\hmW$ and Theorems \ref{thm:epscexperr} and \ref{thm:epsmccexperr}.

\bibliographystyle{siam}
\bibliography{mcmc-asm}

\begin{thebibliography}{10}

\bibitem{Apte2007}
{\sc A.~Apte, M.~Hairer, A.M. Stuart, and J.~Voss}, {\em Sampling the
  posterior: {An} approach to non-{Gaussian} data assimilation}, Physica D:
  Nonlinear Phenomena, 230 (2007), pp.~50 -- 64.
\newblock Data Assimilation.

\bibitem{kde2d}
{\sc Z.~Botev}, {\em kde2d}.
\newblock
  \url{http://www.mathworks.com/matlabcentral/fileexchange/17204-kde2d-data-n-min-xy-max-xy-},
  2007.
\newblock Retrieved Sept.~25, 2015.

\bibitem{brooks2011handbook}
{\sc S.~Brooks, A.~Gelman, G.~Jones, and X.-L. Meng}, {\em Handbook of Markov
  Chain Monte Carlo}, CRC press, Boca Raton, 2011.

\bibitem{Bui2012}
{\sc T.~Bui-Thanh, C.~Burstedde, O.~Ghattas, J.~Martin, G.~Stadler, and L.C.
  Wilcox}, {\em Extreme-scale {UQ} for {Bayesian} inverse problems governed by
  {PDEs}}, in Proceedings of the International Conference on High Performance
  Computing, Networking, Storage and Analysis, SC '12, Los Alamitos, CA, USA,
  2012, IEEE Computer Society Press, pp.~3:1--3:11.

\bibitem{Bui2014}
{\sc T.~Bui-Thanh and M.~Girolami}, {\em Solving large-scale {PDE}-constrained
  {Bayesian} inverse problems with {Riemann} manifold {Hamiltonian Monte
  Carlo}}, Inverse Problems, 30 (2014), p.~114014.

\bibitem{Calvetti2007}
{\sc D.~Calvetti and E.~Somersalo}, {\em Introduction to Bayesian Scientific
  Computing: Ten Lectures on Subjective Computing}, Springer, New York, 2007.

\bibitem{chen1982inequality}
{\sc L.H. Chen}, {\em An inequality for the multivariate normal distribution},
  Journal of Multivariate Analysis, 12 (1982), pp.~306--315.

\bibitem{asm2015}
{\sc P.G. Constantine}, {\em Active Subspaces: Emerging Ideas in Dimension
  Reduction for Parameter Studies}, SIAM, Philadelphia, 2015.

\bibitem{Constantine2014}
{\sc P.~Constantine, E.~Dow, and Q.~Wang}, {\em Active subspace methods in
  theory and practice: Applications to kriging surfaces}, SIAM Journal on
  Scientific Computing, 36 (2014), pp.~A1500--A1524.

\bibitem{constantine2014exploiting}
{\sc P.G. Constantine, M.~Emory, J.~Larsson, and G.~Iaccarino}, {\em Exploiting
  active subspaces to quantify uncertainty in the numerical simulation of the
  {HyShot II} scramjet}, Journal of Computational Physics, 302 (2015), pp.~1 --
  20.

\bibitem{constantine2015computing}
{\sc P.G. Constantine and D.F. Gleich}, {\em Computing active subspaces with
  {Monte Carlo}}, arXiv preprint arXiv:1408.0545v2,  (2015).

\bibitem{constantine2015discovering}
{\sc P.G. Constantine, B.~Zaharatos, and M.~Campanelli}, {\em Discovering an
  active subspace in a single-diode solar cell model}, Statistical Analysis and
  Data Mining: The ASA Data Science Journal,  (2015), pp.~n/a--n/a.

\bibitem{cui2014dimension}
{\sc T.~Cui, K.J.H. Law, and Y.M. Marzouk}, {\em Dimension-independent
  likelihood-informed {MCMC}}, Journal of Computational Physics, 304 (2016),
  pp.~109 -- 137.

\bibitem{Cui2014}
{\sc T.~Cui, J.~Martin, Y.M. Marzouk, A.~Solonen, and A.~Spantini}, {\em
  Likelihood-informed dimension reduction for nonlinear inverse problems},
  Inverse Problems, 30 (2014), p.~114015.

\bibitem{Flath2011}
{\sc H.P. Flath, L.C. Wilcox, V.~Ak\c{c}elik, J.~Hill, B.~van Bloemen~Waanders,
  and O.~Ghattas}, {\em Fast algorithms for {Bayesian} uncertainty
  quantification in large-scale linear inverse problems based on low-rank
  partial {Hessian} approximations}, SIAM Journal on Scientific Computing, 33
  (2011), pp.~407--432.

\bibitem{Flegal2008}
{\sc J.M. Flegal, M.~Haran, and G.L. Jones}, {\em {Markov} chain {Monte Carlo}:
  Can we trust the third significant figure?}, Statistical Science, 23 (2008),
  pp.~pp. 250--260.

\bibitem{Gibbs2002}
{\sc A.L. Gibbs and F.E. Su}, {\em On choosing and bounding probability
  metrics}, International Statistical Review, 70 (2002), pp.~419--435.

\bibitem{Girolami11}
{\sc M.~Girolami and B.~Calderhead}, {\em Riemann manifold {Langevin} and
  {Hamiltonian Monte Carlo} methods}, Journal of the Royal Statistical Society:
  Series B (Statistical Methodology), 73 (2011), pp.~123--214.

\bibitem{golub2013}
{\sc G.H. Golub and C.F.~Van Loan}, {\em Matrix Computations}, The Johns
  Hopkins University Press, Baltimore, 4~ed., 2013.

\bibitem{DRAM2006}
{\sc H.~Haario, M.~Laine, A.~Mira, and E.~Saksman}, {\em {DRAM}: {Efficient}
  adaptive {MCMC}}, Statistics and Computing, 16 (2006), pp.~339--354.

\bibitem{Jefferson2015}
{\sc J.L. Jefferson, J.M. Gilbert, P.G. Constantine, and R.M. Maxwell}, {\em
  Active subspaces for sensitivity analysis and dimension reduction of an
  integrated hydrologic model}, Computers \& Geosciences, 83 (2015), pp.~127 --
  138.

\bibitem{kaipio2005statistical}
{\sc J.~Kaipio and E.~Somersalo}, {\em Statistical and Computational Inverse
  Problems}, Springer, New York, 2005.

\bibitem{Lukaczyk2014}
{\sc T.W. Lukaczyk, P.~Constantine, F.~Palacios, and J.J. Alonso}, {\em Active
  subspaces for shape optimization}, American Institute of Aeronautics and
  Astronautics, 2014.

\bibitem{Martin2012}
{\sc J.~Martin, L.C. Wilcox, C.~Burstedde, and O.~Ghattas}, {\em A stochastic
  {Newton} {MCMC} method for large-scale statistical inverse problems with
  application to seismic inversion}, SIAM Journal on Scientific Computing, 34
  (2012), pp.~A1460--A1487.

\bibitem{Marzouk2009}
{\sc Y.M. Marzouk and H.N. Najm}, {\em Dimensionality reduction and polynomial
  chaos acceleration of {Bayesian} inference in inverse problems}, Journal of
  Computational Physics, 228 (2009), pp.~1862 -- 1902.

\bibitem{spantini2015optimal}
{\sc A.~Spantini, A.~Solonen, T.~Cui, J.~Martin, L.~Tenorio, and Y.M. Marzouk},
  {\em Optimal low-rank approximations of {Bayesian} linear inverse problems},
  arXiv preprint arXiv:1407.3463v2,  (2015).

\bibitem{Stuart2010}
{\sc A.M. Stuart}, {\em Inverse problems: A {Bayesian} perspective}, Acta
  Numerica, 19 (2010), pp.~451--559.

\bibitem{vrugt2009accelerating}
{\sc J.A. Vrugt, C.J.F. Ter~Braak, C.G.H. Diks, B.A. Robinson, J.M. Hyman, and
  D.~Higdon}, {\em Accelerating {Markov} chain {Monte Carlo} simulation by
  differential evolution with self-adaptive randomized subspace sampling},
  International Journal of Nonlinear Sciences and Numerical Simulation, 10
  (2009), pp.~273--290.

\bibitem{numpad}
{\sc Q.~Wang}, {\em Numpad: Numerical prototyping in {Python} assisted by
  automatic differentiation}.
\newblock \url{https://github.com/qiqi/numpad}, 2014.
\newblock Retrieved Nov.~29, 2014.

\end{thebibliography}

\end{document}